\documentclass[12pt,a4paper]{article}

\usepackage{amsbsy,amsfonts,amsmath,amssymb,amsthm}
\usepackage{array}
\usepackage{bm}
\usepackage[utf8]{inputenc}
\usepackage{color}
\usepackage{enumerate}
\usepackage{mathrsfs}
\usepackage{latexsym}
\usepackage[colorlinks=true]{hyperref}
\hypersetup{urlcolor=blue, citecolor=blue, linkcolor=blue}

\definecolor{wineRed}{rgb}{0.7,0,0.3}
\definecolor{grandBleu}{rgb}{0,0,0.8}
\definecolor{darkGreen}{rgb}{0,0.4,0}
\definecolor{lightGreen}{rgb}{0,0.7,0.3}
\definecolor{blueViolet}{rgb}{0.4,0,1.0}
\definecolor{redViolet}{rgb}{0.9,0,0.4}
\definecolor{bloodOrange}{rgb}{0.85,0.05,0}
\definecolor{mycolor}{rgb}{0.8,0,0.2}
\definecolor{}{rgb}{0.8,0,0.2}

\usepackage[textsize=small]{todonotes}
\usepackage{cite}

\newcommand{\bee}{\mathbf e}

\usepackage{comment}

\DeclareMathAlphabet{\mathpzc}{OT1}{pzc}{m}{it}

\usepackage[textsize=small]{todonotes}
\numberwithin{equation}{section}

\theoremstyle{plain}

\newtheorem{lem}{Lemma}

\theoremstyle{definition}
\newtheorem{defn}{Definition}
\newtheorem{thm}{Theorem}
\newtheorem{rem}{Remark}

\newtheorem{ex}{Example}

\def\N{\mathbb{N}}
\def\u{\mathbf{u}}
\def\w{\mathbf{w}}
\def\v{\mathbf{v}}
\def\e{\mathbf{e}}
\def\a{\mathbf{a}}
\def\b{\mathbf{b}}
\newcommand\boeta{\mbox{\boldmath{$\eta$}}}
\newcommand\blambda{\mbox{\boldmath{$\lambda$}}}
\newcommand\bomega{\mbox{\boldmath{$\omega$}}}
\newcommand\bmu{\mbox{\boldmath{$\mu$}}}
\def\R{\mathbb{R}}

\def\Ua{{U_{\varepsilon, \nu, \delta}}}
\def\ua{{{\bf u}_{\varepsilon, \nu, \delta}}}
\def\etaa{{\eta_{\varepsilon, \nu, \delta}}}
\def\ds{\displaystyle}
\def\ts{\textstyle}

\begin{document}
\thispagestyle{plain}
\begin{center}
    \textbf{\Large {
        Existence of {solutions} to a phase-field model \\ of 3D-grain boundary motion governed \\[0.75ex] by {a} regularized 1-harmonic type flow{}\footnotemark[1]}}
\end{center}
    \bigskip
\begin{center}
    {\sc Salvador Moll}\footnotemark[2]
    \\[1ex]
    Department d'An\`{a}lisi Matem\`{a}tica, Universitat de Val\`{e}ncia
    \\
    C/Dr. Moliner, 50, Burjassot, Spain
    \\
    {\ttfamily j.salvador.moll@uv.es}
\end{center}
\begin{center}
    \textsc{Ken Shirakawa}\footnotemark[3]
    \\[1ex]
    {Department of Mathematics, Faculty of Education, Chiba University,
    \\
    1-33, Yayoi-cho, Inage-ku, 263-8522, Chiba, Japan}
    \\
    {\ttfamily sirakawa@faculty.chiba-u.jp}
\end{center}
\begin{center}
    \textsc{and}
\end{center}
\begin{center}
    {\sc Hiroshi Watanabe}\footnotemark[4]
    \\[1ex]
    Division of Mathematical Sciences,
    \\
    Faculty of Science and Technology, Oita University
    \\
    700 Dannoharu, Oita, 870-1192, Japan
    \\
    {\ttfamily hwatanabe@oita-u.ac.jp}
\end{center}

\footnotetext[0]{

$\empty^*$\,AMS Subject Classification
	35K67,  
	35K87, 
	35Q99. 

$\empty^\dag$\,This author is supported by the Spanish MCIU and FEDER project PGC2018-094775-B-I00 and by Conselleria d'Innovaci\'o, Universitats, Ci\`encia i Societat Digital, project AICO/2021/223.

$\empty^\ddag$\,This author is supported by Grant-in-Aid No. 20K03672, JSPS.

$\empty^\S$\,This author is supported by Grant-in-Aid No. 20K03696, 21K03312, JSPS. 
}
\bigskip

\noindent
{\bf Abstract.}
In this paper we propose a quaternion formulation for the orientation variable in the three dimensional {Kobayashi--Warren} model for the dynamics of polycrystals. We obtain existence of solutions to the $L^2$-gradient descent flow of the constrained energy functional via several approximating problems. In particular, {we use} a Ginzburg-Landau type approach and some {extra} regularizations. Existence of solutions to the approximating problems is shown by the use of nonlinear semigroups. Coupled with good a-priori estimates, this leads to successive passages to the limit up to {finally} showing existence of solutions to the proposed model. Moreover, we also obtain a maximum principle for the orientation variable.
\newpage

\section{Introduction}
In \cite{kobayashi2005modeling} and \cite{pusztai2005phase,pusztai2006phase} two very similar models for the dynamics of polycrystals in the three dimensional space were introduced. The authors generalized the existing two dimensional model by {Kobayashi} et al \cite{MR1752970,MR1794359} to the case of 3D-crystals. In essence, the 2D-model is the $L^2$-gradient descent flow of the following energy functional:

 \begin{align*}
     [\eta, \theta] \in [H^1(\Omega)]^2 & \mapsto 
     \ds \frac{1}{2} \int_\Omega |\nabla \eta|^2 \, dx +\int_\Omega G(\eta) \, dx
     \\
     & +\int_\Omega \alpha(\eta) |\nabla \theta|\, dx + \frac{{\kappa^2}}{2}\int_\Omega |\nabla\theta|^2\, dx {\,\in (-\infty, \infty].}
\end{align*}
In the context, the unknowns $ \eta = \eta(t, x) $ and $ \theta = \theta(t, x) $ represent, respectively, the ``orientation order'' and the ``orientation angle'' in a polycrystal. $\eta=1$ corresponds to a completely ordered state while $\eta=0$ corresponds to the state where no meaningful value of mean orientation exists. In the original work \cite{MR1752970},  $G(s)=\frac{1}{2}(1-s)^2$ ensures that the ordered state $\eta=1$ is stable.

In order to generalize the model, one needs to consider orientations in 3D and misorientations, since the term $|\nabla \theta|$ represents the misorientation in a short scale. In 3D, orientations are elements of $SO(3)$, the special orthogonal group in $\R^3$. In \cite{kobayashi2005modeling}, the term $|\nabla\theta|$ is substituted by the corresponding Euclidean norm in $\R^9$; i.e. ${||\nabla P||_{\R^9}}:=\left(\sum_{i,j=1}^3|\nabla p_{i,j}|^2\right)^{\frac{1}{2}}$ for $P=[P]_{i,j}\in SO(3)$. Then, one has to compute the gradient descent flow for the constrained energy, thus ensuring that the solutions for the orientation variable still  belong to $SO(3)$.

In \cite{pusztai2005phase,pusztai2006phase}, instead, a quaternion representation is used for $SO(3)$. Then, by the use of the fact that quaternions can be identified as elements in the unit sphere in $\R^4$, {i.e. $\mathbb{S}^3$,} the authors replaced the term $|\nabla\theta|$ by  the Euclidean norm of the gradient of the quaternion: i.e. $|\nabla{\bf q}|:=\left(\sum_{i=0}^3|\nabla q^i|^2\right)^{\frac{1}{2}}$, for ${\bf q}=(q^0,q^1,q^2,q^3)\in \mathbb{S}^3$.

We take the point of view of \cite{pusztai2005phase,pusztai2006phase} and we consider the following energy functional; constrained to functions with values in the unit sphere of $\R^M$ {with $ 1 < M \in \N $,} i.e. $\mathbb{S}^{M-1}$:
\begin{align}\label{freeEgy}
    [\eta, {\bf u}] & \in L^2(\Omega) \times L^2(\Omega; \R^{M}) \mapsto \mathcal F(\eta,{\bf u})
    \nonumber
    \\
    & :=\left\{ \begin{array}{ll}
        \multicolumn{2}{l}{\ds \frac{1}{2} \int_\Omega |\nabla \eta|^2 \, dx +\int_\Omega G(\eta) \, dx +\int_\Omega \alpha(\eta) |\nabla {\bf u}|\, dx + \frac{{\kappa^2}}{2}\int_\Omega |\nabla {\bf u}|^2\, dx,}
        \\[2ex]
        & \mbox{if }~ \eta\in H^1(\Omega){~\mbox{and}~} {\bf u}\in H^1(\Omega;\mathbb S^{M-1}),
        \\[2ex]
        +\infty, & {\mbox{otherwise;}}
    \end{array} \right.
\end{align}
where
\begin{equation*}
    H^1(\Omega; \mathbb{S}^{M -1}) := \left\{ \begin{array}{l|l}
        {\bf \tilde{u}} \in H^1(\Omega; \R^M) & 
        |{\bf \tilde{u}}| = 1
        ~\mbox{a.e. in $ \Omega $}
    \end{array} \right\}.
\end{equation*}

Now, on any time interval $ (0, T) $ with a finite constant $ T > 0 $, the corresponding $L^2$-gradient descent flow of $\mathcal F$ can be formally written as the following system of PDEs:

\paragraph{\boldmath System (P):}{
    \begin{align*}
        & \begin{cases}
            \partial_t \eta -\mathit{\Delta} \eta +g(\eta) +\alpha'(\eta) |\nabla {\bf u}| = 0 \mbox{ in } (0, T)\times\Omega,
            \\[0.5ex]
            \nabla \eta \cdot {{\bf n}_\Gamma} = 0 \mbox{ on $ (0, T)\times\partial \Omega $,}
            \\[0.5ex]
            \eta(0, x) = \eta_0(x), ~ x \in \Omega;
        \end{cases}
        \\[1ex]
        & \begin{cases}
           \displaystyle \partial_t {\bf u} -\mathrm{div} \left( \alpha(\eta) \frac{\nabla {\bf u}}{|\nabla {\bf u}|}
 +{\kappa^2} \nabla {\bf u}\right) 
            =  \bigl( |\nabla {\bf u}|+\kappa^2|\nabla {\bf u}|^2 \bigr) {\bf u} ~~\mbox{in $ (0,\infty)\times\Omega $,}
            \\[0.5ex] \displaystyle
           \left( \alpha(\eta) \frac{\nabla {\bf u}}{|\nabla {\bf u}|}
 +{\kappa^2} \nabla {\bf u}\right) {{\bf n}_\Gamma} = 0 \mbox{ on $ (0,\infty)\times \partial\Omega $,}
            \\[0.5ex]
            {\bf u}(0, x) = {\bf u}_0(x), ~ x \in \Omega;
        \end{cases}
    \end{align*}}
We point out that the writing is purely formal since there are some undefined terms such as $\displaystyle \frac{\nabla {\bf u}}{|\nabla \bf u|}$. The precise meaning of a solution to the system is given in Section {\ref{sec:mainTh}}.

A natural way to study this type of restricted functionals in the sphere is to relax the constraint via a Ginzburg-Landau approximation (see \cite{MR1269538});
i.e. instead of considering ${\bf u}\in H^1(\Omega;\mathbb{S}^{M-1})$, one lets ${\bf u}\in H^1(\Omega; \R^{M})$ and adds the following term to $\mathcal F$:
$$\frac{1}{4\delta}\int_\Omega(|{\bf u}|^2-1)^2\,dx.$$ After obtaining well posedness {to} the gradient descent flow, 
then the strategy is to let $\delta\to 0^+$ and to show convergence of a subsequence to the corresponding solution to the {system (P).}

We follow exactly this strategy, but, due to the non-differentiability of the Euclidean norm at the origin, we also need to approximate it by a sufficient smooth term; i.e. we replace $$\int_\Omega \alpha(\eta) |\nabla {\bf u}|\,dx\mapstochar\longrightarrow\int_\Omega \alpha(\eta)\sqrt{\varepsilon^2+|\nabla{\bf u}|^2}\, dx,\quad {\rm for \ } \varepsilon>0.$$

Due to technical reasons for a possible future study, we need to perform an extra approximation. One might study the limit problem when $\kappa\to 0^+$ in \eqref{freeEgy}; i.e. the case in which ${\bf u}\in BV(\Omega;SO(3))$. In this case, the term in the energy corresponding to $\ds \int_\Omega \alpha(\eta)|\nabla{\bf u}|$ needs to be replaced by its relaxed functional. This relaxed functional has a jump part that strongly depends on the metric considered in $SO(3)$ (see \cite{giaquinta2008relaxation}). By the considerations stated in Appendix B, and in order to uniquely identify a rotation as an element in $\mathbb{S}^3$, we need to restrict the solutions in the quaternion representation to lie in the open upper hemisphere $\mathbb{S}^3_+:=\{{\bf p}=(p^1,p^2,p^3,p^4)\in \mathbb{S}^3: p^1>0 \}$.

For the sake of generality, we will consider the more general setting ${\bf u}\in\mathbb{S}^{M-1}$ instead of $\mathbb{S}^3$. In Appendix B, we give a maximum principle which ensures that, if the initial datum is in a certain compact subset of $\mathbb{S}^{M-1}_+$, then the solution also does. For the proof of this result, we need a technical restriction; namely, continuity of the solutions. Therefore, we need to perform an extra approximation to $\mathcal F$, by adding the following term to the energy functional:
{
$$
\frac{1}{N+1}\int_\Omega |\nu \nabla {\bf u}|^{N+1}\, dx.
\vspace{-2ex}
$$
}
\medskip

\noindent
We stress that this extra regularization is only a technical tool to prove the maximum principle in Theorem \ref{max_principle} and it is not needed for any of the rest of the results in the present manuscript.

\bigskip

{
The plan of the paper is the following one:

First of all, in Section  \ref{sec:preli}, {we prescribe some notations, and} recall some results about {multi-vectors} that are used in the paper. In Section \ref{sec:mainTh}, we {set up} our main assumptions, and we state the Main Theorem, as the principal result of this paper. In Section \ref{sec:approx_system} we consider the complete energy functional \eqref{freeEgy}, which we call the \emph{free energy,} we define our notion of solution to its $L^2$- gradient descent flow; i.e to the system named by (P)$_{\varepsilon,\nu,\delta}^\kappa$. Then, we prove with the help of an auxiliary convex energy functional, that the system (P)$_{\varepsilon,\nu,\delta}^\kappa$ admits a unique solution for sufficiently smooth initial data (see Theorem \ref{Th.SolvApKWC}). Moreover, stability with respect to the parameter $\nu$ is obtained.

Section \ref{sec:sphere} is devoted to the proof of Main Theorem, i.e. the proof of existence of solution to the $L^2$-gradient descent flow of $\mathcal F$.  First of all, an energy inequality together with the corresponding uniform estimates (in the parameters $\varepsilon$ and $\delta$) are obtained for solutions to (P)$_{\varepsilon,\nu,\delta}^\kappa$. They lead to convergence, first with $\delta\to 0^+$ and up to subsequences, to a solution to the system (P)$_{\varepsilon,\nu}$; i.e, to the gradient descent flow of the restricted energy functional:
$$
\mathcal F_{\varepsilon,\nu}(\eta,{\bf u}):=\left\{ \begin{array}{ll}
    \multicolumn{2}{l}{\ds \frac{1}{2} \int_\Omega |\nabla \eta|^2 \, dx +\int_\Omega G(\eta) \, dx  +\int_\Omega \alpha(\eta) \sqrt{\varepsilon^2+|\nabla {\bf u}|^2}\, dx}
    \\[2ex]
    & \ds \hspace{7ex} +\frac{1}{N+1}\int_\Omega |\nu \nabla {\bf u}|^{N+1}\, dx +\frac{{\kappa^2}}{2}\int_\Omega |\nabla {\bf u}|^2\, dx,
       \\[2ex]
       & \mbox{if }~ \eta\in H^1(\Omega){~\mbox{and}~} {\bf u}\in H^1(\Omega;\mathbb S^{M -1}),
       \\[2ex]
       +\infty, & \mbox{otherwise.}
   \end{array} \right.
$$
Moreover, solutions are shown to be continuous. Therefore, the maximum principle stated and proved in the Appendix applies and we can move further by letting $\nu\to 0^+$. Then, we obtain existence of solutions to the system (P)$_{\varepsilon}$; i.e. to the gradient descent flow of ${\mathcal F_{\varepsilon}}:=\mathcal F_{\varepsilon,0}$. The final step is to let $\varepsilon\to 0^+$, thus obtaining existence of solutions to the gradient descent flow of $\mathcal F$; i.e. to (P). We point out that, since the successive convergences for the orientation variable (with $\nu\to 0^+$ and $\varepsilon\to 0^+$) also hold a.e. in space time, then, the final solutions also satisfy the maximum principle.

Finally, we added two {Appendices}. In the first one, we recall the concept of Mosco convergence and some results related to it that we use in the paper. The second one is devoted to the discussion about the relationship between rotations in $SO(3)$ and their representation as quaternions. It is there where we prove our maximum principle.
}

\section{Preliminaries and assumptions}\label{sec:preli}

\subsection{Abstract notations}
For an abstract Banach space $ X $, we denote by $ \|\cdot\|_{X} $ the norm of $ X $, and by $ \langle \cdot, \cdot \rangle_X $ the duality pairing between $ X $ and its dual $ X' $. In particular, when $ X $ is a Hilbert space, we denote by $ (\cdot,\cdot)_{X} $ the inner product of $ X $. Moreover, when there is no possibility of confusion, we uniformly denote by $ |\cdot| $ the norm of Euclidean spaces, and for any dimension $ d \in \N $, we write the inner product (scalar product) of $ \R^d $, as follows:
    \begin{align*}
        & y \cdot \tilde{y} = \sum_{i = 1}^d y_i \tilde{y}_i, 
        \mbox{ for all } y = [y_1, \dots, y_d], ~ \tilde{y} = [\tilde{y}_1, \dots, \tilde{y}_d] \in \R^d.
    \end{align*}

For any subset $ A $ of a Banach space $ X $, let $ \chi_A : X \longrightarrow \{0, 1\} $ be the characteristic function of $ A $, i.e.:
    \begin{equation*}
        \chi_A: w \in X \mapsto \chi_A(w) := \begin{cases}
            1, \mbox{ if $ w \in A $,}
            \\[0.5ex]
            0, \mbox{ otherwise.}
        \end{cases}
    \end{equation*}

For Banach spaces $ X_1, \dots, X_d $, with $ 1 < d \in \mathbb{N} $, let $ X_1 \times \dots \times X_d $ be the product Banach space endowed with the norm $ \|\cdot\|_{X_1 \times \cdots \times X_d} := \|\cdot\|_{X_1} + \cdots +\|\cdot\|_{X_d} $. However, when all $ X_1, \dots, X_d $ are Hilbert spaces, $ X_1 \times \dots \times X_d $ denotes the product Hilbert space endowed with the inner product $ (\cdot, \cdot)_{X_1 \times \cdots \times X_d} := (\cdot, \cdot)_{X_1} + \cdots +(\cdot, \cdot)_{X_d} $ and the norm $ \|\cdot\|_{X_1 \times \cdots \times X_d} := \bigl( \|\cdot\|_{X_1}^2 + \cdots +\|\cdot\|_{X_d}^2 \bigr)^{\frac{1}{2}} $. In particular, when all $ X_1, \dots,  X_d $ coincide with a Banach space $ Y $, we write:
\begin{equation*}
    [Y]^d := \overbrace{Y \times \cdots \times Y}^{\mbox{$d$ times}}.
\end{equation*}

For a proper, lower semi-continuous (l.s.c.), and convex function $ \Psi : X \to (-\infty, \infty] $ on a Hilbert space $ X $, we denote by $ D(\Psi) $ the effective domain of $ \Psi $. Also, we denote by $\partial \Psi$ the subdifferential of $\Psi$. The subdifferential $ \partial \Psi $ corresponds to a weak differential of convex function $ \Psi $, and it is known as a maximal monotone graph in the product space $ X \times X $. The set $ D(\partial \Psi) := \bigl\{ z \in X \ |\ \partial \Psi(z) \neq \emptyset \bigr\} $ is called the domain of $ \partial \Psi $. We often use the notation ``$ [w_{0}, w_{0}^{*}] \in \partial \Psi $ in $ X \times X $\,'', to mean that ``$ w_{0}^{*} \in \partial \Psi(w_{0})$ in $ X $ for $ w_{0} \in D(\partial\Psi) $'', by identifying the operator $ \partial \Psi $ with its graph in $ X \times X $.
\medskip

Next, for Hilbert spaces $X_1, \cdots, X_d$, with $1<d \in \mathbb{N}$, let us consider a proper, l.s.c., and convex function on the product space $X_1 \times \dots \times X_d$:
\begin{equation*}
\widehat{\Psi}: w = [w_1,\cdots,w_d] \in X_1 \times\cdots\times X_d \mapsto \widehat{\Psi}(w)=\widehat{\Psi}(w_1,\cdots,w_d) \in (-\infty,\infty]. 
\end{equation*}
Besides, for any $i \in \{1, \dots, d\}$, we denote by $\partial_{w_i} \widehat{\Psi}:X_1 \times \cdots \times X_d \to X_i$ a set-valued operator, which maps any $w=[w_1,\dots,w_i,\dots,w_d] \in X_1 \times \dots \times X_i \times \dots \times X_d$ to a subset $ \partial_{w_i} \widehat{\Psi}(w) \subset  X_i $, prescribed as follows:
\begin{equation*}
\begin{array}{rl}
\partial_{w_i}\widehat{\Psi}(w)&=\partial_{w_i}\widehat{\Psi}(w_1,\cdots,w_i,\cdots,w_d)
\\[2ex]
&:= \left\{\begin{array}{l|l}\tilde{w}^* \in X_i & \begin{array}{ll}\multicolumn{2}{l}{(\tilde{w}^*,\tilde{w}-w_i)_{X_i} \le \widehat{\Psi}(w_1,\cdots,\tilde{w},\cdots, w_d)}
\\[0.25ex] 
& \quad -\widehat{\Psi}(w_1,\cdots,w_i,\cdots,w_d), \mbox{ for any $\tilde{w} \in X_i$}\end{array}
\end{array}\right\}.
\end{array}
\end{equation*}
As is easily checked, 
\begin{equation}\label{prodSubDif}
    \partial \widehat{\Psi} \subset \bigl[ \partial_{w_1} \widehat{\Psi} \times \cdots \times \partial_{w_d} \widehat{\Psi} \bigr] \mbox{ in $ [X_1 \times \cdots \times X_d]^2 $,}
\end{equation}
where $ \bigl[ \partial_{w_1} \widehat{\Psi} \times \cdots \times \partial_{w_d} \widehat{\Psi} \bigr] : X_1 \times \cdots \times X_d \longrightarrow 2^{X_1 \times \cdots \times X_d} $ is a set-valued operator, defined as:
\begin{equation*}
    \begin{array}{c}
        \ds \bigl[ \partial_{w_1} \widehat{\Psi} \times \cdots \times \partial_{w_d} \widehat{\Psi} \bigr](w) := \partial_{w_1} \widehat{\Psi}(w) \times \cdots \times \partial_{w_d} \widehat{\Psi}(w) \mbox{ in $ X_1 \times \cdots \times X_d $},
        \\[1.5ex]
        \mbox{for any } w = [w_1, \dots, w_d] \in  D\bigl( \bigl[ \partial_{w_1} \widehat{\Psi} \times \cdots \times \partial_{w_d} \widehat{\Psi} \bigr] \bigr) := D(\partial_{w_1} \widehat{\Psi}) \cap \cdots \cap D(\partial_{w_d} \widehat{\Psi}).
    \end{array}
\end{equation*}
{It} should be noted that the converse inclusion of \eqref{prodSubDif} is not true, in general. 
\subsection{Multi-vectors}\label{ss-multi}

Here we recall some definitions and basic properties about multi-vectors that we need in our analysis. We refer to e.g. \cite[Chapter 1]{Federer} and \cite[Chapter 1]{Darling} for details.

Let $ m \in \N $. The spaces $\Lambda_0(\R^m)$ and $\Lambda_1(\R^m)$ are defined as:
\begin{equation*}
    \Lambda_0(\R^{m}) := \R ~\mbox{and}~ \Lambda_1(\R^m) := \R^m, ~\mbox{{respectively.}}
\end{equation*}
For any integer $2 \le k\leq m$, the $k$-th exterior power of $\R^m$, denoted by $\Lambda_k(\R^m)$, is defined as a set spanned by elements of the form:
\begin{equation*}
\u_{1} \wedge \cdots \wedge \u_{k}, \ \ \u_i \in \R^m, \ i=1, \dots , k,
\end{equation*}
called ``generators'', which are subject to {the} following rules:
\begin{eqnarray}\nonumber
& (a\v+b\w)\wedge \u_{2} \wedge \cdots \wedge \u_{k} = a(\v\wedge \u_{2} \wedge \cdots \wedge \u_{k})+ b(\w\wedge \u_{2} \wedge \cdots \wedge \u_{k});
    & \\[2ex] \nonumber &
\mbox{$\u_{1} \wedge \cdots \wedge \u_{k}$ changes sign if two entries are transposed;}
    & \\[2ex] \label{k-basis}
    & \left. \parbox{12cm}{
            for any basis $ \bigl\{ {\bf e}_1, \dots ,{\bf e}_m \bigr\} $ of $ \R^m $, the set
            \vspace{-1.5ex}
            $$
            \left\{ \begin{array}{l|l}
                {\bf e}_\alpha := {\bf e}_{\alpha_1} \wedge \dots \wedge {\bf e}_{\alpha_k} &
                \alpha = [\alpha_1, \dots, \alpha_k] \in I(k, m)
            \end{array} \right\}
            \vspace{-1.5ex}
            $$
            forms the basis of $ \Lambda_k(\R^m) $, where
            $$
            I(k, m) := \left\{ \begin{array}{l|l}
                \alpha = [\alpha_1, \dots, \alpha_k] \in \mathbb{Z}^k & 1 \leq \alpha_1 < \dots < \alpha_k \leq m
            \end{array} \right\}.
            $$
        } \right\}
& \end{eqnarray}
The elements of $\Lambda_k(\R^m)$ are called multi-vectors (or $k$-vectors), and $\Lambda_k(\R^m)$  is a vector space of dimension {\scriptsize {$\left( \hspace{-0.75ex} \begin{array}{c}  m \\ k \end{array} \hspace{-0.75ex} \right)$}}. 
    Given $k,\ell \in \{ 0, \dots, m \}$ with $k +\ell \le m$, there exists a unique bilinear map $(\blambda, \bmu) \to \blambda \wedge \bmu$ from $\Lambda_k(\R^m) \times \Lambda_\ell(\R^m)$ to $\Lambda_{k +\ell}(\R^m)$, whose effect on generators is
$$
(\u_{1} \wedge \u_{2} \wedge \cdots \wedge \u_{k})  \wedge ( \v_{1} \wedge \v_{2} \wedge \cdots \wedge \v_{\ell}) = \u_{1} \wedge \u_{2} \wedge \cdots \wedge \u_{k}  \wedge \v_{1} \wedge \v_{2} \wedge \cdots \wedge \v_{\ell}.
$$
{Such a map} satisfies
\begin{equation}\label{wedge-comm}
\blambda \wedge \bmu = (-1)^{-k\ell} (\bmu \wedge \blambda) \quad \hbox{for }\ \blambda \in \Lambda_k(\R^m), \ \bmu \in \Lambda_\ell(\R^m).
\end{equation}

The {\it Hodge-star operator}:
\begin{equation*}
    *({}\cdot{}) : \Lambda_k(\R^m) \ni \bm{\lambda} \longrightarrow *\bm{\lambda}\in \lambda_{m -k}(\R^m),
\end{equation*}
is an isomorphism between $\Lambda_k(\R^m)$ and $\Lambda_{m-k}(\R^m)$, which is defined on the basis as:
\begin{equation}
\label{def-hodge}
    \begin{array}{c}
        *(\e_{\alpha_1} \wedge \cdots \wedge \e_{\alpha_k}):= \e_{\alpha_{k+1}} \wedge \cdots \wedge \e_{\alpha_{m}},
        \\[1ex]
        \mbox{for all permutations $ \{ \alpha_1, \dots, \alpha_m \} $ of $ \{1, \dots, m\} $, having positive signature.}
    \end{array}
\end{equation}
In particular, in what follows we will systematically identify $\Lambda_{m-1}(\R^m)$ with $\R^m$ and $\Lambda_m(\R^m)$ with $\R$. We will use the following well known formulas:
\begin{equation}
\label{dobleprod}
*(* \blambda) = (-1)^{k(m-k)} \blambda \quad \mbox{for all } \ \blambda \in \Lambda_k(\R^m),
\end{equation}
(see e.g. \cite[(1.64)]{Darling}) and
\begin{equation}
\label{tripleprod}
    \begin{array}{c}
        \a\wedge *(\b\wedge \bold{c})=(\a\cdot \bold{c}) *\b - (\a\cdot \b) *\bold{c} ~~\mbox{in $ \R^m $ $ (= \Lambda_{m -1}(\R^m)) $,}
        \\[1ex]
        \mbox{for all } \  \a, \b, \bold{c} \in \R^m ~(= \Lambda_1(\R^m))
    \end{array}
\end{equation}
(see e.g.  \cite[Table 1.2]{Darling}).
It follows from (\ref{wedge-comm}), (\ref{tripleprod}), and \eqref{dobleprod} that
\begin{equation}
\label{muuibn}
    \begin{array}{c}
        |\b|^2 \a = (\a \cdot \b) \b - *(*(\a \wedge \b)\wedge \b) ~~\mbox{in $ \R^m $ $ (= \Lambda_1(\R^m)) $,}
        \\[1ex]
        \quad \mbox{for all }\  \a, \b \in \R^m ~(= \Lambda_1(\R^m)).
    \end{array}
\end{equation}

Introducing the inner product on generators $\blambda=\blambda_1\wedge\ldots\wedge \blambda_k$, $\bmu=\bmu_1\wedge\ldots\wedge \bmu_k$, and then extending by linearity to $\Lambda_k(\R^m)$
\begin{equation}\label{inner}
  \langle \blambda,\bmu\rangle_k :={\rm det}\left(\langle \blambda_i,\bmu_j\rangle_{i,j=1}^k\right),
\end{equation}we easily see that \begin{equation}\label{hodge_inner}\langle \blambda,\bmu\rangle_k=\blambda\wedge *\bmu.\end{equation}
Moreover,
\begin{equation}\label{use-as-def}
    |\blambda|_k := \langle \blambda,\blambda\rangle_k^{\frac{1}{2}}=\left(\sum_{\alpha\in I(k, m)}|\lambda_\alpha|^2\right)^\frac{1}{2}\,,\quad\mbox{where}\quad \blambda=\sum_{\alpha\in I(k, m)}\lambda_\alpha \e_\alpha,
\end{equation}
and using \eqref{def-hodge}, it is immediate to see that
\begin{equation}
   \label{normhodge}|*\blambda|_{m-k}=|\blambda|_k\quad {\rm for\ any \ } \blambda\in \Lambda_k(\R^m).
 \end{equation}
Finally, we recall that, given $\blambda\in \Lambda_k(\R^m)$ and $\boeta\in \Lambda_{\ell}(\R^m)$ such that one of them is a generator, then
\begin{equation}
  \label{CSformultivector}|\blambda\wedge\boeta|_{k +\ell}\leq |\blambda|_k|\boeta|_\ell
\end{equation}
(see \cite[pag. 32]{Federer}).

\subsection{Vector valued functions}\label{ss-banach}

Let 
    $X$ be a Banach space  with dual $X'$ and let $ U \subset \R^d$ be a bounded open set endowed with the Lebesgue measure $\mathcal L^d$. 
    A function $u: U \rightarrow X$ is  called {\it simple} if there exist $x_1, \ldots , x_n \in X$ and $U_1, \ldots, U_n$ $\mathcal L^m$-measurable subsets of $U$ such that $u = \sum_{i=1}^n x_i \chi_{U_i}$. The function $u$ is called {\it strongly measurable} if there exists a sequence of simple functions $\{ u_n \}$ such that $\Vert u_n(x) - u(x) \Vert_X \to 0$ as $n\to +\infty$ for almost all $x \in U$. If $1 \leq p < \infty$, then $L^p(U; X)$ stands for the space of (equivalence classes of) strongly measurable functions $u: U \rightarrow X$ with
 $$
    \Vert u \Vert_{L^p(U; X)} := \left( \int_U \Vert u(x) \Vert_X^p \, dx \right)^{\frac{1}{p}}< \infty.$$
 Endowed with this norm, $L^p(U; X)$ is
 is a Banach space. For $p= \infty$, the symbol $L^{\infty}(U; X)$ stands for the space of (equivalence classes of) strongly measurable functions $u: U \rightarrow X$ such that
    $$\Vert u \Vert_{L^\infty(U; X)} := \hbox{ esssup} \bigl\{ \Vert u(x) \Vert_X : \ x \in U \bigr\} < \infty.$$

If $U = (0, T)$ with $ 0 < T \leq \infty $, we write  $L^p(0, T; X) = L^p((0,T); X)$.  For $1 \leq p < \infty$, $L^{p'}(0,T; X')$ ($\frac{1}{p} + \frac{1}{p'} =1$) is isometric to a subspace of
%
%
$(L^p(0,T;X))'$,
with equality if and only if $X'$ has the Radon-Nikod\'ym property (see for instance \cite{MR0453964}).

We consider the vector space $\mathcal D(U;X):= C_0^\infty(U;X)$, endowed with the topology for which a sequence $\varphi_n \to 0$ as $n\to +\infty$ if there exists $K \subset U$ compact such that ${\rm supp} (\varphi_n) \subset K$ for any $n\in\N$ and $D^{\alpha} \varphi_n \to 0$ uniformly on $K$ as $n\to +\infty$ for all multi-index $\alpha$. We denote by $\mathcal D'(U ; X)$ the space of distributions on $U$ with values in $X${;} that is, the set of all linear continuous maps $T:\mathcal D(U;X)\to\R$.
As is well known, $L^p(U; X) \subset  \mathcal D'(U; X)$ through the standard continuous injection. Given $T \in \mathcal D'(U; X)$, the distributional derivative of $T$ is defined by
\begin{equation}\label{def-distributional}
    \langle D_i T, \varphi \rangle:= - \langle T, \partial_i \varphi \rangle\,\quad \mbox{for any $\varphi\in \mathcal D(U;X)$ and any $i\in\{1,\dots,d\}$}.
\end{equation}

\medskip

    \noindent{\it General notations for matrices.} Let $ d, m \in \N ${.} If $A = [ a^{\ell}_k \bigr] = \bigl[ a^{\ell}_k \bigr]_{{1 \leq \ell \leq m}\atop{1 \leq k \leq d}} \in \R^{md}$ is an $m \times d$ matrix, we write
\begin{align*}
    & \begin{cases}
        {\bf a}^\ell = [a^\ell_1, \ldots a_d^\ell] \in \R^{d} ~~\mbox{for $\ell = 1, \dots, m$,}
        \\[1ex]
        \a_k = {}^\mathrm{t}[a^1_k, \ldots , a_k^m \, ] \in \R^{m} ~~\mbox{for $k = 1, \dots, d$.}
    \end{cases}
\end{align*}
If  $B = [b^{\ell}_k] = \bigl[ b^{\ell}_k \bigr]_{{1 \leq \ell \leq m}\atop{1 \leq k \leq d}} \in \R^{md} $ is also an  $m \times d$ matrix, we let
$$
A: B = \sum_{\ell=1}^m \sum_{k=1}^d a^{\ell}_i b^{\ell}_i\quad\mbox{ and }\quad
\vert  A \vert = ( A :  A)^{\frac{1}{2}} = \left( \sum_{\ell=1}^m \sum_{k=1}^d
(a^{\ell}_i)^2 \right)^{\frac{1}{2}}.
$$

    Given $A = \big[ \mathbf a_1, \dots,\mathbf a_m \bigr] \in \R^{d\times m}$ with $ {\bf a}_i \in \R^{d} $, $ i = 1, \dots, m $, and $\b\in\R^d$, we let
\begin{eqnarray*}
A\wedge \b &:=& \left(\mathbf a_1 \wedge \b, \ldots , \mathbf a_m \wedge \b\right),
\\
*(A\wedge \b)&:=& \left(*(\mathbf a_1 \wedge \b), \ldots , *(\mathbf a_m \wedge \b)\right).
\end{eqnarray*}

\subsection{Multi-vector fields}

Let $ d, m \in \N $ be constants of dimensions. Let $U\subset\R^d$ be a bounded open set.
A multi-vector distribution in $U$ is a linear continuous map $\blambda\in \mathcal D'(U; \Lambda_k(\R^m))$ (see \S \ref{ss-banach}). It may be expressed in terms of the basis \eqref{k-basis} as
$$
\blambda=\sum_{\alpha\in I(k, m)}\lambda_\alpha \bee_\alpha, \quad \mbox{with  $\lambda_\alpha\in \mathcal D'(U;\R^m)$ for any $\alpha\in I(k,m)$}.
$$
Then, according to \eqref{def-distributional},
\begin{equation}\label{def-D-multi}
D_i \blambda =\sum_{\alpha\in I(k, m)}D_i\lambda_\alpha \bee_\alpha \quad\mbox{for any $i\in\{1,\ldots,d\}$}.
\end{equation}
 From \eqref{def-D-multi}, the following two identities are easily seen to hold for $k,\ell\in \N$ and $i\in\{1,\dots,d\}$:
\begin{eqnarray}
   \label{distrwedge} \ D_i(\blambda\wedge \boeta)=D_i\blambda\wedge\boeta+\blambda\wedge D_i\boeta
\end{eqnarray}
for  any $\blambda\in L^2(U;\Lambda_k(\R^m))$ such that $D_i\blambda\in L^2(U;\Lambda_k(\R^m))$ and any  $\boeta\in L^2(U;\Lambda_\ell(\R^m))$ such that $D_i\boeta\in L^2(U;\Lambda_\ell(\R^m))$;
\begin{eqnarray}
   \label{hodgederiv} *(D_i\blambda)=D_i(*\blambda)\, && {\rm for \ any \ }\blambda\in \mathcal D'(U;\Lambda_k(\R^m)).
 \end{eqnarray}
For any $k\in \N$, $ [\Lambda_k(\R^m)]^m$ is a Banach space. We use the norm
$$
\|\mathcal A\|_{[\Lambda_k(\R^m)]^m}:=\left(\sum_{i=1}^m |\mathcal{A}_i|_k^2\right)^{\frac{1}{2}}\,,\quad {\rm for \ } \mathcal A=(\mathcal A_1,\ldots,\mathcal A_m)
$$
with $|\cdot|_k$ given by \eqref{use-as-def}.

\section{Main Theorem}\label{sec:mainTh}


We start with setting up the assumptions in the principal part of this paper. The assumptions also fix the notations in the system (P), and in its approximating problems.
\begin{description}
    \item[(A0)]$ 1 < N \in \mathbb{N} $, $ 1 < M \in \mathbb{N}$, $ \kappa > 0 $ and {$ 0 < T < \infty $.} 
    \item[(A1)] $ \Omega \subset \mathbb{R}^N $ is a bounded domain with  Lipschitz boundary $ \Gamma := \partial \Omega $ and  unit outer normal $ {{\bf n}_\Gamma} $. Besides, let:
        \begin{equation*}
            \Omega_T := (0, T) \times \Omega, ~\mbox{and}~ \Gamma_T := (0, T) \times \Gamma,
        \end{equation*}
        \begin{equation*}
            H := L^2(\Omega), ~ {\mathbb{X}} := L^2(\Omega; \R^M), ~ \mathfrak{X} := H \times {\mathbb{X}}.
        \end{equation*}
        \begin{equation*}
            V := H^1(\Omega), ~ {\mathbb{W}} := H^1(\Omega; \R^M), ~ \mathfrak{W} := V \times {\mathbb{W}}.
        \end{equation*}
    \item[(A2)]$ g \in C^1(\R) $ is a fixed Lipschitz function such that
        $ g(0) \leq 0 $ and $ g(1) \geq 0 $, with a potential $ 0 \leq G \in C^2(\R) $, i.e. $ G'(s) = \frac{d}{d s} G(s) = g(s) $ {on $ \R $.
        }
    \item[(A3)]$ 0 < \alpha_0 \in W_\mathrm{loc}^{1, \infty}(\R) $ and $ 0 < \alpha \in C^2(\R) $ are fixed functions, such that:
        \begin{itemize}
            \item $ \alpha'(0) = 0 $, ~ $ \alpha'' \geqq 0 $ on $ \R $, {and $ \alpha $ and $ \alpha \alpha' $ are Lipschitz continuous on $ \R $.}
            \item $ \alpha^* := \inf \alpha_0(\R) \cup \alpha(\R) > 0 $.
        \end{itemize}
    \item[(A4)]For any $ \varepsilon \geq 0 $,  $ f_\varepsilon : \R^{MN} \longrightarrow [0, \infty) $ is a continuous convex function, defined as:
        \begin{equation*}
            f_\varepsilon : W = \bigl[ w_k^\ell \bigr]_{{1 \le \ell \le M}\atop{1 \le k \le N}} \in \R^{MN} \mapsto f_\varepsilon(W) := \sqrt{\varepsilon^2 +|W|^2} \in \R.
        \end{equation*}
    \item[(A5)]For any $ \delta > 0 $, $ \Pi_\delta \in C^2(\R^M) $ is the following 
        Ginzburg--Landau type potential function:
        \begin{equation*}
            \Pi_\delta : {\bf w} \in \R^M \mapsto \Pi_\delta({\bf w}) := \frac{1}{4 \delta} (|{\bf w}|^2 -1)^2 {\in \R.} 
        \end{equation*}
        We let $ \varpi_\delta \in C^1(\R^M; \R^{M}) $ be the gradient of $ \Pi_\delta $, i.e.:
        \begin{equation*}
            \varpi_\delta : {\bf w} \in \R^M \mapsto \varpi_\delta({\bf w}) := \nabla \Pi_\delta({\bf w}) = \frac{1}{\delta} \bigl(|{\bf w}|^2 -1 \bigr) {\bf w} \in \R^M. 
        \end{equation*}
    \item[(A6)] $ U_0 := [\eta_0, {\bf u}_0] \in \mathfrak{W} $ is a fixed pair of initial data.
\end{description}
\begin{rem}\label{Rem.(A4)}
    From (A4), it immediately follows that $ f_\varepsilon : \R^{MN} \longrightarrow [0, \infty) $, $ \varepsilon \geq 0 $, are non-expansive over $ \R^{MN} $. Also,  if $ \varepsilon > 0 $, then $ f_\varepsilon \in C^\infty(\R^{MN}) $, and if $ \varepsilon = 0 $, then the corresponding function $ f_0 $ coincides with the (Euclidean) norm $ \|\cdot\|_{\R^{MN}} $ on $ \R^{MN} $. Additionally, concerning the subdifferentials $ \partial f_\varepsilon $, $ \varepsilon \geq 0 $, we can observe that:
    \begin{align*}
        \partial f_\varepsilon(W ) = & \left\{ \begin{array}{ll}
            \ds \left\{ \frac{W}{\sqrt{\varepsilon^2 +|W|^2}} \right\} ~ \bigl( = [\nabla f_\varepsilon](W) \bigr), & \mbox{if $ \varepsilon > 0 $,}
            \\[3ex]
            \ds \mathrm{Sgn}^{M, N}(W), & \mbox{if $ \varepsilon = 0 $,}
        \end{array} \right.
        \\[1ex]
        &
        \mbox{for all $ W = \bigl[ w_k^\ell \bigr]_{{1 \le \ell \le M}\atop{1 \le k \le N}} \in \R^{MN} $,}
    \end{align*}
    where $ \mathrm{Sgn}^{M, N} : \R^{MN} \longrightarrow 2^{\R^{MN}} $ is the \emph{sign function} on $ \R^{MN} $, i.e.
    \begin{align*}
        \mathrm{Sgn}^{M, N}(W& ) := \left\{ \begin{array}{ll}
            \ds \left\{ \frac{W}{|W|} \right\}, & \mbox{if $ W \ne 0 $,}
            \\[3ex]
            \ds \left\{ \begin{array}{l|l}
                \tilde{W} & \tilde{W} : \tilde{W} \leq 1
            \end{array} \right\}, & \mbox{if $ W = 0 $,}
        \end{array} \right.
        \\[1ex]
        &
        \mbox{for all $ W = \bigl[ w_k^\ell \bigr]_{{1 \le \ell \le M}\atop{1 \le k \le N}} \in \R^{MN} $.}
    \end{align*}
\end{rem}

Next, we define the notion of solution to our system (P).
\begin{defn}\label{DefOfSol(S)}
    A pair of functions $ U := [\eta, {\bf u}] \in L^2(0, T; \mathfrak{X}) $ is called a solution to the system (P), iff.:
    \begin{equation}\label{S00}
        \begin{cases}
            U = [\eta, {\bf u}] \in W^{1, 2}(0, T; \mathfrak{X}) \cap L^\infty(0, T; \mathfrak{W}),
            \\[1ex]
        0 \leq \eta \leq 1 ~\mbox{and}~ {\bf u} \in \mathbb{S}^{M -1}, ~\mbox{a.e. in $ {\Omega_T} $,}
        \end{cases}
    \end{equation}
    \begin{align}\label{S01}
        \bigl( \partial_t \eta(t) & +g(\eta(t)) +\alpha'(\eta(t)) |\nabla {\bf u}(t)|, \varphi \bigr)_{H} +\bigl( \nabla \eta(t), \nabla \varphi \bigr)_{H} = 0,
        \nonumber
        \\
        & \mbox{for any $ \varphi \in V $, a.e. $ t \in (0, T) $, ~subject to $ \eta(0) = \eta_0 $ in $ H $;}
    \end{align}
    and there exist functions $ \mathcal{B}^* \in L^\infty(\Omega_T; \R^{MN}) $ and $ \mu^* \in L^1(0, T; L^1(\Omega)) $, such that:
    \begin{subequations}\label{S02}
    \begin{equation}\label{S02a}
        \begin{cases}
            \mathcal{B}^* \in \mathrm{Sgn}^{M, N}(\nabla {\bf u}) ~\mbox{in $ \R^{MN} $,}              \\[0.5ex]
            \mu^* := (\alpha(\eta)\mathcal{B}^* + {\kappa^2}\nabla {\bf u}) : \nabla {\bf u},
        \end{cases}
        \mbox{a.e. in $ \Omega_T $,}
    \end{equation}
    and
\begin{align}
    & \int_\Omega \partial_t {\bf u}(t) \cdot \bm{\psi}(t) \, dx  +\int_\Omega \alpha(\eta(t)) \mathcal{B}^*(t) : \nabla \bm{\psi}(t) \, dx = \int_\Omega \mu^*(t) {\bf u}(t) \cdot \bm{\psi}(t) \, dx
        \nonumber
        \\
        & \qquad \mbox{for any $ \bm{\psi} \in C^1(\overline{\Omega}; \R^{M}) $, a.e. $ {t \in (0, T)} $, ~subject to $ {\bf u}(0) = {\bf u}_0 $ in $ {\mathbb{X}} $.}
        \label{S02b}
    \end{align}
    \end{subequations}
\end{defn}

Now, the Main Theorem of this paper is stated as follows.
\bigskip

\noindent
\textbf{Main Theorem.}~~
    Let ${ U}_{0}=[\eta_{0}, {\bf u}_{0}] \in {\mathfrak{W}}
$ with 
    ${\bf u}_{0} \in {\mathbb{S}^{M-1}}$ in $\Omega$. Then, the system (P) admits at least one solution ${ U} = [\eta, {\bf u}] \in L^2(0, T;{\mathfrak{X}}) $, such that:
\begin{equation*}
\mathcal{F}({U}(s)) + \int_{0}^{s} \| \partial_{t}{U}(t) \|^{2}_{\mathfrak{X} 
    } dt \le \mathcal{F}({ U}_{0})\ \ \mbox{ for all } s \in (0, T).
\end{equation*}
    Also, concerning the function $ \mathcal{B}^* \in L^\infty(\Omega_T; \R^{MN}) $  in \eqref{S02}, it holds that:
    \begin{equation*}
        \begin{cases}
            \mathrm{div} \bigl( \alpha(\eta) \mathcal{B}^* +\kappa^2 \nabla {\bf u} \bigr) \in L^2(0, T; L^1(\Omega; \R^{M})),
            \\[1ex]
            \mathrm{div} \bigl( \alpha(\eta) \mathcal{B}^* +\kappa^2 \nabla {\bf u} \bigr) \wedge {\bf u} \in L^2(0, T; L^2(\Omega; \Lambda_2(\R^M)).
        \end{cases}
    \end{equation*}
Moreover, if ${\bf u}_{0} \in \mathbb{S}^{M-1}_{+,r}$ in $\Omega$ for some  $r \in (0,1)$, then ${\bf u} \in \mathbb{S}^{M-1}_{+,r}$ a.e. in $\Omega_{T}$, where:
{
\begin{equation}\label{constraint_S}
    \mathbb{S}^{M-1}_{+,r} := \left\{ \begin{array}{l|l}
        {\bf p} =(p_{1}, \ldots, p_{M}) \in \mathbb{S}^{M-1}_{+} & r \le p_{1}
    \end{array} \right\}, ~\mbox{for any $ r \in (0, 1) $.}
\end{equation}
}

\section{Approximating problem}\label{sec:approx_system}

    In this Section, we consider the approximating problem to our system. Let us assume (A0)--(A6), and fix constants $ \varepsilon \geq 0 $, $ \delta \geq 0 $, and $ \nu \geq 0 $. On this basis, the approximating problem consists in the following system of parabolic PDEs, denoted by (P)$_{{\varepsilon, \nu, \delta}}^\kappa$:

\paragraph{\boldmath Problem (P)$_{{\varepsilon, \nu, \delta}}^\kappa$:}{
    \begin{align*}
        & \begin{cases}
            \partial_t \eta -\mathit{\Delta} \eta +g(\eta) +\alpha'(\eta) f_\varepsilon(\nabla {\bf u}) = 0, ~\mbox{in $ \Omega_T $,}
            \\[0.5ex]
            \nabla \eta \cdot {{\bf n}_\Gamma} = 0 \mbox{ on $ \Gamma_T $,}
            \\[0.5ex]
            \eta(0, x) = \eta_0(x), ~ x \in \Omega;
        \end{cases}
        \\[1ex]
        & \begin{cases}
            \partial_t {\bf u} -\mathrm{div} \bigl( \alpha(\eta) \partial f_\varepsilon(\nabla {\bf u}) +{\kappa^2} \nabla {\bf u} {+\nu |\nu\nabla {\bf u}|^{N -1} \nu\nabla {\bf u} }\bigr) +\varpi_\delta({\bf u})  \ni 0 \mbox{ in $ \Omega_T $,}
            \\[0.5ex]
            \bigl( \alpha(\eta) \partial f_\varepsilon(\nabla {\bf u}) +{\kappa^2} \nabla {\bf u} {+\nu |\nu\nabla {\bf u}|^{N -1} \nu\nabla {\bf u} }\bigr) {{\bf n}_\Gamma} {\,\ni\,} 0 \mbox{ on $ \Gamma_T $,}
            \\[0.5ex]
            {\bf u}(0, x) = {\bf u}_0(x), ~ x \in \Omega;
        \end{cases}
    \end{align*}
}
The approximating problem (P)$ _{{\varepsilon, \nu, \delta}}^\kappa $ is derived as a gradient descent flow of a {free energy,} which is defined as:
    \begin{align*}
        \mathcal{F} & _{{\varepsilon, \nu, \delta}}^\kappa :  U := [\eta, {\bf u}] \in \mathfrak{X} \mapsto \mathcal{F}_{{\varepsilon, \nu, \delta}}^\kappa( U) = \mathcal{F}_{{\varepsilon, \nu, \delta}}^\kappa(\eta, {\bf u})
        \\
        & := \Psi_0(\eta) +\Psi_{{\varepsilon, \nu, \delta}}^\kappa(U) = \Psi_0(\eta) +\Psi_{{\varepsilon, \nu, \delta}}^\kappa(\eta, {\bf u}),
    \end{align*}
    with
    \begin{align*}
        \Psi_0 : ~& \eta \in D(\Psi_0) := V \subset H \mapsto  \Psi_0(\eta) := \frac{1}{2} \int_\Omega |\nabla \eta|^2 \, dx +\int_\Omega G(\eta) \, dx,
    \end{align*}
    and
    \begin{align*}
        \Psi _{{\varepsilon, \nu, \delta}}^\kappa &  : U := [\eta, {\bf u}] \in D(\Psi_{{\varepsilon, \nu, \delta}}^\kappa) := {
            \left\{ \begin{array}{l|l}
                [\tilde{\eta}, \tilde{{\bf u}}] \in \mathfrak{W} & \parbox{4cm}{$ \tilde{{\bf u}} \in L^4(\Omega; \R^M) $,  and $ \nu \tilde{{\bf u}} \in W^{1, N +1}(\Omega; \R^{M}) $}
            \end{array} \right\}
        }
        \\
        & \mapsto \Psi_{{\varepsilon, \nu, \delta}}^\kappa( U) = \Psi_{{\varepsilon, \nu, \delta}}^\kappa(\eta, {\bf u}):= \int_\Omega \alpha(\eta) f_\varepsilon(\nabla {\bf u}) \, dx +\frac{{\kappa^2}}{2} \int_\Omega |\nabla {\bf u}|^2 \, dx
        \\
        & \qquad  {+\frac{1}{N +1} \int_\Omega |\nu \nabla {\bf u}|^{N +1} \, dx} +\int_\Omega \Pi_\delta({\bf u}) \, dx \in [0, \infty).
    \end{align*}

\begin{defn}\label{Def.apKWC}
    A pair of functions $ U := [\eta, {\bf u}] \in {L^2(0, T; \mathfrak{X})} $ is called a solution to the system (P)$_{{\varepsilon, \nu, \delta}}^\kappa$, iff.:
    \begin{equation}\label{apKWC00}
        U = [\eta, {\bf u}] \in {W^{1, 2}(0, T; \mathfrak{X}) \cap L^\infty(0, T; \mathfrak{W});}
    \end{equation}
    \begin{align*}
        \bigl( \partial_t \eta(t) & +g(\eta(t)) +\alpha'(\eta(t)) f_\varepsilon(\nabla {\bf u}(t)), \varphi \bigr)_{H} +\bigl( \nabla \eta(t), \nabla \varphi \bigr)_{H} = 0,
        \\
        & \mbox{for any $ \varphi \in V $, a.e. $ {t \in (0, T)} $, ~subject to $ \eta(0) = \eta_0 $ in $ H $;}
    \end{align*}
    and there exists $ \mathcal{B}^* \in L^\infty({\Omega_T}; \mathbb{R}^{MN}) $, such that:
    \begin{align*}
        & \mathcal{B}^* \in \partial f_\varepsilon(\nabla {\bf u}) ~\mbox{in $ \R^{MN} $, a.e. in $ {\Omega_T} $,}
    \end{align*}
    \begin{align*}
        &\left( \partial_t {\bf u}(t) -\frac{1}{\delta} {\bf u}(t), \bm{\psi} \right)_{{\mathbb{X}}} +\int_\Omega \bigl( \alpha(\eta(t)) \mathcal{B}^*(t) +\kappa^2 \nabla {\bf u} \bigr) : \nabla \bm{\psi} \, dx
        \\
        &  \qquad {+\nu \int_\Omega |\nu \nabla {\bf u}(t)|^{N -1} \nu \nabla {\bf u}(t) : \nabla\bm{\psi} \, dx} +\frac{1}{\delta}\int_\Omega |{\bf u}(t)|^2 {\bf u}(t) \cdot \bm{\psi} \, dx = 0,
        \\
        \mbox{for}~ & \mbox{any $ \bm{\psi} \in {W^{1, N +1}(\Omega; \R^M)} $, a.e. $ {t \in (0, T)} $, ~subject to $ {\bf u}(0) = {\bf u}_0 $ in $ {\mathbb{X}} $.}
    \end{align*}
\end{defn}

    Additionally, we note that our system (P)$_{{\varepsilon, \nu, \delta}}^\kappa$ can be reformulated {as} the following Cauchy problem of evolution {equations} in the Hilbert space $ \mathfrak{X} $, denoted by (CP)$_{{\varepsilon, \nu, \delta}}^\kappa$.
\paragraph{\boldmath Cauchy problem (CP)$_{{\varepsilon, \nu, \delta}}^\kappa$:}{
    \begin{align*}
        & \begin{cases}
            { U}'(t) +\partial \Phi_{{\varepsilon, \nu, \delta}}^\kappa({ U}(t)) +{\mathcal{G}_{\delta}^{\kappa}}({ U}(t)) \ni 0 \mbox{ in $ \mathfrak{X} $, ~ $ {t \in (0, T)} $,}
            \\[1ex]
            { U}(0) = { U}_0 \mbox{ in $ \mathfrak{X} $.}
        \end{cases}
    \end{align*}
    }
    In this context, ``\,$'$\,'' denotes the time-derivative ``\,$\frac{d}{dt}$\,'' of an $ \mathfrak{X} $-valued function (in time). For every $ \kappa, {\varepsilon, \nu, \delta} > 0 $, $ \Phi_{{\varepsilon, \nu, \delta}}^\kappa : \mathfrak{X} \longrightarrow [0, \infty] $ is a proper l.s.c. and convex function, defined as:
    \begin{align}
        \Phi_{{\varepsilon, \nu, \delta}}^\kappa\,&  : { U} := [\eta, {\bf u}] \in D(\Phi_{{\varepsilon, \nu, \delta}}^\kappa) { \subset \mathfrak{X}}
        \nonumber
        \\
        & \mapsto \Phi_{{\varepsilon, \nu, \delta}}^\kappa({ U}) = \Phi_{{\varepsilon, \nu, \delta}}^\kappa(\eta, {\bf u}) := \frac{1}{2} \int_\Omega |\nabla \eta|^2 \, dx { +\frac{1}{N +1} \int_\Omega |\nu \nabla {\bf u}|^{N +1} \, dx }
        \nonumber
        \\
        & \qquad +\frac{1}{2} \int_\Omega \left( {\kappa} f_\varepsilon(\nabla {\bf u}) +\frac{1}{{\kappa}} \alpha(\eta) \right)^2  dx +\frac{1}{4\delta} \int_\Omega |{\bf u}|^4 \, dx \in [0, \infty),
        \label{defAppPhi}
    \end{align}
    and $ {\mathcal{G}_{\delta}^{\kappa}} : \mathfrak{X} \longrightarrow \mathfrak{X} $ is a non-monotone perturbation, defined as:
    \begin{align*}
        {\mathcal{G}_{\delta}^{\kappa}}: ~& { U} := [\eta, {\bf u}] \in  \mathfrak{X} \mapsto {\mathcal{G}_{\delta}^{\kappa}}({ U}) = {\mathcal{G}_{\delta}^{\kappa}}(\eta, {\bf u}) := \left[ \begin{array}{c}
            \ds g(\eta) -\frac{1}{{\kappa^2}} \alpha(\eta) \alpha'(\eta),  ~
            \ds -\frac{1}{\delta}{\bf u}
        \end{array} \right] \in \mathfrak{X}.
    \end{align*}

\begin{rem}\label{Rem.CP}
    {Assumptions} (A2) and (A3) guarantee the Lipschitz continuity of the perturbation $ {\mathcal{G}_{\delta}^{\kappa}} $. Hence, the well-posedness of the Cauchy problem (CP)$_{{\varepsilon, \nu, \delta}}^\kappa$ is immediately verified by means of the general theory of nonlinear evolution equations, such as \cite{MR0348562,MR2582280}. Moreover, we will observe the continuous dependence of solution with respect to $ \varepsilon \geq 0 $, $ \delta > 0 $, $ \nu \geq 0 $, and $ \kappa > 0 $, by means of the general theory of operator-convergence as in \cite{MR0773850,MR0298508}.
\end{rem}

    Now, we set the goal of this Section to prove the following Theorems \ref{Th.KWC-equiv-CP} and \ref{Th.SolvApKWC}.
\begin{thm}\label{Th.KWC-equiv-CP}
    The system (P)$_{{\varepsilon, \nu, \delta}}^\kappa$ is equivalent to the Cauchy problem  (CP)$_{{\varepsilon, \nu, \delta}}^\kappa$.
\end{thm}

\begin{thm}\label{Th.SolvApKWC}
    {Let us assume (A0)--(A6), and let us fix the constants $ \varepsilon \geq 0 $, $ \delta > 0 $, $ \nu \geq 0 $. Additionally, we assume that $ {\bf u}_0 \in L^4(\Omega; \R^M) $ and $ \nu {\bf u}_0 \in W^{1, N +1}(\Omega; \R^{M}) $, for the component $ {\bf u}_0 $ of the initial data $ U_0 = [\eta_0, {\bf u}_0] $. Then, the following items hold:}
    \begin{description}
        \item[(I)]The system (P)$_{{\varepsilon, \nu, \delta}}^\kappa$ admits a unique solution {$
            U = [\eta, {\bf u}] \in {L^2(0, T; \mathfrak{X})} $, such that:
            \begin{equation}\label{CPreg}
                \begin{cases}
                    U = [\eta, {\bf u}] \in {W^{1, 2}(0, T; \mathfrak{X}) \cap L^\infty(0, T; \mathfrak{W}),}
                    \\[1ex]
                    \nu {\bf u} \in {L^{N +1}(0, T; W^{1, N +1}(\Omega; \mathbb{R}^{M})).}
                \end{cases}
            \end{equation}
                }
        \item[(II)]let $ \{ \varepsilon_n \}_{n = 1}^\infty \subset [0, \infty) $, $ \{ \delta_n \}_{n = 1}^\infty \subset (0, \infty)  $, {$ \{\nu_n\}_{n = 1}^\infty \subset [0, \infty) $, } and $ \{ \kappa_n \}_{n = 1}^\infty \subset (0, \infty)  $ be sequences of constants, such that:
            \begin{equation}\label{eps-dlt-nu}
                [\kappa_n, \varepsilon_n, \nu_n, \delta_n] \to [\kappa, \varepsilon, \nu, \delta] \mbox{ {in $ \R^4 $,} as $ n \to \infty $.}
            \end{equation}
            Let $  U_0 = [\eta_0, {\bf u}_0] \in \mathfrak{W} $ be the pair of initial data as in (A6), and let {$ \{U_{0, n}\}_{n = 1}^\infty = \{ [\eta_{0, n}, {\bf u}_{0, n}] \}_{n = 1}^\infty \subset \mathfrak{W} \cap L^4(\Omega; \R^M) $ with $ \{ \nu_n \mathbf{u}_{0, n} \} \subset W^{1, N +1}(\Omega; \R^M) $} be a sequence of initial data, such that:
            \begin{align}
                & U_{0, n} = [\eta_{0, n}, {\bf u}_{0, n}] \to  U_0 = [\eta_0, {\bf u}_0] \mbox{ in $ \mathfrak{X} $, and weakly in $ \mathfrak{W} $,}
                \nonumber
                \\
                & {\mbox{and  } \nu_n {\bf u}_{0, n} \to \nu {\bf u}_0 \mbox{ weakly in $ W^{1, N +1}(\Omega; \mathbb{R}^{M}) $, as $ n \to \infty $.}
                \label{cv00}
                }
            \end{align}
            Let $ U = [\eta, {\bf u}] \in {L^2(0, T; \mathfrak{X})} $ be the solution to the system (P)$_{{\varepsilon, \nu, \delta}}^\kappa $. Also, for any $ n \in \N $, let $ { U}_n = [\eta_n, {\bf u}_n] \in {L^2(0, T; \mathfrak{X})} $ be the solution to the system {(P)$_{\varepsilon_n, \nu_n, \delta_n}^{\kappa_n}$} corresponding to the initial data $ {U}_{0, n} = [\eta_{0, n}, {\bf u}_{0, n}] \in \mathfrak{W} $. Then, it holds that:
        \begin{align}
            { U}_n =\,&  [\eta_n, {\bf u}_n] \to { U} = [\eta, {\bf u}] \mbox{ in $ {C([0, T]; \mathfrak{X})} $, in $ {L^2(0, T; \mathfrak{W})} $,}
            \nonumber
            \\
            & \mbox{weakly in $ {W^{1, 2}(0, T; \mathfrak{X})} $, and weakly-$*$ in $ {L^\infty(0, T; \mathfrak{W})} $,}
            \nonumber
            \\
            {\mbox{and }~}& { \nu_nu_n \to \nu {\bf u} \mbox{ in $ {L^{N +1}(0, T; W^{1, N +1}(\Omega; \mathbb{R}^{M}))} $, as $ n \to \infty $.}}
            \label{cv01}
        \end{align}
    \end{description}
\end{thm}


\subsection{Proofs of Theorems \ref{Th.KWC-equiv-CP} and \ref{Th.SolvApKWC}}

For the proofs of Theorems \ref{Th.KWC-equiv-CP} and \ref{Th.SolvApKWC}, we prepare some Lemmas and Remarks.
{
\begin{lem}\label{Lem04}
    Let $ \varepsilon \geq 0 $, $ \delta > 0 $, $ \nu \geq 0 $, and $ \kappa > 0 $ be fixed constants, and let $ \{ \varepsilon_n \}_{n = 1}^\infty \subset [0, \infty) $, $ \{ \delta_n \}_{n = 1}^\infty \subset (0, \infty)  $, $ \{ \nu_n \}_{n = 1}^\infty \subset [0, \infty) $ and $ \{ \kappa_n \}_{n = 1}^\infty \subset (0, \infty)  $ be sequences of constants, as in \eqref{eps-dlt-nu}. Then, the sequence of convex functions $ \{ \Phi_{\varepsilon_n, \nu_n, \delta_n}^{\kappa_n} \}_{n = 1}^{\infty} $ converges to the convex function $ \Phi_{{\varepsilon, \nu, \delta}}^\kappa $ on $ \mathfrak{X} $, in the sense of Mosco (cf. \cite{MR0298508}), as $ n \to \infty $.
\end{lem}
\paragraph{Proof.}{
As is easily checked,
\begin{enumerate}
    \item[$\sharp\,1$)]the sequence $ \bigl\{ {\ts \frac{1}{2} \bigl( {\kappa_n} f_{\varepsilon_n} +\frac{1}{{\kappa_n}} \alpha \bigr)^2} \bigr\}_{n  =1}^{\infty} $ of continuous convex functions:
    \begin{align*}
        & {\ts \frac{1}{2} \bigl( {\kappa_n} f_{\varepsilon_n} +\frac{1}{{\kappa_n}} \alpha \bigr)^2} : [W, \xi] \in \R^{MN} \times \R
        \\
        & \qquad \mapsto \frac{1}{2} \left( {\kappa_n} f_{\varepsilon_n}(W) +\frac{1}{{\kappa_n}} \alpha(\xi) \right)^2 \in [0, \infty), ~ n = 1, 2, 3, \dots,
    \end{align*}
    converges to the continuous convex function:
    \begin{align*}
        & {\ts \frac{1}{2} \bigl( {\kappa_0} f_{\varepsilon_0} +\frac{1}{{\kappa_0}} \alpha \bigr)^2} : [W, \xi] \in \R^{MN} \times \R
        \\
        & \qquad \mapsto \frac{1}{2} \left( {\kappa_0} f_{\varepsilon_0}(W) +\frac{1}{{\kappa_0}} \alpha(\xi) \right)^2 \in [0, \infty),
    \end{align*}
    on $ \R^{MN} \times \R $, in the sense of Mosco, as $ n \to \infty $;
\item[$\sharp\,2$)]the sequence $ \{ \frac{1}{4 \delta_n} \|\cdot\|_{\R^M}^4 \}_{n = 1}^{\infty} $ of continuous convex functions:
    \begin{equation*}
        {\ts\frac{1}{4 \delta_n} \|\cdot\|_{\R^M}^4} : {\bf v} \in \R^M \mapsto \frac{1}{4 \delta_n}\|{\bf v}\|_{\R^M}^4  \in [0, \infty),
    \end{equation*}
    converges to the convex function:
    \begin{equation*}
        {\ts\frac{1}{4 \delta_0} \|\cdot\|_{\R^M}} : {\bf v} \in \R^M \mapsto \frac{1}{4 \delta_0}\|{\bf v}\|_{\R^M}^4 \in [0, \infty),
    \end{equation*}
        on $ \R^M $, in the sense of Mosco, {as $ n \to \infty $;}
    \item[$\sharp\,3$)]the sequence $ \{ \frac{1}{N +1} \|\nu_n \, (\cdot)\|_{\R^M}^{N +1} \}_{n = 1}^{\infty} $ of continuous convex functions:
    \begin{equation*}
        {\ts\frac{1}{N +1} \|\nu_n \, (\cdot)\|_{\R^M}^{N +1}} : {\bf v} \in \R^M \mapsto {\frac{1}{N +1}\|\nu_n {\bf v}\|_{\R^M}^{N +1}}  \in [0, \infty),
    \end{equation*}
    converges to the convex function:
    \begin{equation*}
        {\ts\frac{1}{N +1} \|\nu \, (\cdot)\|_{\R^M}^{N +1}} : {\bf v} \in \R^M \mapsto {\frac{1}{N +1}\|\nu {\bf v}\|_{\R^M}^{N +1}} \in [0, \infty),
    \end{equation*}
        on $ \R^M $, in the sense of Mosco, as $ n \to \infty $.
\end{enumerate}
Taking into account $\sharp\,1$), $\sharp\,2$), {(Fact\,2), and Example \ref{Rem.ExMG01} in Appendix A,} we can verify the condition of lower-bound for $ \{ \Phi_{\varepsilon_n, \nu_n, \delta_n}^{\kappa_n} \}_{n = 1}^\infty $, as follows:
\begin{align*}
    & \varliminf_{n \to \infty} \Phi_{\varepsilon_n, \nu_n, \delta_n}^{\kappa_n}(\eta_n, {\bf u}_n) \geq \varliminf_{n \to \infty} \frac{1}{2} \int_\Omega |\nabla \eta_n|^2 \, dx
    \\
    & \qquad +\varliminf_{n \to \infty} \int_\Omega \frac{1}{2} \left( {\kappa_n} f_{\varepsilon_n} (\nabla {\bf u}_n) +\frac{1}{{\kappa_n}} \alpha(\eta_n) \right)^2 \, dx
    \\
    & \qquad +\varliminf_{n \to \infty}\int_\Omega \frac{1}{4\delta_n} |{\bf u}_n|^4 \, dx +\varliminf_{n \to \infty}\int_\Omega \frac{1}{N +1}\bigl| \nu_n \nabla {\bf u}_n \bigr|^{N +1} \, dx
    \\
    & \geq \frac{1}{2} \int_\Omega |\nabla \eta|^2 \, dx +\int_\Omega \frac{1}{2} \left( {\kappa_0} f_{\varepsilon_0}(\nabla {\bf u}) +\frac{1}{{\kappa_0}} \alpha(\eta) \right)^2 \, dx
    \\
    & \qquad +\int_\Omega \frac{1}{4\delta_0} |{\bf u}|^4 \, dx +\int_\Omega \frac{1}{N +1} \bigl| \nu \nabla {\bf u}  \bigr|^{N +1} \, dx
    \\
    & = \Phi_{{\varepsilon, \nu, \delta}}^\kappa(\eta, {\bf u}),
\end{align*}
for all $ [\eta, {\bf u}] \in \mathfrak{X} $, $ \{ [\eta_n, {\bf u}_n] \}_{n = 1}^\infty \subset \mathfrak{X} $, fulfilling $ [\eta_n, {\bf u}_n] \to [\eta, {\bf u}] $ weakly in $ \mathfrak{X} $ as $ n \to \infty $.
\bigskip

In the meantime, the condition of optimality for $ \{ \Phi_{\varepsilon_n, \nu_n, \delta_n}^{\kappa_n} \}_{n = 1}^\infty $ is observed, via the convergence:
\begin{equation*}
    \lim_{n \to \infty} \Phi_{\varepsilon_n, \nu_n, \delta_n}^{\kappa_n}(\eta, {\bf u}) = \Phi_{{\varepsilon, \nu, \delta}}^\kappa(\eta, {\bf u}),
\end{equation*}
that is verified by means of Lebesgue's dominated convergence theorem (cf. \cite[Theorem 10 on page 36]{MR0492147}). \hfill \qed
}
\begin{rem}\label{Rem_Lem04}
    Let $ \varepsilon \geq 0 $, $ \delta > 0 $, $ \nu \geq 0 $, and $ \kappa > 0 $ be as in Lemma \ref{Lem04}, and let $ \{ \varepsilon_n \}_{n = 1}^\infty \subset (0, \infty) $, $ \{ \delta_n \}_{n = 1}^\infty \subset (0, \infty) $, $ \{ \nu_n \}_{n = 1}^\infty \subset (0, \infty) $, and $ \{ \kappa_n \}_{n = 1}^\infty \subset (0, \infty) $ be the sequences as in \eqref{eps-dlt-nu}. Additionally, let $ \tilde{\eta} \in H $ be a fixed function, and let $ \{ \tilde{\eta}_n \}_{n = 1}^\infty \subset H $ be a sequence such that:
    \begin{equation}\label{tildeEta}
        \tilde{\eta}_n \to \tilde{\eta} ~\mbox{{in $ V $, as $ n \to \infty $.}}
    \end{equation}
    Then, as a corollary of Lemma \ref{Lem04}, we can derive:
    \begin{align}\label{MoscoRL04}
        & \Phi_{\varepsilon_n, \nu_n, \delta_n}^{\kappa_n}(\tilde{\eta}_n, \cdot) \to \Phi_{{\varepsilon, \nu, \delta}}^\kappa(\tilde{\eta}, \cdot) ~\mbox{on $ \mathbb{X} $, in the sense of Mosco, as $ n \to \infty $.}
    \end{align}
    In fact, the condition of lower-bound for the Mosco convergence \eqref{MoscoRL04} will be a straightforward consequence of Lemma \ref{Lem04}. Also, for any $ {\bf u} \in D(\Phi_{{\varepsilon, \nu, \delta}}^\kappa) $, one can easily see from \eqref{defAppPhi}, \eqref{eps-dlt-nu}, and \eqref{tildeEta} that:
    \begin{align*}
        & \lim_{n \to \infty} \Phi_{\varepsilon_n, \nu_n, \delta_n}^{\kappa_n}(\tilde{\eta}_n, {\bf u})  = {\Phi_{{\varepsilon, \nu, \delta}}^\kappa(\tilde{\eta}, {\bf u}).}
    \end{align*}
    This implies that the sequence $ \{ {\bf u}_n = {\bf u} \}_{n = 1}^\infty $ fulfills the condition of optimality condition for \eqref{MoscoRL04}. \hfill \qed
\end{rem}
}
\begin{lem}\label{Lem01}
    Let $ \varepsilon \geq 0 $, $ \delta > 0 $, $ \nu \geq 0 $, and $ \kappa > 0 $ be fixed constants. Then, for any $ \tilde{{\bf u}} \in {\mathbb{W}} $, the functional $ \Phi_{{\varepsilon, \nu, \delta}}^\kappa : \eta \in H \longrightarrow \Phi_{{\varepsilon, \nu, \delta}}^\kappa(\eta, \tilde{{\bf u}}) \in [0, \infty] $ is proper l.s.c. and convex on $ H $, and the subdifferential $ \partial_\eta \Phi_{{\varepsilon, \nu, \delta}}^\kappa(\cdot, \tilde{{\bf u}}) \subset H \times H $ is a single-valued operator, such that:
    \begin{align*}
        D(\partial_\eta & \Phi_{{\varepsilon, \nu, \delta}}^\kappa(\cdot, \tilde{{\bf u}})) = D(\partial \Psi_0) = \left\{ \begin{array}{l|l}
            \tilde{\eta} \in H^2(\Omega)  & \nabla \eta \cdot {{\bf n}_\Gamma} = 0 \mbox{ on $ \Gamma $}
        \end{array} \right\},
        \\
        & \mbox{and }~ \partial_\eta \Phi_{{\varepsilon, \nu, \delta}}^\kappa(\eta, \tilde{{\bf u}}) = -\mathit{\Delta} \eta +\alpha'(\eta)f_\varepsilon(\nabla \tilde{{\bf u}}) +\frac{1}{\kappa^2} \alpha(\eta)\alpha'(\eta), ~~ \forall \eta \in D(\Phi_{{\varepsilon, \nu, \delta}}^\kappa).
    \end{align*}
\end{lem}
\paragraph{Proof.}{
    This Lemma \ref{Lem01} is a straightforward consequence of the standard variational method (cf. \cite{MR0348562,MR2582280}). \qed
}
\begin{lem}\label{Lem02}
    Let $ \varepsilon \geq 0 $, $ \delta > 0 $, {$ \nu \geq 0 $,} and $ \kappa > 0 $ be fixed constants. Then, for any $ \tilde{\eta} \in V $, the functional $ \Phi_{{\varepsilon, \nu, \delta}}^\kappa(\tilde{\eta}, \cdot) : {\bf u} \in {\mathbb{X}} \mapsto \Phi_{{\varepsilon, \nu, \delta}}^\kappa(\tilde{\eta}, {\bf u}) \in [0, \infty] $ is proper l.s.c. and convex on $ {\mathbb{X}} $. Moreover, concerning the subdifferential $ \partial_{\bf u} \Phi_{{\varepsilon, \nu, \delta}}^\kappa(\tilde{\eta}, \cdot) \subset {\mathbb{X}} \times {\mathbb{X}} $, the following two items are equivalent:
    \begin{description}
        \item[\textmd{\it(O)}]$ [{\bf u}, {\bf u}^*] \in \partial_{\bf u} \Phi_{{\varepsilon, \nu, \delta}}^\kappa(\tilde{\eta}; \cdot) $ in $ {\mathbb{X}} \times {\mathbb{X}} $;
        \item[\textmd{\it(I)}]$ [{\bf u}, {\bf u}^*] \in \mathbb{X} \times \mathbb{X} $, $ {\bf u} \in {D(\Phi_{{\varepsilon, \nu, \delta}}^\kappa(\tilde{\eta}, \cdot)) = \left\{ \begin{array}{l|l}
                {\bf u} \in \mathbb{W} \cap L^4(\Omega; \R^{M}) & \nu {\bf u} \in W^{1, N +1}(\Omega; \R^{M})
        \end{array} \right\}} $, and there exists $ W_{\bf u}^* \in L^\infty(\Omega; \R^{MN}) $ such that
    \item[~~~~~~~\textmd{\it(i-a)}]$W_{\bf u}^* \in \partial f_\varepsilon(\nabla {\bf u})$ a.e. in $ \Omega $,
    \item[~~~~~~~\textmd{\it(i-b)}]$ \ds \mathrm{div} \bigl( \alpha(\tilde{\eta}) W_{\bf u}^* +{\kappa^2} \nabla {\bf u} {+\nu |\nu \nabla {\bf u}|^{N -1} \nu \nabla {\bf u}} \bigr) +\frac{1}{\delta} |{\bf u}|^2 {\bf u} \in {\mathbb{X}} $,
    \item[~~~~~~~\textmd{\it(i-c)}]$ \ds \bigl( \alpha(\tilde{\eta}) W_{\bf u}^* +{\kappa^2} \nabla {\bf u} {+\nu |\nu \nabla {\bf u}|^{N -1} \nu \nabla {\bf u}} \bigr) {\bf n}_\Gamma = 0 $ in $ L^2(\Gamma; \R^{M}) $,
    \item[~~~~~~~\textmd{\it(i-d)}]$ \ds {\bf u}^* = -\mathrm{div} \bigl( \alpha(\tilde{\eta}) W_{\bf u}^* +{\kappa^2} \nabla {\bf u} {+\nu |\nu \nabla {\bf u}|^{N -1} \nu \nabla {\bf u}} \bigr) +\frac{1}{\delta} |{\bf u}|^2 {\bf u} $ in $ {\mathbb{X}} $.
    \end{description}
\end{lem}
\paragraph{Proof.}{
{
    When $ \varepsilon > 0 $, this Lemma will be concluded, immediately, by applying the theory of elliptic variational inequalities (cf. \cite{MR0348562,MR0390843,MR2582280}). Additionally, it is observed that the operator $ \partial_{\bf u} \Phi_{{\varepsilon, \nu, \delta}}^\kappa(\tilde{\eta}, \cdot) \subset \mathbb{X} \times \mathbb{X} $ is single-valued, and
    \begin{align}
        & \begin{cases}
            \partial_{\bf u} \Phi_{{\varepsilon, \nu, \delta}}^\kappa(\tilde{\eta}, {\bf u}) = -\mathrm{div} \bigl( \alpha(\tilde{\eta}) [\nabla f_\varepsilon](\nabla {\bf u}) +{\kappa^2} \nabla {\bf u} {+\nu |\nu \nabla {\bf u}|^{N -1} \nu \nabla {\bf u}} \bigr)
            \\
            \qquad \ds +\frac{1}{\delta} |{\bf u}|^2 {\bf u} ~\mbox{in $ \mathbb{X} $,}
        \\[1ex]
        \mbox{subject to}~~\nabla {\bf u} \, {\bf n}_\Gamma = 0 ~\mbox{in $ L^2(\Gamma; \R^{N}) $.}
        \end{cases}
        \label{Lem02-01}
    \end{align}
    Hence, we here focus on the operator $ \partial_{\bf u} \Phi_{0, \nu, \delta}^\kappa(\tilde{\eta}) $ corresponding to the case when $ \varepsilon = 0 $.
    \medskip

    Let us fix arbitrary $ \delta > 0 $, $ \nu \geq 0 $, and $ \kappa > 0 $, and let us define a set-valued map $ {\mathcal{A}_{\nu, \delta}^{\kappa}} : \mathbb{X} \to 2^{\mathbb{X}} $, by putting:
\begin{equation}\label{key1.01}
    D({\mathcal{A}_{\nu, \delta}^{\kappa}}) := \left\{\begin{array}{l|l} {\bf u} \in D(\Phi_{0, \nu, \delta}^\kappa(\tilde{\eta}, \cdot)) & \begin{tabular}{l} \mbox{there exists $ W_{\bf u}^* \in L^{\infty}(\Omega; \R^{MN}) $}
\\[1ex]
\mbox{such that (i-a)--(i-c) hold}. \end{tabular} \end{array}\right\},
\end{equation}
and
\begin{equation}\label{key1.02}
\begin{array}{rcl}
    {\bf u} & \in & D({\mathcal{A}_{\nu, \delta}^{\kappa}}) \subset \mathbb{X} \mapsto {\mathcal{A}_{\nu, \delta}^{\kappa}}{\bf u}
\\[2ex]
    & := &\left\{\begin{array}{l|l} {\bf u}^* \in \mathbb{X}
& \begin{tabular}{l}
    (i-d) holds, for some $ W_{\bf u}^* \in L^{\infty}(\Omega; \R^{MN}) $,
    \\[1ex]
    satisfying (i-a)--(i-c) \end{tabular} \end{array}\right\}.
\end{array}
\end{equation}
Then, the assertion of Lemma \ref{Lem02} can be rephrased as follows:
\begin{equation}\label{key1.03}
    \partial_{\bf u} \Phi_{0, \nu, \delta}^\kappa(\tilde{\eta}, \cdot) = {\mathcal{A}_{\nu, \delta}^{\kappa}} \mbox{ in $ \mathbb{X} \times \mathbb{X} $.}
\end{equation}
This coincidence will be obtained as a consequence of the following {\it Claims $ \# $1--$ \# $2}.
\medskip

\noindent
    {\it Claim $ \# $1: $ {\mathcal{A}_{\nu, \delta}^{\kappa}} \subset \partial_{\bf u} \Phi_{0, \nu, \delta}^\kappa $ in $ \mathbb{X} \times \mathbb{X} $.}

    Let us assume that $ {\bf u} \in D({\mathcal{A}_{\nu, \delta}^{\kappa}}) $ and $ {\bf u}^* \in {\mathcal{A}_{\nu, \delta}^{\kappa}} {\bf u} $ in $ \mathbb{X} $. Then, from (A4) and (i-a)--(i-d), it is inferred that:
\begin{align*}
    \bigl( {\bf u}^*, {\bf z} -{\bf u} \bigr)_{\mathbb{X}} ~& = \bigl< -\mathrm{div} \bigl( \alpha(\tilde{\eta}) W_{\bf u}^* +\kappa^2 \nabla {\bf u} +\nu |\nu \nabla {\bf u}|^{N -1} \nu \nabla {\bf u}), {\bf z} -{\bf u}  \bigr>_{\mathbb{W}}
    \\
    & \qquad +\frac{1}{\delta} \bigl< |{\bf u}|^2 {\bf u}, {\bf z} -{\bf u} \bigr>_{L^4(\Omega; \R^{MN})}
    \\
    & = \int_\Omega \alpha(\tilde{\eta}) W_{\bf u}^* : \nabla ({\bf z} -{\bf u}) \, dx +\kappa^2 \int_\Omega \nabla {\bf u} : \nabla ({\bf z} -{\bf u}) \, dx
    \\
    &  \qquad +\nu \int_\Omega |\nu \nabla {\bf u}|^{N -1} \nu \nabla {\bf u} :  \nabla ({\bf z} -{\bf u}) \, dx +\frac{1}{\delta} \int_\Omega  |{\bf u}|^2 {\bf u} \cdot ({\bf z} -{\bf u}) \, dx
    \\
    & \leq \int_\Omega \alpha(\tilde{\eta}) \bigl( |\nabla {\bf z}| -|\nabla {\bf u}| \bigr) \, dx +\frac{\kappa^2}{2} \int_\Omega \bigl( |\nabla {\bf z}|^2 -|\nabla {\bf u}|^2 \bigr) \, dx
    \\
    & \qquad +\frac{1}{N +1} \int_\Omega \bigl( |\nu \nabla {\bf z}|^{N +1} -|\nu \nabla {\bf u}|^{N +1} \bigr) \, dx +\frac{1}{4\delta} \int_\Omega \bigl( |{\bf z}|^4 -|{\bf u}|^4 \bigr) \, dx
    \\
    & = \frac{1}{2} \int_\Omega \left( \kappa |\nabla {\bf z}| +\frac{1}{\kappa} \alpha(\tilde{\eta}) \right)^2 \, dx -\frac{1}{2} \int_\Omega \left( \kappa |\nabla {\bf u}| +\frac{1}{\kappa} \alpha(\tilde{\eta}) \right)^2 \, dx
    \\
    & \qquad +\frac{1}{N +1} \int_\Omega |\nu \nabla {\bf z}|^{N +1} \, dx -\frac{1}{N +1} \int_\Omega |\nu \nabla {\bf u}|^{N +1} \, dx
    \\
    & \qquad +\frac{1}{4\delta} \int_\Omega |{\bf z}|^4 \, dx -\frac{1}{4\delta} \int_\Omega |{\bf u}|^4 \, dx
    \\
    & = \Phi_{0, \nu, \delta}^\kappa(\tilde{\eta}, {\bf z}) -\Phi_{0, \nu, \delta}^\kappa(\tilde{\eta}, {\bf u}), ~\mbox{for any $ {\bf z} \in D(\Phi(\tilde{\eta}; \cdot)) $.}
\end{align*}
Thus, we have:
\begin{equation*}
    {\bf u} \in D(\partial_{\bf u} \Phi_{0, \nu, \delta}^\kappa(\tilde{\eta}, \cdot)) \mbox{ and } {\bf u}^* \in \partial_{\bf u} \Phi_{0, \nu, \delta}^\kappa(\tilde{\eta}, {\bf u}) \mbox{ in $ \mathbb{X} $},
\end{equation*}
and we can say that:
\begin{equation*}
    {\mathcal{A}_{\nu, \delta}^{\kappa}} \subset \partial_{\bf u} \Phi_{0, \nu, \delta}^\kappa(\tilde{\eta}, \cdot)\mbox{ in $ \mathbb{X} \times \mathbb{X} $.}
\end{equation*}
\medskip

\noindent
    {\it Claim $ \# $2: $ ({\mathcal{A}_{\nu, \delta}^{\kappa}} + \mathcal{I}_{\mathbb{X}}) \mathbb{X} = \mathbb{X} $.}

    Since $ ({\mathcal{A}_{\nu, \delta}^{\kappa}} + \mathcal{I}_{\mathbb{X}}) \mathbb{X} \subset \mathbb{X} $ is trivial, it is sufficient to prove the converse inclusion. Let us take any $ \bf{w} \in \mathbb{X} $. Then, by applying Minty's theorem, we can find a net of functions $ \{ {\bf u}_{\varepsilon} \,|\, 0 < \varepsilon \leq 1 \} \subset \mathbb{W} $, by setting:
\begin{equation*}
    {\left\{ \begin{array}{l|l}
        {\bf u}_\varepsilon := (\partial_{\bf u} \Phi_{{\varepsilon, \nu, \delta}}^\kappa(\tilde{\eta}, \cdot)  + \mathcal{I}_{\mathbb{X}})^{-1}{\bf w} & 0  < \varepsilon \le 1 
    \end{array} \right\}} \mbox{ in $ \mathbb{X} $, }
\end{equation*}
i.e.
\begin{equation}\label{key1.04}
    {\bf w} -{\bf u}_{\varepsilon} \in \partial_{\bf u} \Phi_{{\varepsilon, \nu, \delta}}^\kappa(\tilde{\eta}, {\bf u}_\varepsilon)
 \mbox{ in $ \mathbb{X} $, for any $ 0 < \varepsilon \le 1 $,}
\end{equation}
and in the light of \eqref{Lem02-01}, we can also see that:
\begin{align}\label{key1.05}
    & \int_\Omega \bigl(\alpha(\tilde\eta) [\nabla f_\varepsilon](\nabla {\bf u}_\varepsilon) +\kappa^2 \nabla {\bf u}_\varepsilon +\nu |\nu \nabla {\bf u}_\varepsilon|^{N -1} \nu \nabla {\bf u}_\varepsilon \bigr) : \nabla {\bf z} \, dx +\frac{1}{\delta} \int_\Omega |{\bf u}_\varepsilon|^2 {\bf u}_\varepsilon \cdot {\bf z} \, dx
    \nonumber
    \\
    = & \int_\Omega  ({\bf w} -{\bf u}_\varepsilon) \cdot {\bf z} \, dx, ~\mbox{for all $ {\bf z} \in D(\Phi_{0, \nu, \delta}^\kappa(\tilde{\eta}, \cdot)) $~$ \bigl( = D(\Phi_{{\varepsilon, \nu, \delta}}^\kappa(\tilde{\eta}, \cdot)) \bigr) $, and $ \varepsilon > 0 $.}
\end{align}
In the variational form \eqref{key1.05}, let us take $ {\bf z} = {\bf u}_{\varepsilon} \in D(\Phi_{{\varepsilon, \nu, \delta}}^\kappa(\overline \eta,\cdot)) $. Then, with (A1), (A4)--(A5) and Young's inequality in mind, we deduce that:
\begin{align}\label{key1.06}
    \frac{1}{2} \bigl\|{\bf u}_{\varepsilon} \bigr\|_{\mathbb{X}}^2 +\kappa^2 \bigl\|\nabla {\bf u}_{\varepsilon} & \bigr\|_{L^2(\Omega; \R^{MN})}^2 +\bigl\| \nu \nabla {\bf u}_{\varepsilon} \bigr\|_{L^{N +1}(\Omega; \R^{MN})}^{N +1} +\frac{1}{\delta} \|{\bf u}_\varepsilon\|_{L^4(\Omega; \R^M)}^4
    \nonumber
\\[0ex]
&  \le \frac{1}{2} \bigl\| {\bf w} \bigr\|_{\mathbb{X}}^2, ~\mbox{for any $ 0 < \varepsilon \le 1$.}
\end{align}

On account of \eqref{key1.06}, we find a function $ {\bf u} \in D(\Phi_{{\varepsilon, \nu, \delta}}^\kappa(\tilde{\eta}, \cdot)) $ and a sequence $  1 > \varepsilon_1 > \varepsilon_2 > \varepsilon_3 > \dots > \varepsilon_n \downarrow 0$ as $ n \to \infty $, such that:
\begin{equation}\label{key1.10}
    \begin{cases}
        \u_n := {\bf u}_{\varepsilon_n} \to {\bf u} ~\mbox{in $ \mathbb{X} $, weakly in $\mathbb{W}$, weakly in $ L^4(\Omega; \R^{M}) $,}
        \\
        \qquad \mbox{in the pointwise sense a.e. in $ \Omega $,}
        \\[1ex]
        \nu \nabla {\bf u}_{n} \to \nu \nabla {\bf u} ~\mbox{weakly in $L^{N +1}(\Omega; \R^{MN})$,}
    \end{cases}
    \mbox{ as $ n \to \infty $.}
\end{equation}

Here, in the light of \eqref{key1.10}, \eqref{key1.04} and  Remark \ref{Rem_Lem04}, we can apply Remark \ref{Rem.MG} (Fact\,1) to see that:
\begin{equation*}
    {\bf w} - {\bf u} \in \partial_{\bf u} \Phi_{0, \nu, \delta}^\kappa(\tilde{\eta}, {\bf u}) ~\mbox{in $ \mathbb{X} $,}
\end{equation*}
and
\begin{equation}\label{key1.11}
    \Phi_{\varepsilon_n, \nu, \delta}^\kappa(\tilde{\eta}, {\bf u}_n) \to \Phi_{0, \nu, \delta}^\kappa(\tilde{\eta}, {\bf u}) ~\mbox{as $ n \to \infty $. }
\end{equation}
By virtue of \eqref{defAppPhi}, \eqref{key1.05}--\eqref{key1.11}, $\sharp\,1)$--$\sharp\,3)$, and Remark \ref{Rem.MG} (Fact\,2), we further compute that:
\begin{align}\label{key1.12-1}
    \frac{\kappa^2}{2} \int_{\Omega}  ~& |\nabla {\bf u}|^2 \, dx \leq \frac{\kappa^2}{2} \varliminf_{n \to \infty} \int_{\Omega} |\nabla {\bf u}_n|^2 \, dx \leq \frac{\kappa^2}{2} \varlimsup_{n \to \infty} \int_{\Omega} |\nabla {\bf u}_n|^2 \, dx
    \nonumber
    \\
    & \leq \lim_{n \to \infty} \Phi_{\varepsilon_n, \nu, \delta}^\kappa (\tilde{\eta}, {\bf u}_n) - \varliminf_{n \to \infty} \int_{\Omega} \left( \kappa f_{\varepsilon_n}(\nabla {\bf u}_n) +\frac{1}{\kappa} \alpha(\tilde{\eta}) \right)^2 dx
    \nonumber
    \\
    & \qquad -\varliminf_{n \to \infty} \frac{1}{\delta} \int_{\Omega} \big| {\bf u}_n \big|^4 dx -\varliminf_{n \to \infty} \frac{1}{N +1} \int_{\Omega} \bigl| \nu \nabla {\bf u}_n \bigr|^{N +1} \, dx
    \nonumber
    \\
    & \leq \Phi_{0, \nu, \delta}^\kappa (\tilde{\eta}, {\bf u}) -\int_{\Omega} \left( \kappa |\nabla {\bf u}| +\frac{1}{\kappa} \alpha(\tilde{\eta}) \right)^2 dx
    \nonumber
    \\
    & \qquad -\frac{1}{\delta} \int_{\Omega} \big| {\bf u} \big|^4 dx -\frac{1}{N +1} \int_{\Omega} \bigl| \nu \nabla {\bf u} \bigr|^{N +1} \, dx
    = \frac{\kappa^2}{2} \int_\Omega |\nabla {\bf u}|^2 \, dx.
\end{align}
Having in mind \eqref{key1.10}, \eqref{key1.12-1} and the above calculation and the uniform convexity of $ L^2 $-based topologies, it is deduced that:
\begin{equation}\label{key1.13}
\left\{\begin{array}{l}
    {\bf u}_n \to {\bf u} \mbox{ in $ \mathbb{W} $, in the pointwise sense, a.e. in $ \Omega $,}
\\[0.5ex]
    \nabla {\bf u}_n \to \nabla {\bf u} \mbox{ in $ L^2(\Omega; \R^{MN}) $, in the pointwise sense a.e. in $ \Omega $,}
\end{array} \right.
\mbox{ as $ n \to \infty $.}
\end{equation}
Additionally, invoking the boundedness of $ \Omega $, and Lions's Lemma \cite[Lemma 1.3, in page 12]{MR0259693}, one can derive from \eqref{key1.10} and \eqref{key1.13} that:
\begin{equation*}
    \begin{cases}
        |\u_n|^2 \u_n \to |\u|^2 \u ~\mbox{weakly in $ L^{\frac{4}{3}}(\Omega; \R^{M}) $,}
        \\[1ex]
        |\nu \nabla \u_n|^{N -1} \nu \nabla \u_n \to |\nu \nabla \u|^{N -1} \nu \nabla \u ~\mbox{weakly in $ L^{\frac{N +1}{N}}(\Omega; \R^{MN}) $,}
    \end{cases}
    \mbox{as $ n \to \infty $,}
\end{equation*}
i.e.:
\begin{align}\label{key1.10-02}
    \int_\Omega \nu & |\nu \nabla \u_n|^{N -1} \nu \nabla \u_n : \nabla \mathbf{z} \, dx +\frac{1}{\delta} \int_\Omega |\u_n|^2 \u_n \cdot \mathbf{z} \, dx
    \nonumber
    \\
    & \to \int_\Omega \nu |\nu \nabla \u|^{N -1} \nu \nabla \u : \nabla \mathbf{z} \, dx +\frac{1}{\delta} \int_\Omega |\u|^2 \u \cdot \mathbf{z} \, dx,
    \\
    & \qquad \mbox{for any $ \mathbf{z} \in D(\Phi_{0, \nu, \delta}^\kappa(\tilde{\eta}, \cdot)) $, as $ n \to \infty $.}
    \nonumber
\end{align}

In the meantime, from (A4), one can see that:
\begin{enumerate}
    \item[$\sharp\,4$)]$ \bigl| [\nabla f_{\varepsilon_n}](\nabla \u_n) \bigr| $ ~$ \bigl( = \bigl| \frac{\nabla \u_n}{\sqrt{\varepsilon_n^2 +|\nabla \u_n|^2}} \bigr| \bigr) $~ $ \leq 1$, ~a.e. in $ \Omega $, for any $ n \in \N $,
    \item[$\sharp\,5$)]the sequence of convex functions: 
        $$ \bigl\{ W \in L^2(\Omega; \R^{MN}) \mapsto \bigl\| f_{\varepsilon_n}(W) \bigr\|_{L^1(\Omega)} \in \R \bigr\}_{n = 1}^{\infty} $$
        converges to the convex function:
        $$ W \in L^2(\Omega; \R^{MN}) \mapsto \bigl\| f_0(W) \bigr\|_{L^1(\Omega)} ~\bigl( = \bigl\| W \bigr\|_{L^1(\Omega;\R^{MN})} \bigr),
        $$
        on $ L^2(\Omega; \R^{MN}) $, in the sense of Mosco, as $ n \to \infty $.
\end{enumerate}
\eqref{key1.13} and $ \sharp\,4) $ enables us to say
\begin{align}\label{key1.14}
    [\nabla f_{\varepsilon_n}](\nabla \u_n) \to ~& M_\u ~\mbox{weakly-$*$ in $ L^\infty(\Omega; \R^{MN}) $ as $ n \to \infty $,}
    \nonumber
    \\
    & \mbox{for some $ M_\u \in L^\infty(\Omega; \R^{MN}) $, }
\end{align}
by taking a subsequence if necessary. In view of \eqref{key1.13}, \eqref{key1.14}, $ \sharp\,5) $, Remark \ref{Rem.MG} (Fact\,1), and \cite[Proposition 2.16]{MR0348562}, one can see that:
\begin{equation}\label{key1.15}
M_\u \in \partial f_0 (\nabla \u) \mbox{ a.e. in $ \Omega $.}
\end{equation}

With \eqref{key1.10}, \eqref{key1.10-02}, and \eqref{key1.14} in mind, letting $ n \to \infty $ in \eqref{key1.05} yields that:
\begin{align}\label{key1.16}
    \int_\Omega \bigl(\alpha(\tilde\eta) M_{\bf u} + & \kappa^2 \nabla {\bf u} +\nu |\nu \nabla {\bf u}|^{N -1} \nu \nabla {\bf u} \bigr) : \nabla {\bf z} \, dx +\frac{1}{\delta} \int_\Omega |{\bf u}|^2 {\bf u} \cdot {\bf z} \, dx
    \nonumber
    \\
    = & \int_\Omega  ({\bf w} -{\bf u}) \cdot {\bf z} \, dx, ~\mbox{for all $ {\bf z} \in D(\Phi_{0, \nu, \delta}^\kappa(\tilde{\eta}, \cdot)) $.}
\end{align}
In particular, taking any $ \bm{\varphi}_0 \in C_\mathrm{c}^\infty(\Omega; \R^{M}) $, and putting $ \mathbf{z} = \bm{\varphi}_0 $,
\begin{align*}
    \bigl( \w - \u, \bm{\varphi}_0 \bigr)_{\mathbb{X}} = \int_{\Omega} (\alpha(\tilde \eta) M_\u + & \kappa^2 \nabla \u +\nu |\nu \nabla \u|^{N -1} \nu \nabla \u) : \nabla \bm{\varphi}_0 \, dx +\frac{1}{\delta} \int_\Omega |\u|^2 \u \cdot \bm{\varphi}_0 \, dx,
    \\
    & \mbox{for any $ \bm{\varphi}_0 \in C_\mathrm{c}^\infty(\Omega; \R^{M}) $}
\end{align*}
which implies:
\begin{equation*}
    \begin{array}{c}
    \begin{cases}
        -\mathrm{div} \bigl(\alpha(\tilde \eta) M_\u + \kappa^2 \nabla \u +\nu |\nabla \u|^{N -1} \nabla \u \bigr)+\frac{1}{\delta} |\u|^2 \u = \w - \u \in \mathbb{X},
        \\[1ex]
        -\mathrm{div} \bigl(\alpha(\tilde \eta) M_\u + \kappa^2 \nabla \u +\nu |\nabla \u|^{N -1} \nabla \u \bigr) = \w - \u -\frac{1}{\delta} |\u|^2 \u \in L^{\frac{4}{3}}(\Omega; \R^{M}),
    \end{cases}
        \\[3ex]
    \mbox{in the distributional sense in $ \Omega $.}
    \end{array}
\end{equation*}
Furthermore, applying Green's formula, one can observe that:
\begin{align*}
    \int_\Gamma & \bigl[ \bigl( \alpha(\tilde \eta) M_\u +\kappa^2 \nabla \u +\nu |\nabla \u|^{N -1} \nabla \u \bigr)\mathbf{n}_\Gamma \bigr] \cdot \bm{\varphi} \, d \Gamma
    \\
    & =\int_\Omega \bigl(\alpha(\tilde \eta) M_{\bf u} + \kappa^2 \nabla {\bf u} +\nu |\nu \nabla {\bf u}|^{N -1} \nu \nabla {\bf u} \bigr) : \nabla \bm{\varphi} \, dx  -\int_\Omega \bigl( \w -\u -{\textstyle \frac{1}{\delta}}|\u|^2 \u \bigr) \cdot \bm{\varphi} \, dx
    \\
    & = 0, ~\mbox{for any $ \bm{\varphi} \in C^\infty(\overline{\Omega}; \R^{M}) $.}
\end{align*}
This identity leads to:
\begin{equation}\label{key1.18}
    \bigl[ \bigl( \alpha(\tilde \eta) M_\u +\kappa^2 \nabla \u +\nu |\nabla \u|^{N -1} \nabla \u \bigr)\mathbf{n}_\Gamma \bigr] = 0 ~\mbox{in $ L^2(\Gamma; \R^{M}) $.}
\end{equation}
As a consequence of \eqref{key1.16} and \eqref{key1.18}, we obtain {\it Claim $ \# $2}.
\bigskip

On account of {\it Claims $ \# $1--$ \# $2} and the maximality of {$ \partial_{\bf u} \Phi_{0, \nu, \delta}^{\kappa}(\tilde\eta,\cdot) $} in $ \mathbb{X} \times \mathbb{X} $, we can show the coincidence \eqref{key1.03}, and we conclude this Lemma \ref{Lem02}.
\hfill$ \Box $
}}
\bigskip

Now, we denote by $ \bigl[ \partial_\eta \Phi_{{\varepsilon, \nu, \delta}}^\kappa \times \partial_{\bf u} \Phi_{{\varepsilon, \nu, \delta}}^\kappa \bigr] \subset \mathfrak{X} \times \mathfrak{X} $ a (possibly) set-valued operator, such that:
    \begin{align*}
        & [{ U}, { U}^*] = \bigl[ [\eta, {\bf u}], [\eta^*, {\bf u}^*] \bigr] \in \bigl[ \partial_\eta \Phi_{{\varepsilon, \nu, \delta}}^\kappa \times \partial_{\bf u} \Phi_{{\varepsilon, \nu, \delta}}^\kappa \bigr] \mbox{ in $ \mathfrak{X} \times \mathfrak{X} $,}
    \end{align*}
    iff.
    \begin{align*}
        & \left\{ \hspace{-2ex} \parbox{9.5cm}{
            \vspace{-1ex}
            \begin{itemize}
                \item $ \eta \in D(\partial_\eta \Phi_{{\varepsilon, \nu, \delta}}^\kappa(\cdot, {\bf u})) $, and $ {\bf u} \in D(\partial_{\bf u} \Phi_{{\varepsilon, \nu, \delta}}^\kappa(\eta, \cdot)) $,
                \item $ [\eta^*, {\bf u}^*] \in  \partial_\eta \Phi_{{\varepsilon, \nu, \delta}}^\kappa(\eta, {\bf u}) \times \partial_{\bf u} \Phi_{{\varepsilon, \nu, \delta}}^\kappa(\eta, {\bf u})  $ in $ \mathfrak{X} $.
            \end{itemize}
            \vspace{-1ex}
            } \right.
    \end{align*}
Note that the inclusion $ \partial \Phi_{{\varepsilon, \nu, \delta}}^\kappa \subset \bigl[ \partial_\eta \Phi_{{\varepsilon, \nu, \delta}}^\kappa \times \partial_{\bf u} \Phi_{{\varepsilon, \nu, \delta}}^\kappa \bigr] $ is easily seen, while the converse one is not easy to show.
    \begin{lem}\label{Lem03}
    Let $ \varepsilon \geq 0 $, $ \delta > 0 $, $ \nu > 0 $, and $ \kappa > 0 $ be fixed constants. Then, the following two items hold.
    \begin{description}
        \item[(I)]$ \bigl[ \partial_\eta \Phi_{{\varepsilon, \nu, \delta}}^\kappa \times \partial_{\bf u} \Phi_{{\varepsilon, \nu, \delta}}^\kappa \bigr] $ is semi-monotone in $ \mathfrak{X} \times \mathfrak{X} $, i.e. there exists a positive constant $ R_0 > 0 $, such that $ \bigl[ \partial_\eta \Phi_{{\varepsilon, \nu, \delta}}^\kappa \times \partial_{\bf u} \Phi_{{\varepsilon, \nu, \delta}}^\kappa \bigr] +R_0 \, \mathcal{I} $ is monotone in $ \mathfrak{X} \times \mathfrak{X} $, where $ \mathcal{I} $ is the identity on $ \mathfrak{X} \times \mathfrak{X} $.
        \item[(II)]$ \partial \Phi_{{\varepsilon, \nu, \delta}}^\kappa = \bigl[ \partial_\eta \Phi_{{\varepsilon, \nu, \delta}}^\kappa \times \partial_{\bf u} \Phi_{{\varepsilon, \nu, \delta}}^\kappa \bigr] $ in $ \mathfrak{X} \times \mathfrak{X} $.
    \end{description}
\end{lem}
\paragraph{Proof of Lemma \ref{Lem03} (I).}{
    We set:
    \begin{equation}\label{R0}
        R_0 := 1 +\frac{2}{{\kappa^2}}\|\alpha'\|^2_{L^\infty(\mathbb{R})},
    \end{equation}
    and prove this $ R_0 $ is the required constant.

Let us assume:
\begin{align*}
    & \begin{cases}
        \ds \bigl[ [\eta, {\bf u}], [\eta^*, {\bf u}^*] \bigr] \in \bigl[ \partial_\eta \Phi_{{\varepsilon, \nu, \delta}}^\kappa \times \partial_{\bf u} \Phi_{{\varepsilon, \nu, \delta}}^\kappa \bigr] +R_0 \mathcal{I}
        \\[1ex]
        \ds \bigl[ [\tilde{\eta}, \tilde{{\bf u}}], [\tilde{\eta}^*, \tilde{{\bf u}}^*] \bigr] \in \bigl[ \partial_\eta \Phi_{{\varepsilon, \nu, \delta}}^\kappa \times \partial_{\bf u} \Phi_{{\varepsilon, \nu, \delta}}^\kappa \bigr] +R_0 \mathcal{I}
    \end{cases}
    \mbox{in $ \mathfrak{X} \times \mathfrak{X} $.}
\end{align*}
Then,
\begin{subequations}\label{IA0-00}
\begin{align}
    \bigl( [\eta^*, {\bf u}^*] \, & -[\tilde{\eta}^*, \tilde{{\bf u}}^*], [\eta, {\bf u}] -[\tilde{\eta}, \tilde{{\bf u}}] \bigr)_{\mathfrak{X}}
    \nonumber
    \\[1ex]
    &
    = (\eta^* -\tilde{\eta}^*, \eta -\tilde{\eta})_{H} +({\bf u}^* -\tilde{{\bf u}}^*, {\bf u} -\tilde{{\bf u}})_{\mathbb{X}}
    \nonumber
    \\[1ex]
    & =~ I_1 + I_2 + I_3,
\end{align}
    with
\begin{align}
    I_1 := & \|\nabla(\eta - \tilde{\eta})\|_{[H]^N}^2 +{\kappa^2}\|\nabla ({\bf u} - \tilde{{\bf u}})\|_{[{\mathbb{X}}]^N}^2 +\frac{1}{\delta} \int_\Omega \bigl( |{\bf u}|^2 {\bf u} -|{\bf\tilde{u}}|^2 {\bf\tilde{u}} \bigr) \cdot ({\bf u} -{\bf \tilde{u}}) \, dx
    \nonumber
    \\[1ex]
    & + \int_\Omega \nu^{2}\left( |\nu \nabla {\bf u}|^{N-1}\nabla {\bf u}- |\nu \nabla {\bf \tilde u}|^{N-1}\nabla {\bf \tilde u} \right)\cdot \bigl(\nabla {\bf u} -\nabla \tilde{{\bf u}} \bigr) \, dx
    \nonumber
    \\[1ex]\nonumber
    & +R_0 \bigl( \|\eta - \tilde{\eta}\|_H^2 +\|{\bf u} -\tilde{{\bf u}}\|_{{\mathbb{X}}}^2 \bigr)
    \\[1ex]
    \geq & \|\nabla(\eta - \tilde{\eta})\|_{[H]^N}^2 +{\kappa^2}\|\nabla ({\bf u} - \tilde{{\bf u}})\|_{[{\mathbb{X}}]^N}^2+R_0 \bigl( \|\eta - \tilde{\eta}\|_H^2 +\|{\bf u} -\tilde{{\bf u}}\|_{{\mathbb{X}}}^2 \bigr),
\end{align}
\begin{align}
    I_2:= & \bigl( \alpha'(\eta)f_\varepsilon(\nabla {\bf u}) - \alpha'(\tilde{\eta}) f_\varepsilon(\nabla \tilde{{\bf u}}), \eta-\tilde{\eta} \bigr)_H
    \nonumber
    \\
    & \quad +\frac{1}{{\kappa^2}} \bigl( \alpha(\eta)\alpha'(\eta)-\alpha(\tilde \eta)\alpha'(\tilde \eta), \eta- \tilde \eta \bigr)_H
    \nonumber
    \\
    = & \int_{\Omega} f_\varepsilon(\nabla {\bf u}) \bigl( \alpha'(\eta) - \alpha'(\tilde{\eta}) \bigr)(\eta - \tilde{\eta})\, dx
    \\
    & \quad + \int_{\Omega} \alpha'(\tilde{\eta}) \bigl( f_\varepsilon(\nabla {\bf u})- f_\varepsilon(\nabla \tilde{{\bf u}}) \bigr)(\eta - \tilde{\eta})\, dx
    \nonumber
    \\
    & \ds \quad +\frac{1}{2 {\kappa^2}} \int_{\Omega} \bigl( {\textstyle\frac{d}{d \eta}  [\alpha^2](\eta) -\frac{d}{d \eta} [\alpha^2](\tilde{\eta})}\bigr)(\eta -\tilde{\eta}) \, dx
        \nonumber
        \\
      \geq & - \|\alpha'\|_{L^\infty(\mathbb{R})}\|\eta-\tilde{\eta}\|_H\|\nabla ({\bf u} - \tilde{{\bf u}})\|_{[{\mathbb{X}}]^N}
      \nonumber
\end{align}
by (A3) and H\"older's inequality, and
\begin{align}
    I_3 := &  \int_\Omega \left(\alpha(\eta)[\nabla f_\varepsilon](\nabla {\bf u})-\alpha(\tilde\eta)[\nabla f_\varepsilon](\nabla {\bf \tilde u})\right)\cdot\left(\nabla {\bf u}-\nabla {\bf\tilde u}\right)\, dx\nonumber
    \\ = & \int_\Omega \alpha(\eta)\left([\nabla f_\varepsilon](\nabla{\bf u})-[\nabla f_\varepsilon](\nabla{\bf \tilde u})\right)\cdot\left(\nabla {\bf u}-\nabla {\bf\tilde u}\right)\, dx\nonumber
    \\ + & \int_\Omega (\alpha(\eta)-\alpha(\tilde\eta))[\nabla f_\varepsilon](\nabla {\bf \tilde u})\cdot\left(\nabla {\bf u}-\nabla {\bf\tilde u}\right)\, dx
    \\ \geq &  - \|\alpha'\|_{L^\infty(\mathbb{R})}\|\eta-\tilde{\eta}\|_H \|\nabla ({\bf u} - \tilde{{\bf u}})\|_{[{\mathbb{X}}]^N}\,, \nonumber
\end{align}
by (A3), (A4) and H\"older's inequality.
\end{subequations}
Due to \eqref{R0}, using Young's inequality, the inequalities in \eqref{IA0-00} lead to:
\begin{align*}
    & \bigl( [\eta^*, {\bf u}^*] -[\tilde{\eta}^*, \tilde{{\bf u}}^*], [\eta, {\bf u}] -[\tilde{\eta}, \tilde{{\bf u}}] \bigr)_{\mathfrak{X}} \geq \|\nabla (\eta-\tilde \eta)\|^2_{[H]^N}+\kappa^2\|\nabla({\bf u-\tilde u})\|_{[\mathbb{X}]^N}^2
    \\ & +R_0(\|\eta-\tilde \eta\|^2_H+\|{\bf u-\tilde u}\|^2_{\mathbb{X}})-2\|\alpha'\|_{L^\infty(\mathbb{R})}\|\eta-\tilde{\eta}\|_H \|\nabla ({\bf u} - \tilde{{\bf u}})\|_{[{\mathbb{X}}]^N}
    \\ & \geq \|\nabla (\eta-\tilde \eta)\|^2_{[H]^N}+\kappa^2\|\nabla({\bf u-\tilde u})\|_{[\mathbb{X}]^N}^2+R_0(\|\eta-\tilde \eta\|^2_H+\|{\bf u-\tilde u}\|^2_{\mathbb{X}})
    \\ & -2\frac{\|\alpha'\|^2_{L^\infty(\R)}}{\kappa^2}\|\eta-\tilde \eta\|_H^2-\frac{\kappa^2}{2}\|\nabla({\bf u-\tilde u})\|^2_{[\mathbb{X}]^N}
    \\ & \geq \|\eta - \tilde{\eta}\|_H^2 + \|{\bf u-\tilde u}\|_{\mathbb{X}}^2+\frac{{\kappa^2}}{2}\|\nabla{\bf u}- \tilde{{\bf u}}\|_{[{\mathbb{X}}]^N}^2 \geq 0,
\end{align*}
which implies the (strict) monotonicity of the operator $\bigl[ \partial_\eta \Phi_{{\varepsilon, \nu, \delta}}^\kappa \times \partial_{\bf u} \Phi_{{\varepsilon, \nu, \delta}}^\kappa \bigr] +R_0 \, \mathcal{I}$ in $ \mathfrak{X} \times \mathfrak{X} $.
\qed
}
\paragraph{Proof of Lemma \ref{Lem03} (II).}{
    Let $ R_0 > 0 $ be the constant defined in \eqref{R0}. Then, in the light of the inclusion $ \partial \Phi_{{\varepsilon, \nu, \delta}}^\kappa \subset \bigl[ \partial_\eta \Phi_{{\varepsilon, \nu, \delta}}^\kappa \times \partial_{\bf u} \Phi_{{\varepsilon, \nu, \delta}}^\kappa \bigr]  $ in $ \mathfrak{X} \times \mathfrak{X} $, we can see that:
    \begin{enumerate}
        \item[$\sharp\,6$)]$ \partial \Phi_{{\varepsilon, \nu, \delta}}^\kappa +R_0 \, \mathcal{I} \subset \bigl[ \partial_\eta \Phi_{{\varepsilon, \nu, \delta}}^\kappa \times \partial_{\bf u} \Phi_{{\varepsilon, \nu, \delta}}^\kappa \bigr] +R_0 \, \mathcal{I} $ in $ \mathfrak{X} \times \mathfrak{X} $, i.e. {$ \partial \Phi_{{\varepsilon, \nu, \delta}}^\kappa $ is contained in a monotone graph $ \bigl[ \partial_\eta \Phi_{{\varepsilon, \nu, \delta}}^\kappa \times \partial_{\bf u} \Phi_{{\varepsilon, \nu, \delta}}^\kappa \bigr] +R_0 \, \mathcal{I} $ in $ \mathfrak{X} \times \mathfrak{X} $.}
    \end{enumerate}
    Also, invoking \cite[Theorem 2.10]{MR2582280} and \cite[Corollary 2.11]{MR0348562}, we also have:
    \begin{enumerate}
        \item[$\sharp\,7$)]$ \partial \bigl( \Phi_{{\varepsilon, \nu, \delta}}^\kappa +\frac{R_0}{2} \|\cdot\|_{\mathfrak{X}}^2 \bigr) = \partial \Phi_{{\varepsilon, \nu, \delta}}^\kappa +R_0 \, \mathcal{I} $ in $ \mathfrak{X} \times \mathfrak{X} $, and hence, {$ \bigl[ \partial_\eta \Phi_{{\varepsilon, \nu, \delta}}^\kappa \times \partial_{\bf u} \Phi_{{\varepsilon, \nu, \delta}}^\kappa \bigr] +R_0 \, \mathcal{I} $ coincides with the maximal monotone graph $ \partial \bigl( \Phi_{{\varepsilon, \nu, \delta}}^\kappa +\frac{R_0}{2} \|\cdot\|_{\mathfrak{X}}^2 \bigr) $ in $ \mathfrak{X} \times \mathfrak{X} $.}
    \end{enumerate}

    Then, item (II) is a direct consequence of the above $\sharp\,6$) and $\sharp\,7$). \qed
}

\bigskip

Based on the above Lemmas, Theorems \ref{Th.SolvApKWC} and \ref{Th.KWC-equiv-CP} are proved as follows.

\paragraph{Proof of Theorem \ref{Th.KWC-equiv-CP}.}{
    On account of Lemmas \ref{Lem01} and \ref{Lem02}, one can observe that the system (P)$_{{\varepsilon, \nu, \delta}}^\kappa$ is equivalent to the following Cauchy problem, denoted by ($\widetilde{\mathrm{CP}}$)$_{{\varepsilon, \nu, \delta}}^\kappa$.
\vspace{-2ex}

\paragraph{\boldmath Cauchy problem ($\widetilde{\mathrm{CP}}$)$_{{\varepsilon, \nu, \delta}}^\kappa$:}{
    \begin{align*}
        & \begin{cases}
            { U}'(t) +\bigl[ \partial_\eta \Phi_{{\varepsilon, \nu, \delta}}^\kappa \times \partial_{\bf u} \Phi_{{\varepsilon, \nu, \delta}}^\kappa \bigr]({ U}(t)) +{\mathcal{G}_{\delta}^{\kappa}}({ U}(t)) \ni 0 \mbox{ in $ \mathfrak{X} $, ~ $ {t \in (0, T)} $,}
            \\[1ex]
            { U}(0) = { U}_0 \mbox{ in $ \mathfrak{X} $;}
        \end{cases}
    \end{align*}

}
\noindent
Additionally, in the light of Lemma \ref{Lem03}, the above Cauchy problem ($\widetilde{\mathrm{CP}}$)$_{{\varepsilon, \nu, \delta}}^\kappa$ is equivalent to the Cauchy problem (CP)$_{{\varepsilon, \nu, \delta}}^\kappa$.

Thus, the problems (P)$_{{\varepsilon, \nu, \delta}}^\kappa$, (CP)$_{{\varepsilon, \nu, \delta}}^\kappa$, and ($\widetilde{\mathrm{CP}}$)$_{{\varepsilon, \nu, \delta}}^\kappa$ are equivalent to each other, and this finishes the proof of Theorem \ref{Th.KWC-equiv-CP}. \qed
}

\paragraph{Proof of Theorem \ref{Th.SolvApKWC}.}{
    {First, we show item (I). 
From hypothesis (A6) with $ {\bf u}_0 \in L^4(\Omega; \R^{M}) $ and $ \nu {\bf u}_0 \in W^{1, N +1}(\Omega; \R^{M}) $ we get
~${ U}_0 = [\eta_0, {\bf u}_0] \in D(\Phi_{{\varepsilon, \nu, \delta}}^\kappa)$.
    Hence, by applying the general theories of evolution equations, e.g., \cite[{Theorem~4.1, p.~158}]{MR2582280}, \cite[Theorem~3.6 {and Proposition~3.2}]{MR0348562}, \cite[Section~2]{MR1720614} and \cite[Theorem~1.1.2]{Kenmochi81}, we immediately have the existence and uniqueness of solution $ { U} = [\eta, {\bf u}] \in {L^2(0, T; \mathfrak{X})} $ to (CP)$_{{\varepsilon, \nu, \delta}}^\kappa$ (and therefore for (P)$_{{\varepsilon, \nu, \delta}}^\kappa$, by Theorem \ref{Th.KWC-equiv-CP}), such that:
\begin{equation}\label{CPreg01}
    { U} \in {W^{1, 2}(0, T; \mathfrak{X}) \cap L^\infty(0, T; \mathfrak{W}),} \mbox{ and } \Phi_{{\varepsilon, \nu, \delta}}^\kappa ({U}) \in {L^{\infty}(0, T).}
\end{equation}
{\eqref{defAppPhi} and \eqref{CPreg01} imply the regularity {stated} in \eqref{CPreg}. Additionally,} for any constant $ 0 < T < \infty $, there exists a positive constant {$ C_1 $, independent of $ U_0 $, such that:} 
\begin{equation}\label{MT1.2}
\|{ U}'\|_{L^2(0,T;\mathfrak{X})}^2 + \sup_{t \in (0,T)} \Phi_{{\varepsilon, \nu, \delta}}^\kappa({ U}(t)) \leq C_1 \left(1 + \|{ U}_0\|_{\mathfrak{X}}^2 +\Phi_{{\varepsilon, \nu, \delta}}^\kappa({ U}_0) \right).
\end{equation}

Next, we verify item (II). Let us fix any bounded open interval $ I \subset {(0, T)} $. Then, taking into account \eqref{eps-dlt-nu}--\eqref{cv01}, \eqref{defAppPhi}, and \eqref{MT1.2}, one can observe that $ \{ \Phi_{\varepsilon_n, \nu_n, \delta_n}^{\kappa_n} ({ U}_{0, n}) \}_{n = 1}^\infty $ is bounded, and furthermore,
\begin{align*}
    & \begin{cases}
        \ds\{ { U}_n \}_{n = 1}^\infty = \{ [\eta_n, {\bf u}_n] \}_{n = 1}^\infty  \mbox{ is bounded in $ W^{1,2}(I; \mathfrak{X}) \cap L^{\infty}(I; \mathfrak{W}) $,}
\\[0.5ex]
        \ds\{ \nu_{n} \nabla {\bf u}_{n} \}_{n = 1}^\infty \mbox{ is bounded in $ L^{\infty}(I; L^{N +1}(\Omega; \R^{MN})) $.}
    \end{cases}
\end{align*}
Therefore, applying the Poincar\'{e}--Wirtinger inequality, and general compactness theories, such as Ascoli's theorem (cf. \cite[Corollary 4]{MR0916688}), we find subsequences of \linebreak $ \{ [\kappa_n, \varepsilon_n, \nu_n, \delta_n] \}_{n = 1}^\infty $ and $ \{ { U}_n \}_{n = 1}^\infty = \{ [\eta_n, {\bf u}_n] \}_{n = 1}^\infty $ (not relabeled), such that:
\begin{subequations}\label{MT2.1}
\begin{align}\label{MT2.1a}
    { U}_n \to { U} ~~& \mbox{in $ C(\overline{I}; \mathfrak{X}) $, weakly in $ W^{1,2}(I; \mathfrak{X}) $, and weakly-$*$ in $ L^{\infty}(I; \mathfrak{W}) $,}
\end{align}
together with
\begin{align}\label{MT2.1b}
    & \alpha(\eta_n) \to \alpha(\eta) \mbox{ in $ L^{2}(I; H) $,}
\end{align}
and
\begin{align}\label{MT2.1c}
    & \nu_n {\bf u}_n \to \nu {\bf u} \mbox{ weakly in $ L^{N +1}(I; W^{1, N +1}(\Omega; \R^{M})) $, as $ n \to \infty $.}
\end{align}
\end{subequations}
Also, by  Lemma \ref{Lem04} and (Fact\,2),  we have that a sequence $ \{ \widehat{\Phi}_n^I \}_{n = 1}^\infty $ of proper, l.s.c., and convex functions on $ L^2(I; \mathfrak{X}) $, defined as:
        \begin{equation*}
            \tilde{{ U}} \in L^2(I; \mathfrak{X}) \mapsto \widehat{\Phi}_n^I(\tilde{{ U}}) := \left\{ \begin{array}{ll}
                \multicolumn{2}{l}{\ds \int_I \Phi_{\varepsilon_n, \nu_n, \delta_n}^{\kappa_n}(\tilde{{ U}}(t)) \, dt,}
                    \\[1ex]
                    & \mbox{ if $ \Phi_{\varepsilon_n, \nu_n, \delta_n}^{\kappa_n}(\tilde{{ U}}) \in L^1(I) $,}
                    \\[2.5ex]
                    \infty, & \mbox{ otherwise,}
                \end{array} \right. \mbox{for $ n = 1, 2, 3, \dots $;}
        \end{equation*}
        converges to a proper, l.s.c., and convex function $ \widehat{\Phi}^I $ on $ L^2(I; \mathfrak{X}) $, defined as:
        \begin{equation*}
            \tilde{{ U}} \in L^2(I; \mathfrak{X}) \mapsto \widehat{\Phi}^I(\tilde{{ U}}) := \left\{ \begin{array}{ll}
                \multicolumn{2}{l}{\ds \int_I \Phi_{{\varepsilon, \nu, \delta}}^\kappa(\tilde{{ U}}(t)) \, dt, \mbox{ if $ \Phi_{{\varepsilon, \nu, \delta}}^\kappa(\tilde{{ U}}) \in L^1(I) $,}}
                    \\[2ex]
                    \infty, & \mbox{ otherwise;}
                \end{array} \right.
        \end{equation*}
        on $ L^2(I; \mathfrak{X}) $, in the sense of Mosco, as $ n \to \infty $.

From \eqref{MT2.1} and (Fact\,1), it is seen that
\begin{equation}\label{MT2.3}
    [{U}, -{U}' - {\mathcal{G}_\delta^\kappa}({U})] \in \partial \widehat{\Phi}^I \mbox{ in $ L^2(I;\mathfrak{X}) \times L^2(I;\mathfrak{X}) $,}
\end{equation}
and
\begin{equation}\label{MT2.4}
\widehat{\Phi}_{n}^I({ U}_n) \to \widehat{\Phi}^I({ U}) \mbox{ as $ n \to \infty $.}
\end{equation}
Note that \eqref{cv00}, \eqref{MT2.3}, and (Fact\,1) enable us to say that $ { U} = [\eta, {\bf u}] $ is a solution to the Cauchy problem (CP)$_{{\varepsilon, \nu, \delta}}^\kappa$. So, {due to uniqueness} of solutions, the convergences in \eqref{MT2.1} hold for the whole sequence and not only for subsequences.
\bigskip

In the meantime, from (A3), (A4), and \eqref{MT2.1}, it is inferred that:
\begin{subequations}\label{MT2.6}
\begin{align}\label{MT2.6a}
    & \lim_{n \to \infty} \frac{1}{\kappa_n^2} \int_I \int_\Omega \alpha(\eta_n)^2  \, dx dt =  \frac{1}{\kappa^2} \int_I \int_\Omega \alpha(\eta)^2 \, dx dt,
\end{align}
\begin{align}\label{MT2.6b}
    \varliminf_{n \to \infty} & \frac{1}{2} \int_I \int_\Omega |\nabla \eta_n|^2 \, dx dt \geq \frac{1}{2} \int_I \int_\Omega |\nabla \eta|^2 \, dx dt,
\end{align}
\begin{align}\label{MT2.6c}
    & \varliminf_{n \to \infty} \frac{1}{\delta_n} \int_I \int_\Omega |{\bf u}_n|^{4} \, dx dt \geq \frac{1}{\delta} \int_I \int_\Omega |{\bf u}|^{4} \, dx dt,
\end{align}
\begin{align}\label{MT2.6d}
    \varliminf_{n \to \infty} & \frac{1}{N +1} \int_I \int_\Omega |\nu_n \nabla {\bf u}_n|^{N +1} \, dx dt \geq \frac{1}{N +1} \int_I \int_\Omega |\nu \nabla {\bf u}|^{N +1} \, dx dt,
\end{align}
\begin{align}\label{MT2.6e}
    \varliminf_{n \to \infty} \frac{\kappa_n^2}{2} \int_I \int_\Omega f_{\varepsilon_n}({\nabla} & {\bf u}_n)^2 \, dx dt  
\geq \frac{\kappa^2}{2} \int_I \int_\Omega f_\varepsilon(\nabla {\bf u})^2 \, dx dt,
\end{align}
and
\begin{align}\label{MT2.6f}
    \varliminf_{n \to \infty} &  \int_I \int_\Omega \alpha(\eta_n) f_{\varepsilon_n}(\nabla {\bf u}_n) \, dx dt
    \nonumber
    \\
    & \geq \varliminf_{n \to \infty} \int_I \int_\Omega \alpha(\eta) f_\varepsilon(\nabla {\bf u}_n) \, dx dt
    \nonumber
    \\
    & \qquad -\varlimsup_{n\to\infty}\|\alpha(\eta_n) -\alpha(\eta)\|_{L^2(I; H)} \sup_{n \in \N} \left( \int_I \int_\Omega \bigl( \varepsilon_n^2 +|\nabla {\bf u}_n|^2 \bigr) \, dx dt \right)^{\frac{1}{2}}
    \nonumber
    \\
    & \qquad -\|\alpha(\eta)\|_H \varlimsup_{n \to \infty} \left( \int_I \int_\Omega |\varepsilon_n -\varepsilon|^2 \, dx dt \right)^{\frac{1}{2}}
    \nonumber
    \\
    & \geq \int_I \int_\Omega \alpha(\eta) f_\varepsilon(\nabla {\bf u}) \, dx dt\,.
\end{align}
\end{subequations}
From \eqref{MT2.4} and \eqref{MT2.6}, it follows that:
\begin{align}\label{MT2.7}
    & \begin{cases}
        \ds \lim_{n \to \infty} \bigl\| \nabla \eta_n \bigr\|_{L^2(I; [H]^N)}^2 = \bigl\| \nabla \eta \bigr\|_{L^2(I; [H]^N)}^2,
        \\[1ex]
        \ds \lim_{n \to \infty} \bigl\| \nabla {\bf u}_n \bigr\|_{L^2(I; L^2(\Omega; \R^{MN}))}^2 = \bigl\| \nabla {\bf u} \bigr\|_{L^2(I; L^2(\Omega; \R^{MN}))}^2,
        \\[1ex]
        \ds \lim_{n \to \infty} \bigl\| \nu \nabla {\bf u}_n \bigr\|_{L^{N +1}(I; L^{N +1}(\Omega; \R^{MN}))}^{N +1} = \bigl\| \nu \nabla {\bf u} \bigr\|_{L^{N +1}(I; L^{N +1}(\Omega; \R^{MN}))}^{N +1}.
    \end{cases}
\end{align}

Now, taking into account \eqref{MT2.1} and \eqref{MT2.7}, and applying the uniform convexity of $ L^2 $-based topologies, one can observe that:
\begin{equation}\label{MT2.7-1}
    \begin{cases}
        \eta_n \to \eta ~\mbox{in $ L^2(I; V) $,}~ {\bf u}_n \to {\bf u} ~\mbox{in $ L^2(I; \mathbb{X}) $,}
        \\[1ex]
        \nu_n \nabla {\bf u}_n \to \nu \nabla {\bf u} ~\mbox{in $ L^{N +1}(I; L^{N +1}(\Omega; \R^{MN})) $,}
    \end{cases}
    \mbox{as $ n \to \infty $.}
\end{equation}

Thus, convergences \eqref{cv00} and \eqref{cv01} are concluded as direct consequences of \eqref{MT2.1} and \eqref{MT2.7-1}. \hfill \qed
\medskip

    }
}

\section{Proof of Main Theorem.}\label{sec:sphere}

In this Section, we show the proof of Main Theorem. To see this, we take the limit in (P)$_{{\varepsilon, \nu, \delta}}^\kappa$  as $\delta$, $\nu$, $\varepsilon\to 0^+$, respectively. 

At first, we derive an energy inequality for {$ \mathcal{F}_{\varepsilon, \nu, \delta}^{\kappa} $} and a priori estimates for the approximate solutions.

\subsection{Energy inequality and a priori estimates.}
\begin{lem}
Let ${ U}_{0} = [\eta_{0}, {\bf u}_{0}] \in 
    \mathfrak{W}$, with ${\bf u}_0\in L^4(\Omega;\R^M)$, $\nu {\bf u}_0\in L^{N+1}(\Omega;\R^M)$ and $\Ua := [\etaa, \ua]$ be a solution to \rm{(P)$_{{\varepsilon, \nu, \delta}}^\kappa$}. 
    Then, for any $T>0$, $\Ua$ satisfies the following energy inequality:
\begin{equation}\label{ene-inq-1}
    \mathcal{F}_{{\varepsilon, \nu, \delta}}^\kappa({ U}_{{\varepsilon, \nu, \delta}}(s)) + \int_{0}^{s} \|\partial_{t}\Ua(t)\|^{2}_{\mathfrak{X}} \, dt \le \mathcal{F}_{{\varepsilon, \nu, \delta}}^\kappa({ U}_{0})\ \ \mbox{ {for all} } s \in [0,T].
\end{equation}

\noindent
Moreover, it follows that
\begin{equation}\label{app.sol-1}
\left\{ \begin{array}{l}
\ds \Ua \in W^{1,2}({0, T} ; \mathfrak{X}) \cap L^{\infty}({0, T} ; \mathfrak{W}), \vspace{3mm}\\
    \ds \kappa^{2} \nabla \ua \in L^{\infty}({0, T} ; L^{2}(\Omega ; \R^{MN})
),\ \ \nu \nabla \ua \in L^{\infty}({0, T}; {L^{N+1}}(\Omega, \R^{MN})),
\end{array} \right.
\end{equation}
and
\begin{equation}\label{app.sol-2}
    \begin{cases}
        \delta^{- \frac{1}{2}} \bigl\||\ua|^{2} - 1\bigr\|_{L^{\infty}(0,T; H)} \le {\mathcal{F}_{{\varepsilon, \nu, \delta}}^\kappa({ U}_{0})},
        \\[1ex]
        \bigl\|\alpha(\etaa)[\nabla f_{\varepsilon}](\nabla \ua)\bigr\|_{L^{\infty}({0, T}; L^{\infty}(\Omega; \R^{MN}))} \le {\mathcal{F}_{{\varepsilon, \nu, \delta}}^\kappa({ U}_{0})}.
    \end{cases}
\end{equation}
\end{lem}
\paragraph{Proof.}{
Multiplying both sides of the equation of $\ua$ in (P)$_{{\varepsilon, \nu, \delta}}^\kappa$ by $\partial_{t}\ua$, we have
\begin{align}\label{lem5-1}
\|\partial_{t}\ua\|_{
{\mathbb{X}}}^{2} & + \frac{d}{dt} \Psi_{{\varepsilon, \nu, \delta}}^\kappa(\etaa, \ua)  - \int_{\Omega} \alpha'(\etaa) \partial_{t}\etaa f_{\varepsilon}(\nabla \ua) dx = 0.
\end{align}
Also, multiplying both side of the equation of $\etaa$ in (P)$_{{\varepsilon, \nu, \delta}}^\kappa$ by $\partial_{t}\etaa$, it follows that
\begin{align}\label{lem5-2}
\|\partial_{t}\etaa\|_{H}^{2} + \frac{d}{dt} \Psi_{0}(\etaa) + \int_{\Omega} \alpha'(\etaa) \partial_{t}\etaa f_{\varepsilon}(\nabla \ua) dx = 0.
\end{align}
Adding up \eqref{lem5-1}--\eqref{lem5-2}, and integrating from $0$ to $s$, we get \eqref{ene-inq-1}.
\bigskip

The regularity and the a priori estimates are given immediately by Theorem \ref{Th.SolvApKWC} and (\ref{ene-inq-1}).


\begin{lem}
If $|{\bf u}_{0}| \le 1$ a.e. in $\Omega$, then solutions $\ua$ to (P)$_{{\varepsilon, \nu, \delta}}^\kappa$ satisfy $|\ua| \le 1$ a.e. in $\Omega$.
\end{lem}
\paragraph{Proof.}{
    {We define:}
\begin{equation*}
\chi(z) := \frac{(z-1)_{+}}{z} = \left\{ \begin{array}{ll}
0 \ \ & \mbox{ for } 0 \le z \le 1, \vspace{3mm}\\
\ds \frac{z-1}{z}\ \ & \mbox{ for } z > 1.
\end{array}\right.
\end{equation*}
Then, $\chi$ is a nonnegative monotone increasing function on $[0, \infty)$. Multiplying both sides of
the equation of $\ua$ in (P)$_{{\varepsilon, \nu, \delta}}^\kappa$ by $\ua \chi(|\ua|)$, we have
\begin{align*}
& \frac{1}{2} \frac{d}{dt} \int_{\{ |\ua| > 1\}} 
(|\ua| - 1)^{2} dx + \int_{\{ |\ua| > 1\}} \varpi_{\delta}(\ua) \cdot \ua \chi(|\ua|) dx \nonumber \\[1ex]
    = & - \int_{\{ |\ua| > 1\}} (\alpha(\etaa) [\nabla f_{\varepsilon}](\nabla \ua) \\[1ex]
& \hspace{30mm} + {\kappa^2} \nabla \ua + {\nu | \nu \nabla \ua|^{N-1} \nu \nabla \ua}) : \nabla \left( \ua 
\chi(|\ua|) \right) dx \\[1ex]
= & - \int_{\{ |\ua| > 1\}} (\alpha(\etaa) [\nabla f_{\varepsilon}](\nabla \ua) : \nabla \ua  + {\kappa^2} |\nabla \ua|^{2} + |\nu \nabla \ua|^{N+1}) \chi(|\ua|) dx \\[1ex]
    & - \int_{\{ |\ua| > 1\}} \left( \frac{\alpha(\etaa)}{\sqrt{\varepsilon^{2} + |\nabla \ua|^{2}}} + {\kappa^2} + \nu^2|\nu \nabla \ua|^{N-1}\right)  \frac{|\nabla |\ua|^{2} |^{2}}{4|\ua|^{3}} dx \le 0.
\end{align*}
Here,
\begin{align*}
    \nabla \ua : \,&  \frac{\ua \nabla |\ua|}{|\ua|^{2}} = \frac{\nabla \ua : (\ua ((\nabla \ua)^{t} \ua))}{|\ua|^{3}}
    \\[1ex]
    & = \frac{\bigl|(\nabla \ua)^{t}\ua \bigr|^{2}}{|\ua|^{3}} = \frac{\bigl|\nabla |\ua|^{2} \bigr|^{2}}{4|\ua|^{3}}.
\end{align*}
Therefore,
\begin{equation*}
\int_{\{|\ua(t)| > 1\}} (|\ua(t)| - 1)^{2} dx \le \int_{\{|\ua(0)| > 1\}} (|\ua(0)| - 1)^{2} dx = 0.
\end{equation*}
Hence, $|\ua| \le 1$ a.e. in $\Omega$. \qed
}

\bigskip
At the end of this subsection, we introduce some compactness results which can be proved as \cite[Theorem 2.1]{MR1254812} and \cite[Lemma 9]{MR1970719} by using the fact that the operator 
{${\rm div}\left(\alpha(\etaa) [\nabla f_{\varepsilon}](\nabla \cdot) +{\kappa^2} (\nabla \cdot) +\nu^{N +1} |\nabla \cdot|^{N} (\nabla \cdot) \right)$}
is uniformly elliptic.

\begin{lem}[{cf.\cite[Lemma 2.2]{MR2466164}}]\label{lem-cpt-1}
    Let $\varepsilon > 0$, $\nu > 0$ and $\kappa > 0$ be fixed. Let $\{{\bf w}_{\varepsilon, \nu , \delta}\}_{\delta>0}$ be bounded in $W^{1,2}({0, T}; {\mathbb{X}}) \cap L^{\infty}({0, T}; {\mathbb{W}})$, $\{\eta_{\varepsilon, \nu , \delta}\}_{\delta>0}$ be bounded in $W^{1,2}({0, T}; {H}) \cap L^{\infty}({0, T}; {V})$, and $\{{\bf f}_{\varepsilon, \nu, \delta} \}_{\delta>0}$ be bounded in $L^{1}({0, T}; L^{1}(\Omega;\R^{M}))$ uniformly in $\delta$, respectively. Moreover, $0 \le \eta_{\varepsilon,\nu, \delta} \le 1$ a.e. in $\Omega_T$ and ${\bf w}_{\varepsilon, \nu, \delta}$ satisfies the equation: for $\delta > 0$,
\begin{align*}
    \partial_{t} {\bf w}_{\varepsilon, \nu, \delta} - \mathrm{div} \bigl(\alpha(\eta_{\varepsilon, \nu, \delta}) [\nabla f_{\varepsilon}](\nabla  \,& {\bf w}_{\varepsilon, \nu, \delta}) + {\kappa^2} \nabla {\bf w}_{\varepsilon, \nu, \delta} + {\nu |\nu \nabla {\bf w}_{\varepsilon, \nu, \delta}|^{N-1} \nu \nabla {\bf w}_{\varepsilon, \nu, \delta}} \bigr)
    \\
    & = {\bf f}_{\varepsilon, \nu, \delta} ~~ \mbox{a.e. in } \Omega_T,
\end{align*}
in the sense of distributions. Then $\{{\bf w}_{\varepsilon, \nu, \delta}\}_{\delta>0}$ is precompact in $L^{q}({0, T} ; W^{1, q}(\Omega ; \R^{M}))$ for all $1 \le q < 2$.

\end{lem}

\begin{lem}[{cf.\cite[Lemma 2.2]{MR2466164}}]\label{lem-cpt-2}
    Let $\varepsilon > 0$ and $\kappa > 0$ be fixed. Let $\{{\bf w}_{\varepsilon, \nu}\}_{\nu>0}$ be bounded in $W^{1,2}({0, T}; {\mathbb{X}}) \cap L^{\infty}({0, T}; {\mathbb{W}})$, $\{\eta_{\varepsilon, \nu}\}_{\nu>0}$ be bounded in $W^{1,2}({0, T}; {H}) \cap L^{\infty}({0, T}; {V})$, and $\{{\bf f}_{\varepsilon, \nu} \}_{\nu>0}$ be bounded in $L^{1}({0, T}; L^{1}(\Omega;\R^{M}))$ uniformly in $\nu$, respectively. Moreover, $0 \le \eta_{\varepsilon,\nu} \le 1$ a.e. in $\Omega_T$ and ${\bf w}_{\varepsilon, \nu}$ satisfies the equation: for $\nu > 0$,
\begin{equation*}
    \partial_{t} {\bf w}_{\varepsilon, \nu} - \mathrm{div} (\alpha(\eta_{\varepsilon, \nu}) [\nabla f_{\varepsilon}](\nabla {\bf w}_{\varepsilon, \nu}) + {\kappa^2} \nabla {\bf w}_{\varepsilon, \nu} + {\nu |\nu \nabla {\bf w}_{\varepsilon, \nu}|^{N-1} \nu \nabla {\bf w}_{\varepsilon, \nu}} ) = {\bf f}_{\varepsilon, \nu}~~ \mbox{a.e. in } \Omega_T,
\end{equation*}
in the sense of distributions. Then $\{{\bf w}_{\varepsilon, \nu}\}_{\nu>0}$ is precompact in $L^{q}({0, T} ; W^{1, q}(\Omega ; \R^{M}))$ for all $1 \le q < 2$.

\end{lem}

\begin{lem}[{cf.\cite[Lemma 2.3]{MR2466164}}]\label{lem-cpt-3}
    Let $\kappa > 0$ be fixed. Let $\{{\bf w}_{\varepsilon}\}_{\varepsilon>0}$ be bounded in \linebreak $W^{1,2}({0, T}; {\mathbb{X}}) \cap L^{\infty}({0, T}; {\mathbb{W}})$, $\{\eta_{\varepsilon}\}_{\varepsilon>0}$ be bounded in $W^{1,2}({0, T}; {H}) $ $ \cap L^{\infty}({0, T}; {V})$, and $\{{\bf f}_{\varepsilon} \}_{\varepsilon>0}$ be bounded in $L^{1}({0, T}; L^{1}(\Omega;\R^{M}))$ uniformly in $\varepsilon$, respectively. Moreover, $0 \le \eta_{\varepsilon} \le 1$ a.e. in $\Omega_T$ and ${\bf w}_{\varepsilon}$ satisfies the equation: for $\varepsilon>0$,
\begin{equation*}
    \partial_{t} {\bf w}_{\varepsilon} - \mathrm{div} (\alpha(\eta_{\varepsilon}) [\nabla f_{\varepsilon}](\nabla {\bf w}_{\varepsilon}) + {\kappa^2} \nabla {\bf w}_{\varepsilon} ) = {\bf f}_{\varepsilon} ~~ \mbox{a.e. in } \Omega_T,
\end{equation*}
in the sense of distributions. Then $\{{\bf w}_{\varepsilon}\}_{\varepsilon>0}$ is precompact in
$L^{1}({0, T} ; W^{1, 1}(\Omega ; \R^{M}))$.

\end{lem}

\subsection{The limit as $\delta \to 0$.}

In this subsection we study the limit problem, as $\delta\to 0$, in (P)$_{{\varepsilon, \nu, \delta}}^\kappa$, assuming $0<\varepsilon<\varepsilon_0$; i.e. we solve the following problem. Note that the dependence on $ \kappa $ is not anymore needed or used. Therefore, we remove it.

\paragraph{\boldmath Problem (P)$_{\varepsilon, \nu}$:}
{
    \begin{align*}
        & \begin{cases}
            \partial_t \eta_{\varepsilon, \nu} -\mathit{\Delta} \eta_{\varepsilon, \nu} +g(\eta_{\varepsilon, \nu}) +\alpha'(\eta_{\varepsilon, \nu}) f_\varepsilon(\nabla {\bf u}_{\varepsilon, \nu}) = 0 \mbox{ in } \Omega_T,
            \\[0.5ex]
            \nabla \eta_{\varepsilon, \nu} \cdot {{\bf n}_\Gamma} = 0 \mbox{ on $ \Gamma_T $,}
            \\[0.5ex]
            \eta_{\varepsilon, \nu}(0, x) = \eta_0(x), ~ x \in \Omega;
        \end{cases}
        \\[1ex]
        & \begin{cases}
            \partial_t {\bf u}_{\varepsilon, \nu} -\mathrm{div} \bigl( \alpha(\eta_{\varepsilon, \nu}) {\,[\nabla f_\varepsilon]}(\nabla {\bf u}_{\varepsilon, \nu}) +{\kappa^2} \nabla {\bf u}_{\varepsilon, \nu} + {\nu |\nu \nabla {\bf u}_{\varepsilon, \nu}|^{N-1} \nu \nabla {\bf u}_{\varepsilon, \nu}} \bigr) 
            = \mu_{\varepsilon, \nu} {\bf u}_{\varepsilon, \nu} \mbox{ in $ \Omega_T $,}
            \\[0.5ex]
            \bigl( \alpha(\eta_{\varepsilon, \nu}) [\nabla f_\varepsilon](\nabla {\bf u}_{\varepsilon, \nu}) +{\kappa^2} \nabla {\bf u}_{\varepsilon, \nu} + {\nu |\nu \nabla {\bf u}_{\varepsilon, \nu}|^{N-1} \nu \nabla {\bf u}_{\varepsilon, \nu}} \bigr) {{\bf n}_\Gamma} = 0 \mbox{ on $ \Gamma_T $,}
            \\[0.5ex]
            {\bf u}_{\varepsilon, \nu}(0, x) = {\bf u}_0(x), ~ x \in \Omega;
        \end{cases}
    \end{align*}
}
together with:
\begin{equation*}
    \mu_{\varepsilon, \nu} := (\alpha(\eta_{\varepsilon, \nu})[\nabla f_{\varepsilon}](\nabla {\bf u}_{\varepsilon, \nu}) + {\kappa^2}\nabla {\bf u}_{\varepsilon, \nu} + { \nu |\nu \nabla {\bf u}_{\varepsilon, \nu}|^{N-1} \nu \nabla {\bf u}_{\varepsilon, \nu}}) : \nabla {\bf u}_{\varepsilon, \nu}, ~\mbox{a.e. in $ \Omega_T $.}
\end{equation*}

\begin{thm}
\label{thm-1}
Let ${ U}_{0}=[\eta_{0}, {\bf u}_{0}] \in \mathfrak{W}$ with  $|{\bf u}_{0}| = 1$ in $\Omega$. Then, there exists ${ U}_{\varepsilon, \nu} = [\eta_{\varepsilon, \nu}, {\bf u}_{\varepsilon, \nu}] \in C({[0,T]};\mathfrak{X})$ such that ${ U}_{\varepsilon, \nu}$ satisfies (P)$_{\varepsilon, \nu}$ in the sense of distributions and
\begin{equation}\label{app.sol-3}
\left\{ \begin{array}{l}
\ds { U}_{\varepsilon, \nu} \in W^{1,2}({0, T} ; \mathfrak{X}) \cap L^{\infty}({0, T} ; \mathfrak{W}), \vspace{3mm}\\
    {\kappa^2} \nabla {\bf u}_{\varepsilon, \nu} \in L^{\infty}({0, T} ; L^{2}(\Omega;\R^{MN})), \ \ \nu \nabla {\bf u}_{\varepsilon, \nu} \in L^{\infty}({0, T}; L^{N+1}(\Omega ; \R^{MN})), \vspace{3mm}\\
\ds 0 \le \eta_{\varepsilon, \nu} \le 1\ \ \mbox{ a.e. in } \Omega_{T},\ \ |{\bf u}_{\varepsilon, \nu}| = 1\ \ \mbox{  a.e. in } \Omega_{T},
\end{array}\right.
\end{equation}
and
\begin{equation}\label{ene-inq-2}
    \mathcal{F}_{\varepsilon, \nu}({ U}_{\varepsilon, \nu}(s)) + \int_{0}^{s} \| \partial_{t}{ U}_{\varepsilon, \nu}(t) \|^{2}_{\mathfrak{X}} dt \le \mathcal{F}_{\varepsilon, \nu}({ U}_{0})\leq \mathcal F_{\varepsilon_0,\nu}(U_0)\ \ \mbox{ {for all} } s \in [0,T],
\end{equation}
\begin{align}
& \bigl\| \alpha(\eta_{\varepsilon, \nu})[\nabla f_{\varepsilon}](\nabla {\bf u}_{\varepsilon, \nu}) \bigr\|_{L^{\infty}(
\Omega_{T}; \R^{MN})} \le C, \label{app.sol-4} \\[1ex]
    & \left\|\mbox{div} \left( \begin{array}{l}
        \alpha(\eta_{\varepsilon, \nu})[\nabla f_{\varepsilon}](\nabla {\bf u}_{\varepsilon, \nu}) + {\kappa^2}\nabla {\bf u}_{\varepsilon, \nu} 
        \\
        \hspace{14ex}+ {\nu |\nu \nabla {\bf u}_{\varepsilon, \nu}|^{N-1} \nu \nabla {\bf u}_{\varepsilon, \nu}} 
    \end{array} \right) \right\|_{L^{2}({0, T};L^{1}(\Omega;\R^{M}))} < C, \label{app.sol-5}\\[1ex]
    & \left\| \mbox{div} \left(\left( \begin{array}{l} 
        \alpha(\eta_{\varepsilon, \nu})[\nabla f_{\varepsilon}](\nabla {\bf u}_{\varepsilon, \nu}) + {\kappa^2}\nabla {\bf u}_{\varepsilon, \nu}
        \\
        \hspace{12ex}+{\nu |\nu \nabla {\bf u}_{\varepsilon, \nu}|^{N-1} \nu \nabla {\bf u}_{\varepsilon, \nu}} 
    \end{array} \right) \wedge {\bf u}_{\varepsilon, \nu} \right) \right\|_{{L^{2}(\Omega_{T};\Lambda_{2}(\R^{M}))}} < {C,} \label{app.sol-6}
\end{align}
where the constant $C>0$ is independent of $\varepsilon$.
\end{thm}
\paragraph{Proof.}{
We observe that, since $|{\bf u}_{0}| = 1$, $\Pi_{\delta}({\bf u}_{0}) = 0$. Hence,
\begin{equation*}
    \mathcal F_{{\varepsilon, \nu, \delta}}^{\kappa}({ U}_{0}) =\mathcal F_{\varepsilon, \nu}({ U}_{0})< \mathcal F_{\varepsilon_0, \nu}({ U}_{0})=:C <+\infty.
\end{equation*}
Therefore, having in mind the energy inequality (\ref{ene-inq-1}), the uniform estimates (\ref{app.sol-1})--(\ref{app.sol-2}) and the maximum principle $|\ua| \le 1$ on $\Omega_{T}$ and $0 \le \etaa \le 1$ a.e. in $\Omega_{T}$, we obtain that there exists a subsequence $\{U_{\varepsilon,\nu,\delta_n}\}_{n}$ and a function ${ U}_{\varepsilon, \nu} \in C({[0,T]}; \mathfrak{X}) \cap W^{1,2}({0, T}; \mathfrak{X}) \cap L^{\infty}({0, T}; \mathfrak{W})$ such that
\begin{equation}\label{U-conv-1}
\left\{ \begin{array}{ll}
U_{\varepsilon,\nu,\delta_n} \to { U}_{\varepsilon, \nu} \ \ & \mbox{ in } C({[0,T]}; \mathfrak{X}), \vspace{3mm}\\
\ds & \mbox{ weakly in } W^{1,2}({0, T}; \mathfrak{X}), \vspace{3mm}\\
\ds & \mbox{ weakly-}\ast \mbox{ in } L^{\infty}({0, T}; \mathfrak{W}), \vspace{3mm}\\
\ds \nu {\bf u}_{\varepsilon, \nu, \delta_n} \to \nu {\bf u}_{\varepsilon, \nu} \ \ & \mbox{ weakly-}\ast \mbox{ in } L^{\infty}({0, T}; W^{1,N+1}(\Omega; \R^{M})), \vspace{3mm}\\
\ds |{\bf u}_{\varepsilon, \nu, \delta_n}| \to 1 & \mbox{ strongly in } L^{2}({0, T}; {H}),
\end{array} \right.
\end{equation}
as $n\to \infty$ by the Aubin type compactness theorem \cite[Corollary 4]{MR0916688}. Then, we also see that
\begin{equation}\label{u=1}
|{\bf u}_{\varepsilon, \nu}| = 1\ \ \mbox{ a.e. in } \Omega_{T}.
\end{equation}

By $|\ua| \le 1$ in $\overline{\Omega_{T}}$, the Lebesgue dominated convergence theorem implies that
\begin{equation}\label{u-conv-1}
{\bf u}_{\varepsilon, \nu, \delta_n} \to {\bf u}_{\varepsilon}\ \ \mbox{ strongly in } L^{r}({0, T};L^{r}(\Omega;\R^{M}))\ \ \mbox{ as } n\to\infty,
\end{equation}
for all $r \in [1, \infty)$.

\medskip
By $|\ua| \le 1$  in $\overline{\Omega_{T}}$ again, we can show that
\begin{equation*}
\int_{0}^{T} \int_{\Omega} |\varpi_{\delta}(\ua)| dx dt \le \frac{1}{\delta} \int_{0}^{T}\int_{\Omega} (1 - |\ua|^{2})^{2} dxdt + \frac{1}{\delta} \int_{0}^{T}\int_\Omega (1 - |\ua|^{2}) |\ua|^{2} dxdt.
\end{equation*}
By (\ref{app.sol-2}), the first term of right-hand side is uniformly bounded with respect to $\delta$. Multiplying both sides of the equation of $\ua$ by $\ua$, integrating by parts and applying H\"older inequality, and \eqref{ene-inq-2}--\eqref{app.sol-6}, we see that
\begin{align*}
&  \frac{1}{\delta} \int_{0}^{T} \int_{\Omega} (1-|\ua|^{2}) |\ua|^{2} dxdt \\[1ex]
\le & \ds \|\partial_{t}\ua\|_{L^{2}({0, T};
{\mathbb{X}})} \|\ua\|_{L^{2}({0, T};
{\mathbb{X}})} \\[1ex]
& \ds + \|\alpha(\etaa) [\nabla f_{\varepsilon}](\nabla \ua)\|_{L^{\infty}({0, T};L^{\infty}(\Omega; \R^{MN}))} \|\nabla \ua\|_{L^{1}({0, T};L^{1}(\Omega; \R^{MN}))} \\[1ex]
    & \ds + \frac{1}{2} \|{\kappa} \nabla \ua\|^{2}_{L^{2}({0, T};L^{2}(\Omega;\R^{MN}))} + {\|\nu \nabla {\bf u}\|_{L^{N+1}({0, T}; L^{N+1}(\Omega; \R^{MN}))}^{N+1}} < \infty.
\end{align*}
Hence, ${\bf f}_{{\varepsilon, \nu, \delta}} := - \varpi_{\delta}(\ua)$ is uniformly bounded with respect to $\delta$ in $L^{1}({0, T}; L^{1}(\Omega; \R^{M}))$. 
Therefore, we can apply Lemma \ref{lem-cpt-1} to get
\begin{equation*}
\nabla {\bf u}_{\varepsilon, \nu, \delta_n} \to \nabla {\bf u}_{\varepsilon, \nu} \ \ \mbox{ strongly in } L^{q}({0, T} ; L^{q}(\Omega; \R^{MN}))\ \ \mbox{ as } n\to\infty, \mbox{ for all } q \in [1,2).
\end{equation*}
The above convergence, (\ref{app.sol-2}) and (\ref{U-conv-1}) imply that
\begin{align}\label{u-conv-2}
    & ~~\alpha(\eta_{\varepsilon, \nu, \delta_n}) [\nabla f_{\varepsilon}](\nabla {\bf u}_{\varepsilon, \nu, \delta_n}) \to \alpha(\eta_{\varepsilon, \nu}) [\nabla f_{\varepsilon}](\nabla {\bf u}_{\varepsilon, \nu}) 
    \nonumber
    \\ 
    & \mbox{ weakly-}\ast \mbox{ in } L^{\infty}({0, T} ; L^{\infty}(\Omega ; \R^{MN})), \mbox{ as  $n\to \infty$.}
\end{align}

\medskip
Next, we use wedge product technique. Taking the wedge product of the equation of $\ua$ with $\ua$,  it follows that
\begin{align*}
\partial_{t} \ua \wedge \ua - \mbox{div}( & (\alpha(\etaa)[\nabla f_{\varepsilon}](\nabla \ua) \nonumber \\[1ex]
& + {\kappa^2}\nabla {\bf u}_{{\varepsilon, \nu, \delta}} {+ \nu |\nu \nabla {\bf u}_{{\varepsilon, \nu, \delta}}|^{N-1} \nu \nabla {\bf u}_{{\varepsilon, \nu, \delta}}}) \wedge \ua) {= 0,}
\end{align*}
by (\ref{distrwedge}). Setting $\delta = \delta_{n}$ and letting $n\to\infty$, we get
\begin{align}\label{eq-u-1}
\int_{0}^{T} \int_{\Omega} \{ & \langle(\partial_{t}{\bf u}_{\varepsilon, \nu} \wedge {\bf u}_{\varepsilon, \nu}), \bomega\rangle_2 + {\sum_{i=1}^{N}} \langle((\alpha(\eta_{\varepsilon, \nu})[\nabla f_{\varepsilon}](
{\partial_{x_{i}}} {\bf u}_{\varepsilon, \nu}) \nonumber\\[1ex]
& + {\kappa^2} 
{\partial_{x_{i}}} {\bf u}_{\varepsilon, \nu}{+ \nu |\nu \nabla {\bf u}_{\varepsilon, \nu}|^{N-1} \nu 
{\partial_{x_{i}}} {\bf u}_{\varepsilon, \nu}}) \wedge {\bf u}_{\varepsilon, \nu}), 
{\partial_{x_{i}}} \bomega\rangle_2 \} dx dt = 0,
\end{align}
for $\bomega \in L^{\infty}({0, T} ; W^{1,N+1}(\Omega, \Lambda_{2}(\R^{M}))
)$ by (\ref{U-conv-1}), (\ref{u-conv-1}) and (\ref{u-conv-2}).

Taking $\bomega = ({\bf u}_{\varepsilon, \nu} \wedge 
{\bm{\psi}})$ in (\ref{eq-u-1}) for $
{\bm{\psi}} \in C^{1}(\overline{\Omega_{T}}{; \R^{M}})$ and using \eqref{hodge_inner} and \eqref{distrwedge}, we have
\begin{align*}
\int_{0}^{T} \int_{\Omega} \{ & (\partial_{t}{\bf u}_{\varepsilon, \nu} \wedge {\bf u}_{\varepsilon, \nu} \wedge *( {\bf u}_{\varepsilon, \nu} \wedge 
{\bm{\psi}})) + {\sum_{i=1}^{N}}(\alpha(\eta_{\varepsilon, \nu})[\nabla f_{\varepsilon}](
{\partial_{x_{i}}} {\bf u}_{\varepsilon, \nu}) \\[1ex]
& + {\kappa^2} 
{\partial_{x_{i}}} {\bf u}_{\varepsilon, \nu}{+ \nu |\nu \nabla {\bf u}_{\varepsilon, \nu}|^{N-1} \nu 
{\partial_{x_{i}}} {\bf u}_{\varepsilon, \nu}}) \wedge 
{\partial_{x_{i}}} ( {\bf u}_{\varepsilon, \nu} \wedge *({\bf u}_{\varepsilon, \nu} \wedge 
{\bm{\psi}})) \} dx dt = 0.
\end{align*}
Noting (\ref{u=1}), we see that
\begin{equation*}\label{note-1}
\begin{array}{l}
\ds \partial_{t}{\bf u}_{\varepsilon, \nu} \cdot {\bf u}_{\varepsilon, \nu} = 0,\ \ \vspace{3mm}\\
\ds (\alpha(\eta_{\varepsilon, \nu})[\nabla f_{\varepsilon}](\nabla {\bf u}_{\varepsilon, \nu}) + {\kappa^2}\nabla {\bf u}_{\varepsilon, \nu}{+ \nu |\nu \nabla {\bf u}_{\varepsilon, \nu}|^{N-1} \nu \nabla {\bf u}_{\varepsilon, \nu}})^{t} {\bf u}_{\varepsilon, \nu} = 0,\ \
\end{array}
\mbox{ a.e. in } \Omega_{T}.
\end{equation*}
According to \eqref{tripleprod}, we obtain
\begin{equation*}
{\bf u}_{\varepsilon, \nu} \wedge *( {\bf u}_{\varepsilon, \nu} \wedge 
{\bm{\psi}}) =({\bf u}_{\varepsilon, \nu} \cdot 
{\bm{\psi}})* {\bf u}_{\varepsilon, \nu} -*
{\bm{\psi}}.
\end{equation*}
Having the previous equality in mind, we can see that
\begin{align}
& \int_{0}^{T} \int_{\Omega} \{ \partial_{t}{\bf u}_{\varepsilon, \nu} \cdot 
{\bm{\psi}} + (\alpha(\eta_{\varepsilon, \nu})[\nabla f_{\varepsilon}](\nabla {\bf u}_{\varepsilon, \nu}) + {\kappa^2}\nabla {\bf u}_{\varepsilon, \nu}{+ \nu |\nu \nabla {\bf u}_{\varepsilon, \nu}|^{N-1} \nu \nabla {\bf u}_{\varepsilon, \nu}}) : \nabla 
{\bm{\psi}} \} dx dt \nonumber \\[1ex]
= & \int_{0}^{T} \int_{\Omega} 
\mu_{\varepsilon, \nu} ({\bf u}_{\varepsilon, \nu} \cdot 
{\bm{\psi}}) dxdt, \label{u-weak-1}
\end{align}
for $\bm{\psi} \in C^{1}(\overline{\Omega_{T}}{; \R^{M}})$ by \eqref{k-basis}, \eqref{inner}, \eqref{hodge_inner}, \eqref{hodgederiv}.


On the other hand, it follows that
\begin{align*}
&\|\alpha'(\eta_{{\varepsilon, \nu, \delta}})[\nabla f_{\varepsilon}](\nabla {\bf u}_{{\varepsilon, \nu, \delta}})\|_{L^{\infty}({0, T} ; L^{\infty}(\Omega{; \R^{MN}}))} \\[1ex]
& \le \frac{\|\alpha'\|_{C([0,1])}}{\alpha^*} \|\alpha(\eta_{{\varepsilon, \nu, \delta}})[\nabla f_{\varepsilon}](\nabla {\bf u}_{{\varepsilon, \nu, \delta}})\|_{L^{\infty}({0, T} ; L^{\infty}(\Omega{; \R^{MN}}))} < \infty.
\end{align*}
Hence, we have
\begin{align}\label{u-conv-3}
    \alpha'( \eta & _{\varepsilon, \nu, \delta_n})[\nabla f_{\varepsilon}](\nabla {\bf u}_{{\varepsilon, \nu, \delta}}) \to \alpha'(\eta_{\varepsilon, \nu})[\nabla f_{\varepsilon}](\nabla {\bf u}_{\varepsilon, \nu}) 
    \nonumber
    \\ 
    & \mbox{ weakly-}\ast \mbox{ in } L^{\infty}({0, T} ; L^{\infty}(\Omega ; \R^{MN})),
\end{align}
as $n\to\infty$. By (\ref{U-conv-1}) and (\ref{u-conv-3}), it follows that
\begin{equation*}
\int_{0}^{T}(\partial_{t}\eta_{\varepsilon, \nu}(t) + g(\eta_{\varepsilon, \nu}(t)) + \alpha'(\eta_{\varepsilon, \nu}(t))f_{\varepsilon}(\nabla {\bf u}_{\varepsilon, \nu}(t)), \varphi)_{H}dt + \int_{0}^{T}(\nabla \eta_{\varepsilon, \nu}(t), \nabla \varphi)_{H}dt = 0,
\end{equation*}
for any $\varphi \in L^{2}({0, T}; V) \cap L^{\infty}(\Omega_{T})$.

{
Finally, we confirm (\ref{ene-inq-2})-(\ref{app.sol-6}). By (\ref{ene-inq-1}) and (\ref{U-conv-1}), the energy inequality (\ref{ene-inq-2}) immediately holds. Moreover, the estimates (\ref{app.sol-4})-(\ref{app.sol-6}) also follow by (\ref{ene-inq-2}), (\ref{u-weak-1}) and (\ref{eq-u-1}), respectively.
}

}

\subsection{The limit as $\nu \to 0$.}

Now, the aim is to solve the limit problem in (P)$_{\varepsilon,\nu}$, as $\nu\to 0$ assuming $0<\varepsilon<\varepsilon_0$; i.e. we solve:
\paragraph{\boldmath Problem (P)$_{\varepsilon}$: }{
    \begin{align*}
        & \begin{cases}
            \partial_t \eta_{\varepsilon} -\mathit{\Delta} \eta_{\varepsilon} +g(\eta_{\varepsilon}) +\alpha'(\eta_{\varepsilon}) f_\varepsilon(\nabla {\bf u}_{\varepsilon}) = 0 \mbox{ in } \Omega_T,
            \\[0.5ex]
            \nabla \eta_{\varepsilon} \cdot {{\bf n}_\Gamma} = 0 \mbox{ on $ \Gamma_T $,}
            \\[0.5ex]
            \eta_{\varepsilon}(0, x) = \eta_0(x), ~ x \in \Omega;
        \end{cases}
        \\[1ex]
        & \begin{cases}
            \partial_t {\bf u}_{\varepsilon} -\mathrm{div} \bigl( \alpha(\eta_{\varepsilon}) {\,[\nabla f_\varepsilon]}(\nabla {\bf u}_{\varepsilon}) +{\kappa^2} \nabla {\bf u}_{\varepsilon} \bigr) 
            = \mu_{\varepsilon} {\bf u}_{\varepsilon} \mbox{ in $ \Omega_T $,}
            \\[0.5ex]
            \bigl( \alpha(\eta_{\varepsilon}) [\nabla f_\varepsilon](\nabla {\bf u}_{\varepsilon}) +{\kappa^2} \nabla {\bf u}_{\varepsilon} \bigr) {{\bf n}_\Gamma} = 0 \mbox{ on $ \Gamma_T $,}
            \\[0.5ex]
            {\bf u}_{\varepsilon}(0, x) = {\bf u}_0(x), ~ x \in \Omega;
        \end{cases}
    \end{align*}
}
together with:
\begin{equation*}
    \mu_{\varepsilon} := (\alpha(\eta_{\varepsilon})[\nabla f_{\varepsilon}](\nabla {\bf u}_{\varepsilon}) + {\kappa^2}\nabla {\bf u}_{\varepsilon} ) : \nabla {\bf u}_{\varepsilon}, ~\mbox{a.e. in $ \Omega_T $.}
\end{equation*}

\begin{thm}[]
Let ${ U}_{0}=[\eta_{0}, {\bf u}_{0}] \in {\mathfrak{W}}
$ with  {${\bf u}_{0} \in \mathbb{S}^{M-1}$} in $\Omega$. Then, there exists ${ U}_{\varepsilon} = [\eta_{\varepsilon}, {\bf u}_{\varepsilon}] \in C({[0,T]};{\mathfrak{X}}
)$ such that ${ U}_{\varepsilon}$ satisfies (P)$_{\varepsilon}$ in the sense of distributions and
\begin{equation}\label{app.sol-7}
\left\{ \begin{array}{l}
\ds { U}_{\varepsilon} \in W^{1,2}({0, T} ; {\mathfrak{X}}
) \cap L^{\infty}({0, T} ; {\mathfrak{W}}
), \vspace{3mm}\\
    {\kappa^2} \nabla {\bf u}_{\varepsilon} \in L^{\infty}({0, T} ; L^{2}(\Omega;\R^{MN})), \vspace{3mm}\\
\ds 0 \le \eta_{\varepsilon} \le 1\ \ \mbox{ a.e.  } \Omega_{T},\ \
    {\bf u}_{\varepsilon} \in {\mathbb{S}^{M-1}}\ \ \mbox{  {a.e. in $ \Omega_T $,}}
\end{array}\right.
\end{equation}
and
\begin{equation*}\label{ene-inq-3}
\mathcal{F}_{\varepsilon}({ U}_{\varepsilon}(s)) + \int_{0}^{s} \|\partial_{t}{ U}_{\varepsilon}(t)\|^{2}_{\mathfrak{X}
    } dt \le \mathcal{F}_{\varepsilon}({ U}_{0})\ \ \mbox{ {for all} } s \in [0,T],
\end{equation*}
\begin{align}
& \bigl\|\alpha(\eta_{\varepsilon})[\nabla f_{\varepsilon}](\nabla {\bf u}_{\varepsilon})\bigr\|_{L^{\infty}(
{\Omega_{T}}; \R^{MN}))} \le C, \label{app.sol-8} \\[1ex]
& \bigl\|\mbox{div}(\alpha(\eta_{\varepsilon})[\nabla f_{\varepsilon}](\nabla {\bf u}_{\varepsilon}) + {\kappa^2}\nabla {\bf u}_{\varepsilon} )\bigr\|_{L^{2}({0, T};L^{1}(\Omega;\R^{M}))} < C, \label{app.sol-9}\\[1ex]
& \bigl\|\mbox{div}((\alpha(\eta_{\varepsilon})[\nabla f_{\varepsilon}](\nabla {\bf u}_{\varepsilon}) + {\kappa^2}\nabla {\bf u}_{\varepsilon}) \wedge {\bf u}_{\varepsilon})\bigr\|_{L^{2}(
{\Omega_{T};\Lambda_2(\R^M)})
)} < C, \label{app.sol-10}
\end{align}
where $C$ does not depend on $\varepsilon$. Moreover, if ${\bf u}_{0} \in \mathbb{S}^{M-1}_{+,r}$ in $\Omega$ for $r \in (0,1)$, then ${\bf u}_{\varepsilon} \in \mathbb{S}^{M-1}_{+,r}$ a.e. in $\Omega_{T}$.
\end{thm}

\paragraph{Proof.}{
By the same argument as that in the proof of Theorem \ref{thm-1}, there exists a subsequence $\{{ U}_{\varepsilon, \nu_n}\}_{n}$ and a function ${ U}_{\varepsilon} \in C({[0,T]}; {\mathfrak{X}}
) \cap W^{1,2}({0, T}; {\mathfrak{X}}
) \cap L^{\infty}({0, T}; {\mathfrak{W}}
)$ such that
\begin{equation}\label{U-conv-4}
\left\{ \begin{array}{ll}
{ U}_{\varepsilon, \nu_n} \to { U}_{\varepsilon} \ \ & \mbox{ in } C({[0,T]}; {\mathfrak{X}}
), \vspace{3mm}\\
\ds & \mbox{ weakly in } W^{1,2}({0, T}; {\mathfrak{X}}
), \vspace{3mm}\\
\ds & \mbox{ weakly-}\ast \mbox{ in } L^{\infty}({0, T}; {\mathfrak{W}}
), \vspace{3mm}\\
\ds {\nu (| \nu \nabla{\bf u}_{\varepsilon, \nu}|^{N-1} \nu \nabla {\bf u}_{\varepsilon, \nu}) \to 0} \ \ & \mbox{ in } {L^{2}({0, T}; L^{2}(\Omega; \R^{MN}))}, \vspace{3mm}\\
\ds |{\bf u}_{\varepsilon, \nu}| \to 1 & \mbox{ strongly in } L^{2}({0, T}; {H}),
\end{array} \right.
\end{equation}
as $\nu \to 0$. Then, we also see that
\begin{equation}\label{u=1-2}
{\bf u}_{\varepsilon} \in {\mathbb{S}^{M-1}}\ \ \mbox{ a.e. in } \Omega_{T}.
\end{equation}
After that, the proof follows the same way as Theorem \ref{thm-1} using Lemma \ref{lem-cpt-2} instead of Lemma \ref{lem-cpt-1}. In fact, we have

\begin{align}\label{eq-u-2}
\int_{0}^{T} \int_{\Omega} \{ {\langle}(\partial_{t}{\bf u}_{\varepsilon} \wedge {\bf u}_{\varepsilon}), \bomega{\rangle_{2}} + {\sum_{i=1}^{N}} \langle((\alpha(\eta_{\varepsilon})[\nabla f_{\varepsilon}](
{\partial_{x_{i}}} {\bf u}_{\varepsilon}) + {\kappa^2}
{\partial_{x_{i}}} {\bf u}_{\varepsilon}) \wedge {\bf u}_{\varepsilon}), 
{\partial_{x_{i}}} \bomega \rangle_2\} dx dt = 0,
\end{align}
for $\bomega \in L^{\infty}({0, T} ; {H^1(\Omega; \Lambda_2(\R^M))})$, 
 and
\begin{align*}\label{eq-u-3}
\int_{0}^{T} \int_{\Omega} \{ \partial_{t}{\bf u}_{\varepsilon} \cdot 
{\bm{\psi}} + (\alpha(\eta_{\varepsilon})[\nabla f_{\varepsilon}](\nabla {\bf u}_{\varepsilon}) + \kappa^2 \nabla {\bf u}_{\varepsilon}) : \nabla 
{\bm{\psi}} \} dx dt = \int_{0}^{T} \int_{\Omega} \mu_{\varepsilon} ({\bf u}_{\varepsilon} \cdot 
{\bm{\psi}}) dxdt,
\end{align*}
for $
{\bm{\psi}} \in C^{1}(\overline{\Omega_{T}}{;\R^{M}})$.

The rest of the proof is the same as that of Theorem \ref{thm-1} and we omit the details.

Finally, let ${\bf u}_{0} \in \mathbb{S}^{M-1}_{+,r}$ for $r \in (0,1)$. %
Therefore, Theorem \ref{max_principle} (maximum principle) implies that
\begin{equation*}\label{u=1-1.5}
{\bf u}_{\varepsilon, \nu} \in \mathbb{S}^{M-1}_{+,r} \ \ \mbox{ a.e. in } \Omega_{T},
\end{equation*}
which, by (\ref{U-conv-4}) implies that ${\bf u}_{\varepsilon} \in \mathbb{S}^{M-1}_{+,r}$ a.e. in $\Omega_{T}$.
}
}

\subsection{The limit as $\varepsilon \to 0$.}\label{subsec:final}

Finally, we solve the initial system (P) by letting $\varepsilon\to 0$. As $ \varepsilon \to 0 $, the limit problem is formulated as follows.

\paragraph{\boldmath Problem (P) :}{
    \begin{align*}
        & \begin{cases}
            \partial_t \eta -\mathit{\Delta} \eta +g(\eta) +\alpha'(\eta) |\nabla {\bf u}| = 0 \mbox{ in } \Omega_T,
            \\[0.5ex]
            \nabla \eta \cdot {{\bf n}_\Gamma} = 0 \mbox{ on $ \Gamma_T $,}
            \\[0.5ex]
            \eta(0, x) = \eta_0(x), ~ x \in \Omega;
        \end{cases}
        \\[1ex]
        & \begin{cases}
            \partial_t {\bf u} -\mathrm{div} \bigl( 
            \alpha(\eta)\mathcal{B} +{\kappa^2} \nabla {\bf u} \bigr) 
            = \mu {\bf u} \mbox{ in $ \Omega_T $,}
            \\[0.5ex]
            \bigl( 
            \alpha(\eta)\mathcal{B} +{\kappa^2} \nabla {\bf u} \bigr) {{\bf n}_\Gamma} = 0 \mbox{ on $ \Gamma_T $,}
            \\[0.5ex]
            {\bf u}(0, x) = {\bf u}_0(x), ~ x \in \Omega;
        \end{cases}
    \end{align*}
}
together with:
\begin{equation*}
    \mathcal{B} \in \mathrm{Sgn}^{M, N}(\nabla {\bf u}), ~\mbox{and}~ \mu := (
    \alpha(\eta)\mathcal{B} + {\kappa^2}\nabla {\bf u}) : \nabla {\bf u}, ~\mbox{a.e. in $ \Omega_T $.}
\end{equation*}

Now, to prove the Main Theorem, it is sufficient to verify the following Theorem.

\begin{thm}\label{thm-5}
Let ${ U}_{0}=[\eta_{0}, {\bf u}_{0}] \in {\mathfrak{W}}
$ with 
    ${\bf u}_{0} \in {\mathbb{S}^{M-1}}$ in $\Omega$. Then,
        there exist ${ U} = [\eta, {\bf u}] \in L^2(0,T;{\mathfrak{X}}) $  and $ \mathcal{B} \in L^\infty(\Omega_T; \R^{MN}) $ such that
\begin{equation*}
\left\{ \begin{array}{l}
\ds { U} \in W^{1,2}({0, T} ; {\mathfrak{X}}
) \cap L^{\infty}({0, T} ; {\mathfrak{W}}
), \vspace{3mm}\\
    {\kappa^2} \nabla {\bf u} \in L^{\infty}(0,T; L^{2}(\Omega;\R^{MN})), \vspace{3mm}\\
\ds 0 \le \eta \le 1\ \ \mbox{ a.e. in } \Omega_{T},\ \ {\bf u} \in {\mathbb{S}^{M-1}}\ \ \mbox{  a.e. in } \Omega_{T},
\end{array}\right.
\end{equation*}
and
\begin{equation*}\label{ene-inq-4}
\mathcal{F}({ U}(s)) + \int_{0}^{s} \|\partial_{t}{ U}(t)\|^{2}_{\mathfrak{X} 
    } dt \le \mathcal{F}({ U}_{0})\ \ \mbox{ {for all} } s \in [0,T].
\end{equation*}

    \noindent
    Also, there exists a constant $ C > 0 $, independent of $ \kappa $, such that:
\begin{align*}
    & \bigl\|\alpha(\eta)\mathcal{B} \bigr\|_{L^{\infty}(\Omega_T; \R^{MN})} \le C, 
\\[1ex]
        & \bigl\| \mbox{div}\bigl(\alpha(\eta) \mathcal{B} + {\kappa^2}\nabla {\bf u}\bigr) \bigr\|_{L^{2}(0,T; L^{1}(\Omega;\R^{M}))} < C, 
\\[1ex]
        & \bigl\|\mbox{div}\bigl( \alpha(\eta) \mathcal{B} + {\kappa^2}\nabla {\bf u}) \wedge {\bf u} \bigr) \bigr\|_{L^{2}(
{\Omega_{T};\Lambda_{2}(\R^{M})})} < C. 
\end{align*}
    Moreover, if ${\bf u}_{0} \in \mathbb{S}^{M-1}_{+,r}$ in $\Omega$ for  $r \in (0, 1)$, then 
    ${\bf u} \in \mathbb{S}^{M-1}_{+,r}$ a.e. in $\Omega_{T}$, where $ \mathbb{S}^{M -1}_{+, r}$ is the subset of $ \mathbb{S}^{M -1} $ given in \eqref{constraint_S}.
\end{thm}
\paragraph{Proof.}{
Let ${ U}_{\varepsilon}$ be the solution to (P)$_{\varepsilon}$ constructed in the previous Section. Due to (\ref{app.sol-7}), (\ref{app.sol-8})-(\ref{app.sol-10}), (\ref{u=1-2}), there exists a subsequence $\{ { U}_{\varepsilon_n}\}_{n}$ and a function ${ U} \in C({[0,T]}; {\mathfrak{X}} 
) \cap W^{1,2}({0, T}; {\mathfrak{X}} 
) \cap L^{\infty}({0, T}; {\mathfrak{W}} 
)$ such that
\begin{equation}\label{U-conv-5}
\left\{ \begin{array}{ll}
{ U}_{\varepsilon_n} \to { U} \ \ & \mbox{ in } C({[0,T]}; {\mathfrak{X}}
), \vspace{3mm}\\
\ds & \mbox{ weakly in } W^{1,2}({0, T}; {\mathfrak{X}} 
), \vspace{3mm}\\
\ds & \mbox{ weakly-}\ast \mbox{ in } L^{\infty}({0, T}; {\mathfrak{W}} 
), \vspace{3mm}\\
\ds |{\bf u}_{\varepsilon_n}| \to 1 & \mbox{ strongly in } L^{2}({0, T}; {H}),
\end{array} \right.
\end{equation}
as $n\to\infty$ by the Aubin type compactness theorem \cite[Corollary 4]{MR0916688}. Then, we also see that
\begin{equation*}\label{u=1-3}
{\bf u} \in {\mathbb{S}^{M-1}} \ \ \mbox{ a.e. in } \Omega_{T}.
\end{equation*}

By $|{\bf u}| \le 1$ in $\overline{\Omega_{T}}$, the Lebesgue dominated convergence theorem implies that
\begin{equation*}\label{u-conv-4}
{\bf u}_{\varepsilon_n} \to {\bf u}\ \ \mbox{ strongly in } L^{r}({0, T};L^{r}(\Omega;\R^{M}))\ \ \mbox{ as } n\to\infty,
\end{equation*}
for all $r \in [1, \infty)$.

Next, we set ${\bf f}_{\varepsilon} := \mu_{\varepsilon} {\bf u}_{\varepsilon}$. Then, ${\bf f}_{\varepsilon}$ is uniformly bounded in $L^{1}({0, T}; L^{1}(\Omega; \R^{M}))$ for $\varepsilon$. By Lemma \ref{lem-cpt-3}, it follows that
\begin{equation*}
\nabla {\bf u}_{\varepsilon_n} \to \nabla {\bf u}\ \ \mbox{ strongly in } L^{1}({0, T}; L^{1}(\Omega; \R^{MN}))\ \ \mbox{ as } n\to\infty.
\end{equation*}
The above convergence and (\ref{app.sol-8}) imply that there exists
    $ \mathcal{B} \in L^\infty(\Omega_T; \R^{MN}) $ such that
\begin{align*}
    & \alpha(\eta_{\varepsilon_n})[\nabla f_{\varepsilon_n}](\nabla {\bf u}_{\varepsilon_n}) \to \alpha(\eta) \mathcal{B}\ \ \mbox{ weakly-}\ast\ \mbox{ in }
        L^\infty(\Omega_T; \R^{MN}),
\\[1ex]
    & \alpha(\eta_{\varepsilon_n})[\nabla f_{\varepsilon_n}](\nabla {\bf u}_{\varepsilon_n})  + {\kappa^2} \nabla {\bf u}_{\varepsilon_n} \to \alpha(\eta)\mathcal{B} + {\kappa^2} \nabla {\bf u} \ \ \mbox{ weakly in }
        L^2(\Omega_T; \R^{MN}),
\end{align*}
as $n\to\infty$. In addition, we see that
\begin{equation*}
\alpha(\eta)\mathcal{B} : \nabla {\bf u} = \alpha(\eta) |\nabla {\bf u}| \ \ \mbox{ a.e. in } \Omega_{T}.
\end{equation*}
Letting $n\to\infty$ in (\ref{eq-u-2}), we get
\begin{equation}\label{eq-u-4}
\int_{0}^{T} \int_{\Omega} \{ {\langle}(\partial_{t}{\bf u} \wedge {\bf u}), \bomega\rangle_2 + {\sum_{i=1}^{N}} \langle((\alpha(\eta){\mathcal{B}_{i}} + {\kappa^2} 
{\partial_{x_{i}}} {\bf u}) \wedge {\bf u}) , 
{\partial_{x_{i}}} \bomega\rangle_2 \} dx dt = 0,
\end{equation}
for {$\bomega \in C^{1}(\overline{\Omega_{T}};\Lambda_2(\R^M))$}.
{Here, $\mathcal{B} = (\mathcal{B}_{i})$, $\mathcal{B}_{i} \in L^{\infty}(\Omega_{T}; \R^{M})$.}
Since $\mathcal{B} \in 
 L^\infty(\Omega_T; \R^{MN}) $, ${\kappa^2} \nabla {\bf u} \in L^{\infty}({0, T}; L^{2}(\Omega; \R^{MN}))$, the density argument implies that (\ref{eq-u-4}) holds for any $\bomega \in L^{\infty}({0, T} ; W^{1,2}(\Omega; {\Lambda_2(\R^M)}))$. Taking $\bomega = ({\bf u} \wedge 
{\bm{\psi}})$ for $
{\bm{\psi}} \in C^{1}(\overline{\Omega_{T}}{; \R^{M}})$ in (\ref{eq-u-4}), we have
\begin{equation*}
\int_{0}^{T} \int_{\Omega} \{ (\partial_{t}{\bf u} \wedge {\bf u}) \wedge {\ast}({\bf u} \wedge 
{\bm{\psi}}) + {\sum_{i=1}^{N}} ((\alpha(\eta){\mathcal{B}_{i}} + {\kappa^2} 
{\partial_{x_{i}}} {\bf u} ) \wedge {\bf u})\wedge 
{\partial_{x_{i}}} ({\ast}({\bf u} \wedge 
{\bm{\psi}})) \} dx dt = 0.
\end{equation*}
In a similar way to Section 5.2,  $\mu=(\alpha(\eta)\mathcal{B} + {\kappa^2} \nabla {\bf u}) : \nabla {\bf u} \in L^{1}({0, T}; {L^{1}(\Omega)})$ and the following equation is satisfied:
\begin{equation*}
\int_{0}^{T} \int_{\Omega} \{ \partial_{t}{\bf u} \cdot 
{\bm{\psi}} + (\alpha(\eta)\mathcal{B} + {\kappa^2}\nabla {\bf u}) : \nabla 
{\bm{\psi}} \} dx dt = \int_{0}^{T} \int_{\Omega} \mu {{\bf u}} \cdot 
{\bm{\psi}} dxdt.
\end{equation*}

The rest of the proof is the same as that of Theorem \ref{thm-1} and we omit it. Finally, in the case that ${\bf u}_{0} \in \mathbb{S}^{M-1}_{+,r}$ a.e. in $\Omega_{T}$, by Theorem \ref{thm-5} and \eqref{U-conv-5} we easily get the conclusion.
}


\section*{Appendix}

In this appendix we recall the notion of Mosco--convergence and some results about it that we use in the paper.

\begin{defn}[Mosco-convergence: cf. \cite{MR0298508}]\label{Def.Mosco}
    Let $ X $ be an abstract Hilbert space. Let $ \Psi : X \rightarrow (-\infty, \infty] $ be a proper, l.s.c., and convex function, and let $ \{ \Psi_n \}_{n = 1}^\infty $ be a sequence of proper, l.s.c., and convex functions $ \Psi_n : X \rightarrow (-\infty, \infty] $, $ n = 1, 2, 3, \dots $.  Then, it is said that $ \Psi_n \to \Psi $ on $ X $, in the sense of Mosco, as $ n \to \infty $, iff. the following two conditions are fulfilled:
\begin{description}
    \item[(\hypertarget{M_lb}{M1}) The condition of lower-bound:]$ \ds \varliminf_{n \to \infty} \Psi_n(\check{w}_n) \geq \Psi(\check{w}) $, if $ \check{w} \in X $, $ \{ \check{w}_n  \}_{n = 1}^\infty \subset X $, and $ \check{w}_n \to \check{w} $ weakly in $ X $, as $ n \to \infty $.
    \item[(\hypertarget{M_opt}{M2}) The condition of optimality:]for any $ \hat{w} \in D(\Psi) $, there exists a sequence \linebreak $ \{ \hat{w}_n \}_{n = 1}^\infty  \subset X $ such that $ \hat{w}_n \to \hat{w} $ in $ X $ and $ \Psi_n(\hat{w}_n) \to \Psi(\hat{w}) $, as $ n \to \infty $.
\end{description}
    As well as, if the sequence of convex functions $ \{ \widehat{\Psi}_\varepsilon \}_{\varepsilon \in \Xi} $ is labeled by a continuous argument $\varepsilon \in \Xi$ with a range $\Xi \subset \mathbb{R}$ , then for any $\varepsilon_{0} \in \Xi$, the Mosco-convergence of $\{ \widehat{\Psi}_\varepsilon \}_{\varepsilon \in \Xi}$, as $\varepsilon \to \varepsilon_{0}$, is defined by those of subsequences $ \{ \widehat{\Psi}_{\varepsilon_n} \}_{n = 1}^\infty $, for all sequences $\{ \varepsilon_n \}_{n=1}^{\infty} \subset \Xi$, satisfying $\varepsilon_{n} \to \varepsilon_{0}$ as $n \to \infty$.
\end{defn}

\begin{rem}\label{Rem.MG}
    Let $ X $, $ \Psi $, and $ \{ \Psi_n \}_{n = 1}^\infty $ be as in Definition~\ref{Def.Mosco}. Then, the following hold.
\begin{description}
    \item[(\hypertarget{Fact5}{Fact\,1})](cf. \cite[Theorem 3.66]{MR0773850} and \cite[Chapter 2]{Kenmochi81}) Let us assume that
    \begin{equation}\label{Mosco01}
    \Psi_n \to \Psi \mbox{ on $ X $, in the sense of  Mosco, as $ n \to \infty $,}
    \vspace{-1ex}
\end{equation}
and
\begin{equation*}
\left\{ ~ \parbox{10cm}{
$ [w, w^*] \in X \times X $, ~ $ [w_n, w_n^*] \in \partial \Psi_n $ in $ X \times X $, $ n \in \N $,
\\[1ex]
$ w_n \to w $ in $ X $ and $ w_n^* \to w^* $ weakly in $ X $, as $ n \to \infty $.
} \right.
\end{equation*}
Then, it holds that:
\begin{equation*}
[w, w^*] \in \partial \Psi \mbox{ in $ X \times X $, and } \Psi_n(w_n) \to \Psi(w) \mbox{, as $ n \to \infty $.}
\end{equation*}
    \item[(\hypertarget{Fact6}{Fact\,2})](cf. \cite[Lemma 4.1]{MR3661429} and \cite[Appendix]{MR2096945}) Let $ d \in \mathbb{N} $, and let $  S \subset \R^d $ be a bounded open set. Then, under the Mosco-convergence as in \eqref{Mosco01}, a sequence $ \{ \widehat{\Psi}_n^S \}_{n = 1}^\infty $ of proper, l.s.c., and convex functions on $ L^2(S; X) $, defined as:
        \begin{equation*}
            w \in L^2(S; X) \mapsto \widehat{\Psi}_n^S(w) := \left\{ \begin{array}{ll}
                    \multicolumn{2}{l}{\ds \int_S \Psi_n(w(t)) \, dt,}
                    \\[1ex]
                    & \mbox{ if $ \Psi_n(w) \in L^1(S) $,}
                    \\[2.5ex]
                    \infty, & \mbox{ otherwise,}
                \end{array} \right. \mbox{for $ n = 1, 2, 3, \dots $;}
        \end{equation*}
        converges to a proper, l.s.c., and convex function $ \widehat{\Psi}^S $ on $ L^2(S; X) $, defined as:
        \begin{equation*}
            z \in L^2(S; X) \mapsto \widehat{\Psi}^S(z) := \left\{ \begin{array}{ll}
                    \multicolumn{2}{l}{\ds \int_S \Psi(z(t)) \, dt, \mbox{ if $ \Psi(z) \in L^1(S) $,}}
                    \\[2ex]
                    \infty, & \mbox{ otherwise;}
                \end{array} \right.
        \end{equation*}
        on $ L^2(S; X) $, in the sense of Mosco, as $ n \to \infty $.
\end{description}
\end{rem}
\begin{ex}[Example of Mosco-convergence]\label{Rem.ExMG01}
    Let $ \{ f_\varepsilon \}_{\varepsilon \geq 0} \subset C(\R^d) $ be the sequence of non-expansive convex functions, as in (A4). Then, for any $ \varepsilon_0 \geq 0 $, the uniform estimate:
        \begin{align*}
            |f_\varepsilon(W) -f_{\tilde{\varepsilon}}(\tilde{W})| \leq & \bigl\| [\varepsilon, W] -[\tilde{\varepsilon}, \tilde{W}] \bigr\|_{\R^{1 +MN}} \leq |\varepsilon -\tilde{\varepsilon}| +\bigr\|W -\tilde{W} \bigr\|_{\R^{MN}},
            \nonumber
            \\
            & \mbox{for all $ \varepsilon, \tilde{\varepsilon} \geq 0  $ and $ y, \tilde{y} \in \R^d $,}
    \end{align*}
    immediately implies:
    \begin{equation*}
        f_\varepsilon \to f_{\varepsilon_0} \mbox{ uniformly on $ \mathbb{R}^{MN} $, as $ \varepsilon \to \varepsilon_0 $,}
    \end{equation*}
    and this uniform convergence leads to:
    \begin{equation*}
        f_\varepsilon \to f_{\varepsilon_0} \mbox{ on $ \mathbb{R}^{MN} $, in the sense of Mosco, as $ \varepsilon \to \varepsilon_0 $.}
    \end{equation*}
\end{ex}
\begin{rem}\label{Mosco00}
    Generally, for a sequence of proper l.s.c. and convex functions on a Hilbert space, the uniform convergence implies the convergence in the sense of Mosco.
\end{rem}

\section*{Appendix B}
In this Appendix, we give some insight and results related to the group of rotations $SO(3)$ as a Riemannian manifold that we use throughout the paper.

First of all,
according to Euler's theorem, any rotation can be characterized by an angle of rotation $w\in [0, 2\pi[$ and an axis of rotation ${\bf n}\in \mathbb{S}^2$. Moreover, the range of the angle can be taken into $[0,\pi]$ (note that $(w,{\bf n})$ and $(-w,{-\bf n})$ correspond to the same rotation) and in case $w\neq \pi$ then, this representation is unique. It is well known that $SO(3)$ is diffeomorphic to the real projective space $\mathbb{P}^3(\R)$. However, we will use instead that $SO(3)$ is diffeomorphic to a quotient group in the unit hypersphere $\mathbb{S}^3$ in $\R^4$.

We construct the following surjective homeomorphism between the group of unit quaternions (which is isomorphic to the unit special group $SU(2)$) and $SO(3)$:
$${\bf q}=(q^0,q^1,q^2,q^3)\mapsto {\bf q}(\bullet){\bf q}^{-1}.$$ It is easy to see that this is a two--to--one homeomorphism and that for any rotation there are two unit quaternions associated to it: ${\bf q}$ and $-{\bf q}$. Moreover, the following relation between quaternions and axis-angle parametrization is obtained:
$$q^0=\pm|\cos(\frac{w}{2})|\,,\quad q^i=\pm |\sin(\frac{w}{2})|n_i\,, i=1,\ldots,3.$$
Therefore, since the group of unit quaternions is isomorphic to $\mathbb{S}^3$, we conclude that $SO(3)$ is diffeomorphic to the quotient group $\mathbb{S}^3/\sim$  where $\sim$ is the equivalence class of antipodal points in the hypersphere; i.e. ${\bf p}\sim \bf{q}$ iff ${\bf p}+\bf{q}=0$. Furthermore, we see that points in the equator (i.e. of the form $(0,p_1,p_2,p_3)$) correspond to the angle of rotation $w=\pi$. Then, it is obtained that rotations with angle of rotation $w\in [0,\pi[$ can be uniquely identified as points in the open upper hemisphere $\mathbb{S}^3_+$; i.e. $(p_0,p_1,p_2,p_3)$ with $0<p_0\leq 1$. We also denote by ${\bf p_0}$ the north pole on the sphere.

\medskip

In the paper, we will restrict ourselves to initial rotations ${\bf P_0}$ such that the angle of rotation $w_{{\bf P}_0}$ satisfies $0\leq w_{{\bf P}_0}<\pi$. Therefore, we can equivalently work with initial data ${\bf {\bf u}}_0\in \mathbb{S}^3_+$ and with the geometry of $\mathbb{S}^3_+$, which is much simpler and it has been used much wider in variational constraint problems than that of $SO(3)$. However, in order that the solution to the flow still belongs to $\mathbb{S}^3_+$ and therefore it can be identified with a rotation, we need the following result:

\begin{thm}(Maximum principle)\label{max_principle}
  Suppose that ${\bf {\bf u}}_0\in \overline{B_g({\bf p}_0; R)}$, with $R<\frac{\pi}{2}$ (equivalently $w_{{\bf P}_0}\in [0,r]$ with $r<\pi$). Then, the solution to (P)$_{\varepsilon, \nu}$ satisfies
  $$ {\bf u}_{\varepsilon,\nu}\in \overline{B_g({\bf p}_0; R)}\,,\quad {\rm a.e. \ in \ } \Omega.$$
\end{thm}
\begin{proof}
    The proof follows the proof of \cite[Lemma 4]{MR3941105}. We proceed by contradiction. Let $T^* = \inf\{t\in [0, T[ : {\bf u}(t,\Omega)\subset \hspace{-2.25ex} /~ \overline{B_g({\bf p}_0;R)}\}$. Due to
continuity of ${\bf {\bf u}}$, there is a $\delta > 0$ such that ${\bf {\bf u}}(t;\Omega)\subset B_g({\bf p}_0;\frac{\pi}{2})$ for $t\in [0;T^*+\delta[$.

    We now take the equation for ${\bf u}_{\varepsilon,\nu}$ in (P)$_{\varepsilon, \nu}$ and we take the projection $\pi_{{\bf u}_{\varepsilon,\nu}}$ from $\mathbb{S}^{M-1}$ to $T_{{\bf u}_{\varepsilon,\nu}}\mathbb{S}^{M-1}$. Noting that $\pi_{{\bf u}_{\varepsilon,\nu}}(\mu_{\varepsilon,\nu} {\bf u}_{\varepsilon,\nu}) = 0$, we get
\begin{equation}\label{eq.proj.}
\partial_t{\bf u}=\pi_{\u}{\rm div}(Z),
\end{equation}
 with
$$\u=\u_{\varepsilon,\nu}\,, \quad Z:= \alpha(\eta_{\varepsilon, \nu}) {\,[\nabla f_\varepsilon]}(\nabla {\bf u}_{\varepsilon, \nu}) +{\kappa^2} \nabla {\bf u}_{\varepsilon, \nu}+\nu |\nu \nabla {\bf u}_{\varepsilon, \nu}|^{N-1} \nu \nabla {\bf u}_{\varepsilon, \nu}.$$

     We choose on $B_g({\bf p}_0;\frac{\pi}{2})$ a polar coordinate system ${\bf p}\mapsto(p^r; p^{\theta_1},\ldots, p^{\theta_{M-2}})$ centered at ${\bf p}_0$.
     Next, we compute the second equation in (P)$_{\varepsilon, \nu}$ for the radial coordinate:


The metric in the polar coordinates around the north pole is the following one: 
        $${g=\begin{pmatrix}
                                                                                                                                1 & 0 & 0  & \cdots & 0 \\
                                                                                                                                0 & \sin^2(r) & 0  & \cdots & 0 \\
                                                                                                                                0 & 0 & \sin^2(r)\sin^2(\theta_1)  & \cdots &0
                                                                                                                                \\ \vdots & \vdots & \vdots & \ddots & \vdots   \\
                                                                                                                                0 & 0 & 0   & \cdots & \sin^2(r)\sin^2(\theta_1)\cdot\ldots\cdot\sin^2(\theta_{M-3})
        \end{pmatrix}.}$$
Therefore, the Christoffel symbols of the second kind for the variable $r$ are 
        $${\Gamma^r=-\frac{\sin(2r)}{2}\begin{pmatrix}
                                                                                            0 & 0 & 0 & \ldots & 0 \\
                                                                                            0 & 1 & 0 & \ldots & 0 \\
                                                                                            0 & 0 & \sin^2(\theta_1) & \ldots & 0 \\
                                                                                            \vdots & \vdots & \vdots & \ddots & 0 \\
                                                                                            0 & 0 & 0 & \ldots & \sin^2(\theta_1)\cdot\ldots\cdot\sin^2(\theta_{M-3})
        \end{pmatrix}.}$$
Next, we note that   (see \cite{MR164306})\begin{equation*}\label{eells-sampson}\pi_{\bf {\bf u}}({\rm div} {Z})^i={\rm div} {\bf Z}^i+\sum_{j,k,l}\Gamma^i_{j,k}({\bf {\bf u}}){\bf u}^j_{x^\ell}{\bf Z}^k_l\,, \ i=r,\theta_1,\ldots,\theta_{M-3}\,,\quad {\rm for \ any \ }  Z\in W^{1,1}(\Omega;{\mathbb S}^{N-1}).\end{equation*}

Thus, we get that the equation for the radial coordinate in \eqref{eq.proj.} is the following one: $${u}^r_t={\rm div} {\bf Z}^r-\frac{\sin{2{u}^r}}{2}\left(\frac{1}{\sqrt{\varepsilon^2+|\nabla {\u}|^2}}+\kappa^2+\nu^{N+1}|\nabla {\u}|^{N-1}\right)\times $$$$\times \left(|\nabla u^{\theta_1}|^2+\sum_{i=1}^{M-3}\sin^2({u}^{\theta_1})\cdot\ldots\cdot\sin^2(u^{\theta_i})|\nabla u^{\theta_{i+1}}|^2\right).$$

Therefore, since ${u}^r\in [0,\frac{\pi}{2}] $, $${u}^r_t\leq {\rm div} {\bf Z}^r .$$
Thus,
$$\frac{1}{2}\frac{d}{dt}\int_\Omega ({u}^r-R)_+^2=\int_\Omega {u}^r_t({u}^r-R)_+\leq \int_{\Omega}{\rm div} {\bf Z}^r ({u}^r-R)_+=-\int_{\Omega\cap\{{ u}^r>R\}}{\bf Z}^r: \nabla {u}^r$$$$=-\int_{\Omega\cap\{{ u}^r>R\}}|\nabla u^r|^2\left(\frac{1}{\sqrt{\varepsilon^2+|\nabla {\bf u}|^2}}+\kappa^2+\nu^{N+1}|\nabla {\bf u}|^{N-1}\right)\leq 0,$$ and this finishes the proof. 

\end{proof}


\end{document}